\documentclass[11pt]{amsart}
\usepackage[top=3cm,bottom=2cm,left=3cm,right=3cm,marginparwidth=1.75cm]{geometry}
\usepackage[format=plain, font=footnotesize]{caption}
\usepackage{amsmath, amsthm, nccmath, amssymb, mathtools, hyperref}
\hypersetup{
	colorlinks=true,       % false: boxed links; true: colored links
	linkcolor=orange,        % color of internal links
	citecolor=purple,        % color of links to bibliography
	filecolor=magenta,     % color of file links
	urlcolor=blue
}
\usepackage[shortlabels]{enumitem}
\usepackage{graphicx}
\usepackage[dvipsnames]{xcolor}
\usepackage[frozencache,cachedir=minted-cache]{minted}
 
 \usepackage{listings}
\usepackage{url}
\usepackage{cleveref}
\usepackage{systeme}
\usepackage{array}
\usepackage{bigints}
\usepackage{comment}
\usepackage{tikz-cd}
\usepackage{subcaption} % for sub figures! 

\newtheorem{theorem}{Theorem}
\numberwithin{theorem}{section}
\newtheorem{proposition}[theorem]{Proposition}

\newtheorem{corollary}[theorem]{Corollary}
\newtheorem{conjecture}[theorem]{Conjecture}

\theoremstyle{definition}
\newtheorem{definition}[theorem]{Definition}
\newtheorem{remark}[theorem]{Remark}
\newtheorem{example}[theorem]{Example}

\usepackage{listings}
\lstdefinelanguage{code}{
basicstyle=\small\ttfamily,
alsoletter=",
classoffset=1,
keywords={gb, eliminate, saturate, diff, degree, flatten, apply, tensor, product},
keywordstyle={\color{teal}},
classoffset=2,
morekeywords={from, to, list, terms, toList, entries, for, end, if, return},
keywordstyle=\color{blue!70!green!80!white},
classoffset=3,
morekeywords={QQ},
keywordstyle={\color{frenchplum}},
classoffset=4,
morekeywords={ideal, matrix, gens},
keywordstyle={\color{teal}},
xleftmargin=1.5cm,
xrightmargin=1em,
columns=fullflexible,
keepspaces=true,
stepnumber=1,
numbers=none,
captionpos=b,
showspaces=false,
frame=none
}

\definecolor{weborange}{RGB}{255,165,0}
\definecolor{darkgray}{rgb}{.4,.4,.4}
\lstdefinelanguage{cocoa}
{
  basicstyle=\small\ttfamily,
  commentstyle=\color{red!90!black},       % comment style
  stringstyle=\color{green!70!black},      % string literal style
  morecomment=[l]{//},
  morecomment=[l]{--},
  morecomment=[s]{/*}{*/},
  morestring=[b]{"},
  classoffset=1,
  morekeywords={
    define,enddefine,if,endif,for,endfor,
    use,in,then,else,elif,return,
    and,or,
    break,continue,ciao,do,exit,
    ImportByValue,ImportByRef,importbyvalue,importbyref,
    in,isin,IsIn,
    on,opt,PrintLn,println,print,
    quit,ref,return,step,
    toplevel,TopLevel,
    then,to,
    use
  },
  keywordstyle=\color{blue!70!green!80!white},
  classoffset=2,
  morekeywords={
    Ideal,Mat,       Not,Record,Error,
    ideal,mat,matrix,not,record,error,submodule
  },
  keywordstyle=\color{purple!50!white},
  classoffset=3,
  morekeywords={
    TRUE,FALSE,True,False,true,false,
    Lex,Xel,DegLex,DegRevLex,
    PosTo,ToPos,Null
  },
  keywordstyle=\color{brown},
  xleftmargin=1cm,
  xrightmargin=1em,
  columns=fullflexible,
  keepspaces=true,
  stepnumber=1,
  numbers=none,
  captionpos=b,
  showspaces=false,
  showstringspaces=false,
  frame=none
}
\usepackage{mathtools}

\newcommand{\rmi}{{I}}
\newcommand{\rmidist}{\mathcal{I}}
\newcommand{\er}{Erd\"os-R\'enyi} 
\newcommand{\ideal}[1]{\left<{#1}\right>}

\DeclareMathOperator{\lcm}{lcm}
\DeclareMathOperator{\ann}{ann}
\DeclareMathOperator{\ass}{Ass}
\DeclareMathOperator{\pd}{pd}

\DeclareMathOperator{\depth}{depth}

\newcommand{\KK}{\mathbb{K}}
\newcommand{\thelcm}{{\mathbf m}}

\title{Computation of dominant ideals}
\author[WICA II Probabilistic Models Group]{
Anna Maria Bigatti,
Nursel Erey, 
Selvi Kara, 
Augustine O'Keefe, 
Sonja Petrovi\'c,
Pierpaola Santarsiero,
and
Janet Striuli}\thanks{
This project began during the Women in Commutative Algebra (WICA) II meeting at the Centro Internazionale per la Ricerca Matematica in Trento, Italy in October 2023. WICA II was generously supported by NSF grant DMS–2324929.
SK is supported by  NSF grant DMS-2418805.  
SP is partially supported by the Simons Foundation through the Travel Support for Mathematicians Gift \# 854770 and DOE Grant \#1010629. 
PS acknowledges partial support by the Deutsche Forschungsgemeinschaft (DFG, German Research
Foundation)–Projektnummer 445466444 and by the European Union under
NextGenerationEU. PRIN 2022, Prot. 2022E2Z4AK}

\date{}

\begin{document}

\begin{abstract}
We consider the problem of determining whether a monomial ideal is dominant. This property is critical for determining for which monomial ideals the Taylor resolution is minimal. We first analyze dominant ideals with a fixed least common multiple of generators using combinatorial methods. Then, we adopt a probabilistic approach via the \er\ type model, examining both homogeneous and non-homogeneous cases. This model offers an efficient alternative to exhaustive enumeration, allowing the study of dominance through small random samples, even in high-dimensional settings.
\end{abstract}
\maketitle

\section{Introduction}

Interactions between combinatorics and commutative algebra have a rich and long history. A more recent history uses fundamental ideas from probabilistic combinatorics to compute and study average behaviors of objects in commutative algebra. 
This is the setting of our study. In particular, the algebraic object we are interested in is a minimal free resolution of monomial ideals; the combinatorial property we focus on is when the ideal is  \emph{dominant} (Definition~\ref{def:dominant}); and the  computational frameworks we consider focus on two types of questions: how to enumerate and understand all dominant ideals under some constraints (Section~\ref{sec:lcm}), and how to probabilistically sample from the set using models from the literature (Section~\ref{sec:random sampling}).   
Let us explain how these notions are related. 

Monomial resolutions---the free resolutions of monomial ideals---can be quite complex to construct despite dealing with the ``simplest" class of ideals. The most effective approach to constructing monomial resolutions typically relies on combinatorial objects such as simplicial complexes or CW-complexes. These objects \emph{support} the resolution, allowing its structure to be inferred from combinatorial data \cite{EK90, Mer10, BPS98, BS98, MSY00, Nov00, NPS02}. While there are cases where no appropriate CW-complex can support a resolution \cite{Vel08}, this combinatorial approach remains one of the most powerful tools available. A fundamental goal in this area is to develop free resolutions that \emph{preserve minimality}, as maintaining minimality is essential for extracting algebraic invariants from these ideals.

\begin{figure}[h]
    \centering
    \includegraphics[width=0.6\textwidth]{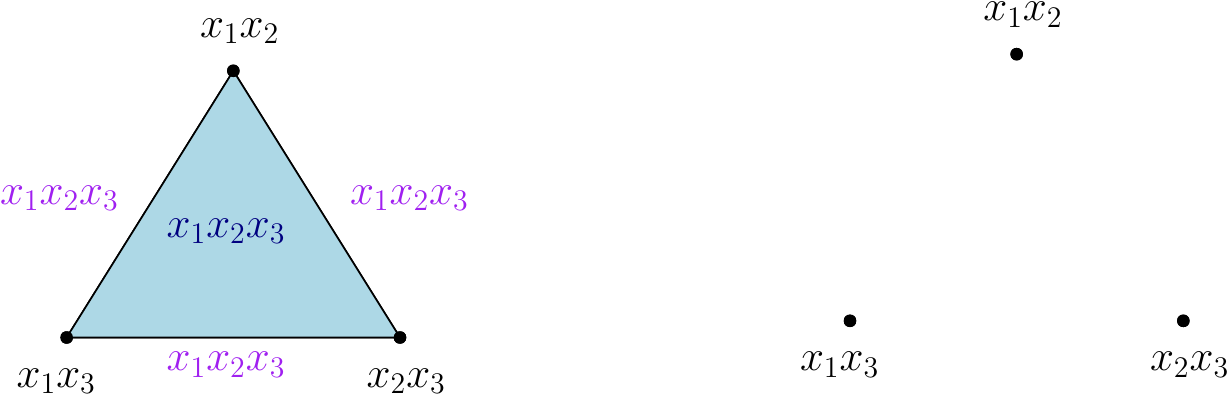}
     \captionsetup{singlelinecheck=off}
    \caption[.]{
Taylor complex (left) and Scarf complex (right) of $I=(x_1x_2, x_2x_3, x_1x_3)$. The Taylor resolution of $I\subset S=\KK[x_1,x_2,x_3]$ being $$ S^1
    \xleftarrow[]{\begin{psmallmatrix}x_1^2& x_1x_3& x_2x_3\end{psmallmatrix}}S^3
\xleftarrow[]{\begin{psmallmatrix}-x_3& -x_2x_3& 0\\ x_1&0&-x_2\\ 0&x_1^2&x_1\end{psmallmatrix}}S^3\xleftarrow[]{\begin{psmallmatrix}x_2\\ -1\\ x_1\end{psmallmatrix}}S^1.
    $$}\label{fig:TaylorScarf}
\end{figure}

A well-known—but highly non-minimal—construction is the \emph{Taylor resolution}, introduced in \cite{Tay66}. Given a monomial ideal $I$ with $r$ generators, its Taylor resolution is supported on an $(r-1)$-simplex known as the \emph{Taylor complex}. The vertices of this complex correspond to the generators of $I$, and its faces are labeled by the least common multiples of the incident vertices. The resolution itself mirrors the boundary complex of the Taylor complex, with monomial coefficients associated with these face labels. However, a key issue is that the \emph{repetition of labels across faces} indicates redundancies in the resolution, leading to  \emph{non-minimality}.
Various strategies exist for trimming the Taylor complex to remove these redundancies while still supporting a valid resolution. One such approach involves removing all faces that introduce repeated labels, resulting in the \emph{Scarf complex} \cite{BPS98}. However, this method does not always yield a resolution, which remains a significant challenge; see Figure~\ref{fig:TaylorScarf} for an example.

Although the non-minimality of the Taylor resolution is widely acknowledged, its extent has largely been discussed qualitatively rather than quantitatively. 
One of the goals of this work is to provide a quantitative perspective on this issue. To this end, we study dominant ideals (See Definition~\ref{def:dominant}):  monomial ideals each of whose minimal generators has a variable with exponent larger than in any other; in other words,  each generator ``dominates" in some variable. 
 \cite{Dominant_Taylor} shows that dominant ideals are precisely the monomial ideals for which the Taylor resolution is minimal.  This complete characterization of when the Taylor resolution attains minimality highlights dominant ideals as a structurally significant subclass of monomial ideals.

 \medskip 

 Given this characterization, a natural question arises: how frequently do dominant ideals occur among all monomial ideals? Since dominant ideals represent the only cases where the Taylor resolution is minimal, understanding their prevalence allows us to measure how often the Taylor resolution is minimal. The main motivation of this work is to quantify what is meant by the Taylor resolution being ``highly non-minimal" by studying the frequency of dominant ideals under various constraints, such as fixing the number of variables or the degrees of monomials in a polynomial ring. By doing so, we try to provide a somewhat concrete measure of how rare or common it is for a monomial ideal to have a minimal Taylor resolution, offering a more precise perspective on the scale of the Taylor resolution’s non-minimality.

This line of work aligns with recent efforts to incorporate probabilistic approaches in commutative algebra (\cite{Asymp_syzygies, dLPSSW, Average_Min_Res, Asymptotic_degree}). Related work in \cite{Average_Min_Res} investigates the \emph{average behavior of minimal free resolutions of monomial ideals} under a probabilistic model.
For example, their results suggest that, under a natural probability distribution, minimal free resolutions of monomial ideals almost always have {maximal length}, reinforcing the commutative algebra folklore idea that Cohen-Macaulayness is a ``rare" property. 

\medskip

Dominant ideals achieve maximal length in a free resolution. Therefore, we also study how a dominant sequence is related to the projective dimension in Section~\ref{section:max pdim}.  We give two examples that show that more structure on the ideal is required to give upper or lower bounds.

\medskip 

In the spirit of collaboration upon which the WICA workshops are founded, we use two major computer algebra systems in our computations, namely, \texttt{Macaulay2} \cite{M2} and \texttt{CoCoA} \cite{CoCoA}. For reproducibility of the simulation results, our code is made available here. Code that constructs dominant ideals is provided in Section~\ref{sec:basics and computation}, while supplementary code to produce examples and data is included in the Appendix.

%DIFDELCMD < %%%
 \subsection{The basics: dominant ideals and computation} \label{sec:basics and computation}

\begin{definition}[Dominance]\label{def:dominant}
Let $S=\mathbb{K}[x_1,\dots ,x_n]$. 
For any monomial $u$, the $x_i$-\textit{degree} of $u$ is the highest power of $x_i$ dividing $u$. We will denote the $x_i$-degree of $u$ by $\deg_{x_i}u$. Let $L$ be a set of monomials in $S$. An element $v\in L$ is called a \textit{dominant monomial} in $L$ if there exists a variable $x_i$ such that  $\deg_{x_i}v>\deg_{x_i}w$ for all $w\in L\setminus\{v\}$. In this case, we say that $v$ is dominant in $x_i$ and $x_i$ is a \textit{dominant variable} for $v$. If every monomial of $L$ is dominant in $L$, then we call $L$ a \textit{dominant set}. 

For any monomial ideal $I$, we will denote its minimal monomial generating set by $G(I)$. We say $I$ is a \textit{dominant ideal} if $G(I)$ is a dominant set.
\end{definition}

\begin{remark}\label{rem:max-numgens-for dominant}
 Note that a dominant ideal $I$ in $S$ can have at most $n$ minimal generators. If there is only one minimal generator, then the ideal is trivially dominant. 
\end{remark}

\begin{remark}\label{rem:empty-is-dominant}
We consider the zero ideal, with $L=\emptyset$, to be  dominant. 
\end{remark}

\medskip 

Before giving some examples, we offer two snippets of code in  \texttt{Macaulay2} \cite{M2} in  \texttt{CoCoA} \cite{CoCoA} 
implementing the definition of dominant ideals.
The predicate function (\texttt{isDominant} in \texttt{M2} or \texttt{IsDominant} in \texttt{CoCoA}) takes as input an ideal $I\subset R$ and returns a boolean value telling whether $I$ is dominant or not. 
If~$I$ is not a monomial ideal, then the function gives an error. Otherwise, it first deals with the trivial cases in which $G$ has either cardinality $\le1$ or $>n$. 
Finally, it applies the actual definition, verifying that all generators in $G$ are dominant.

\begin{lstlisting}[language=code, caption={Macaulay2 predicate function for dominant monomial ideals.}, label=M2dominantF]
isDominant = (I)->(
    if isMonomialIdeal(ideal mingens I) == false then(
	print "I must be a monomial ideal";
	return false;
	)else(
	R := ring I;
	G := flatten entries mingens I;
	if #G<=1 then return true;
	if #G>#(gens R) then return false else(
            lcmG := lcm(G);
	    for g in G do
                if lcmG == lcm(delete(g,G)) then return false;
	    );
	);
    return true;
    )
\end{lstlisting}

\goodbreak

\begin{lstlisting}[language=cocoa,caption={CoCoA predicate function for dominant monomial ideals.},label=COCOAdominantF]
define IsDominant(I)
  if not(AreGensMonomial(I))
    then ERROR("gens of I must be monomials");
  endif;
  R := RingOf(I);
  G := MinGens(I);
  if len(G) <= 1 then return true; endif;
  if len(G) > NumIndets(R) then return false; endif;
  lcmG := lcm(G);
  foreach g in G do
    if lcmG = lcm(diff(G,[g])) then return false; endif;
  endforeach;
  return true;
enddefine;
\end{lstlisting}

To illustrate the implementation of the functions in listings \ref{M2dominantF} and \ref{COCOAdominantF}, we provide two examples that show how to they can be use to detect dominance for a monomial ideal.
\begin{example}[M2 and CoCoA-5 implementation]
After loading the function \texttt{isDominant}, we define a polynomial ring $R$ in three indeterminates $x,y,z$ and we consider the following ideals: $I=( x^2y,xz^3,yz )$,  $J=( x^2y,xz^3,y^z )$, and $K=( x^2,y^2,x^3zy) $.
\begin{lstlisting}[language=code]
R = QQ[x,y,z];
I = ideal(x^2*y, x*z^3, y*z);
isDominant(I)
J = ideal(x^2*y, x*z^3, y^2*z);
isDominant(J)
K = ideal(x^2, y^2, x^3*z*y^2);
isDominant(K) 
\end{lstlisting}
The syntax in CoCoA-5 is fairly similar:
\begin{lstlisting}[language=cocoa]
use R ::= QQ[x,y,z];
I := ideal(x^2*y, x*z^3, y*z);
IsDominant(I); --> false
J := ideal(x^2*y, x*z^3, y^2*z);
IsDominant(J); --> true
K := ideal(x^2, y^2, x^3*z*y^2);
IsDominant(K); --> true
\end{lstlisting}
The output immediately reveals that $I$ is not dominant while $J$ is. Notice that the ideal $K$ is also dominant, since it is minimally generated by $x^2$ and $y^2$.
\end{example}

%-----------------------------------------------------------
\section{Monomial ideals with maximum projective dimension} \label{section:max pdim}
 %-----------------------------------------------------------

We are interested in detecting how data about dominant subsets of the minimal generating set    can be used to study invariants of the monomial ideal.  In particular, we focus  on projective dimension. 

By Hilbert's Syzygy Theorem, the projective dimension of a monomial ideal in $n$ variables can be at most $n$. Alesandroni \cite{A} characterized all monomial ideals which have the maximum projective dimension $n$. In this section, we determine when a monomial ideal has an associated prime of a given height. As a consequence, we recover Alesandroni's result.

Adopt the notation from Definition~\ref{def:dominant}. In particular, recall that  $G(I)$ denotes the minimal generating set of the monomial ideal $I\in\mathbb K[x_1,\dots,x_n]$. 
\begin{theorem}\label{thm:associated primes}
    For any monomial ideal $I\subset S=\mathbb{K}[x_1,\dots ,x_n]$ and $1\leq k \leq n$, the following statements are equivalent:
    \begin{enumerate}
        \item There is a dominating set $L=\{u_1,\dots ,u_k\}\subset G(I)$ such that $x_{i_j}$ is a dominant variable for $u_j$ for each $j=1,\dots, k$ and $L$ satisfies the property that
        for every $w\in G(I)\setminus L$ there exists $t\in \{1,\dots ,k\}$ such that $\deg_{x_{i_t}}w\geq \deg_{x_{i_t}}u_t$.
        \item $I$ has an associated prime of height $k$.
    \end{enumerate}
\end{theorem}
\begin{proof}
    (1) $\implies$ (2). Let $v=\lcm(u_1,\dots ,u_k)/(x_{i_1}x_{i_2}\dots x_{i_k})$ and let $\overline{v}$ be the image of $v$ in $S/I$. We claim that $\ann(\overline{v})=(x_{i_1},\dots ,x_{i_k})$. Since $L$ is a dominating set, for any $j\in\{1,\dots ,k\}$ the monomial $u_j$ divides $x_{i_j}v$. Since $L\subset G(I)$, it follows that $(x_{i_1},\dots ,x_{i_k})\subseteq \ann(\overline{v})$.
    To show the reverse containment, we use the fact that the associated prime ideals of a monomial ideal are generated by subsets of variables, see \cite[Corollary 1.3.9]{HHbook}. Assume for a contradiction there is a variable $x_r\notin \{x_{i_1}, \dots ,x_{i_k}\}$ with $x_r\in \ann(\overline{v})$. Then there is an element in $G(I)$, say $z$, which divides $x_rv$. For any $1\leq j \leq k$, the $x_{i_j}$-degree of $x_rv$ is less than that of $u_j$. Therefore, $z\in G(I)\setminus L$. Then by the assumption, there is a $t\in \{1,\dots ,k\}$ such that $\deg_{x_{i_t}}z\geq \deg_{x_{i_t}}u_t$. This implies that the $x_{i_t}$-degree of $x_rv$ is at least the $x_{i_t}$-degree of $u_t$, which is a contradiction.

    (2) $\implies$ (1). Suppose that $I$ has an associated prime of height $k$. After relabeling the indices, we may assume by \cite[Corollary 1.3.9]{HHbook} that $(x_1,\dots ,x_k)\in \ass(S/I)$. By \cite[Corollary 1.3.10]{HHbook} there is a monomial $v$ such that $(x_1,\dots ,x_k)=I:(v)$. By \cite[Proposition 1.2.2]{HHbook} there exists $u_1,\dots ,u_k\in G(I)$ such that
    \[G(I:(v))=\Big\{\frac{u_j}{\gcd(u_j,v)}:  j=1,\dots ,k\Big\}.\]
    Without loss of generality, suppose that $u_j/\gcd(u_j,v)=x_j$ for every $j=1,\dots ,k$. This shows that $L=\{u_1,\dots ,u_k\}$ is a dominating set and $x_j$ is a dominant variable for $u_j$. Moreover it follows that 
    \[\deg_{x_j}u_j=\deg_{x_j}v+1\]
    for every $j\in \{1,\dots ,k\}$. Now let $w\in G(I)\setminus L$. Since $w/\gcd(w,v)\in I$, the variable $x_t$ divides $w/\gcd(w,v)$ for some $t\in\{1,\dots ,k\}$. This implies $\deg_{x_t}w>\deg_{x_t}v$ and thus
    \[\deg_{x_t}w\geq \deg_{x_t}v+1=\deg_{x_t}u_t\]
    as desired.
\end{proof}

For two monomials $u$ and $v$ in $S$, we say $u$ \textit{strongly divides} $v$ if for every variable $x_i$ which divides $u$ one has $\deg_{x_i}v>\deg_{x_i}u$.

\begin{corollary}\cite[Corollary 5.3]{A} \label{cor:max_pdim}
    For any monomial ideal $I\subset S=\mathbb{K}[x_1,\dots ,x_n]$, the following statements are equivalent:
    \begin{enumerate}
        \item $\pd(S/I)=n$.
        \item $G(I)$ contains a dominant set $L$ of cardinality $n$, such that no element in $G(I)$ strongly divides $\lcm(L)$.
    \end{enumerate}
\end{corollary}
\begin{proof}
    By the Auslander-Buchsbaum formula $\pd(S/I)=n$ if and only if $\depth(S/I)=0$. The latter is equivalent to the fact that the maximal homogeneous ideal belongs to the associated primes of $S/I$. In this case, by Theorem~\ref{thm:associated primes} there is a dominating set $L=\{u_1,\dots ,u_n\}$ such that $x_i$ is a dominant variable for $u_i$ and for every $w\in G(I)\setminus L$ there exists $t\in \{1,\dots ,n\}$ such that $\deg_{x_t}w\geq \deg_{x_t}u_t$. It is clear that no element of $L$ strongly divides $\lcm(L)$ because the $x_i$-degrees of $L$ and $u_i$ are equal. On the other hand, if $w\in G(I)\setminus L$ and $\deg_{x_t}w\geq \deg_{x_t}u_t$, then $\deg_{x_t}w\geq \deg_{x_t}\lcm(L)$ and thus $w$ does not strongly divide $\lcm(L)$ either.
\end{proof}

In Corollary \ref{cor:max_pdim}, we observe that the projective dimension reaches its maximum value of $n$ when there is a dominant subset of cardinality $n$. We thus pose the natural general question: \emph{Is the maximum size of all dominant subsets always equal to the projective dimension?} The following example demonstrates that the answer is no, indicating that the strong divisibility condition in Corollary \ref{cor:max_pdim} (2) is essential and cannot be omitted.

\begin{example}
    Let $I = (a^3b^2, b^3c^2, a^2c^3, abc) \subset S = \mathbb{K}[a,b,c]$ and $L = \{a^3b^2, b^3c^2, a^2c^3\}$. Observe that $L$ is a dominant subset of $G(I)$ with maximum size, yet  $\text{proj dim } S/I = 2$.
\end{example}

The example above is also treated by Alessandroni  as the ideal $I$ is semi-dominant, and therefore the scarf complex is a projective resolution. The projective dimension of $I = (a^3b^2, b^3c^2, a^2c^3)$ is $3$ and the projective resolution is given by the Taylor complex. 
One might hope for a general ideal $I$ that the maximal size of a dominant subset of generators would provide an upper bound for $\text{proj dim }S/I$. The following example shows that this is not the case. 

\begin{example} Let $I=(ab, ac, bd, cd) \subset \mathbb{K}[a,b,c,d]$. Then $\text{proj dim} (S/I)=3$ but the maximum size of a dominant subset of the generators is $2$, e.g. $L=\{ab,cd\}$.
\end{example}

There are some classes of ideals for which maximum size of dominant subsets does equal the projective dimension. For example, suppose that $L$ is a dominant subset of $G(I)$ of maximal size $k$, so that the dominant variables are $x_1,\dots,x_k$. Then the projective dimension would be bounded above by $k$ if $x_i\nmid m$ for any $m\in G(I)\setminus L$, $i=1,\dots,l$. 
To see why this bound holds, note that localizing at the non-dominant variables would eventually result in the localized $I$ having the form in Corollary~\ref{cor:max_pdim}. Localization is exact, therefore it does not change the projective dimension.

%-----------------------------------------------
\section{Dominant ideals with a given least common multiple}\label{sec:lcm}
%-----------------------------------------------

The least common multiple of the minimal generators is crucial when determining if an ideal is dominant (see the two functions {\tt isDominant()} in the previous section). 
This section investigates dominant ideals with a given least common multiple (lcm). Fixing the lcm  will enable enumeration of dominant ideals using a combinatorial formula. 

\begin{definition}
We define the \textbf{lcm of a monomial ideal} $I \subset R$ to be the lcm of the minimal set of monomial generators of~$I$.
\end{definition}
\begin{example}
    The ideal $I = ( x^2, \;x^2z, \;y^3)$
    is minimally generated by $\{x^2,\; y^3\}$, hence
    is dominant, and its lcm is $\thelcm=x^2y^3$. 
\end{example}

\begin{example}
    The ideal $I = ( x^2y^3, \;x^2z^4)$ is dominant, and its lcm is $\thelcm=x^2y^3z^4$.
    The first generator is dominant in $y$ and the second in $z$.
    No generator is dominant in $x$.
\end{example}
We enumerate dominant ideals with a given lcm $\thelcm$ up to 5 variables (see Subsection~\ref{sec:enumerating-dominant-lcm}). These formulas are interesting because they allow us to compare the frequency of dominant ideals with the number of all monomial ideals with lcm $\thelcm$. We close this section with a set of open problems.

\subsection{Enumerating all dominant ideals with a given least common multiple}
\label{sec:enumerating-dominant-lcm}

The complexity of the formula increases very quickly as the number of variables grows.
We tackle just the first few cases. The formulae we compute were actually checked by computations in CoCoA-5, as described in Section~\ref{sec:computing-formula-all-dominant-lcm}.

Recall that, for dominance to occur, the number of minimal generators must not exceed the number of variables (see Remark~\ref{rem:max-numgens-for dominant}).
Therefore, we divide the enumeration by counting the number of dominant ideals with lcm =$\thelcm$, while the number of minimal generators varies from 0 to $n$.

The first case we consider is quite easy.

\begin{proposition}
    	 The number of dominant ideals in $\KK[x,y]$ with $\thelcm=x^{m_1}y^{m_2}$ is $1+m_1m_2$.
\end{proposition}
\begin{proof}
Consider a monomial ideal $I\subset R$ with lcm $\thelcm$.

Clearly, if $I$ has 1 minimal generator then $I = \langle \thelcm\rangle$.

If $I$ has 2 minimal generators, then it must have the following form
$$
( x^{m_1}y^j,\;x^iy^{m_2} ), \hbox{ where } 0\leq i< m_1, \; 0\leq j< m_2.
$$	
Therefore there are $m_1\cdot m_2$ dominant ideals with 2 minimal generators and lcm $=\thelcm$.
Hence, the total number of dominant ideals with lcm $=\thelcm$ is $2+m_1m_2$.	
\end{proof}

Next we consider the case of 3 variables.
Here have more cases to consider, so we introduce a partition based on a system of labels for families of sub-cases, which emphasize the variables whose exponents is non-maximal, and therefore may vary in their respecting range $0,...,m_i$.
This way of partitioning is the key for counting all dominant ideals with a given lcm.

Before stating and proving the theorem,
we illustrate through some examples the partition for the dominant ideals $\KK[x,y,z]$ with lcm = $x^2y^3z^4$,  having at least 2 generators.

We also introduce in these examples the concept of the \textbf{footprint} of an ideal with a given lcm, encoding the \textbf{non-maximal variables} for the generators, which is the key for partitioning and counting the dominant ideals.

\begin{definition}
Let $I$ be a monomial ideal in $K[x_1,...,x_n]$ with lcm $\thelcm= x_1^{m_1}\cdots x_n^{m_n}$, and let $G$ be 
the minimal set of monomial generators of $I$.
We define the \textbf{footprint} of a monomial $g=x_1^{g_1}\cdots x_n^{g_n}\in G$ as the product of the $x_i$'s such that $g_i<m_i$. 
Then, the \textbf{footprint} of $I$ is the set of footprints of the minimal generators.
\end{definition}

\begin{example}\label{ex:dominant3-with3gens}
For a dominant ideal $I$ with three generators in $\KK[x,y,z]$ each generator must have maximal exponent in one, and only one, variable.  Then the footprint of $I$ must be $\{yz, xz, xy\}$, indicating the ``free'' non-maximal variables for each generator:
\\\phantom.\hfil
$( x^2y^{j_1}z^{k_1}, \;\; x^{i_1}y^3z^{k_2}, \;\;x^{i_1}y^{j_2}z^4 ) $

The possible generators are graphically represented as the blue dots in the picture, one for each face.
There are $(3\cdot4)(2\cdot4)(2\cdot3)=m_1^2m_2^2m_3^2 = 576$ many choices.

\begin{figure}[H]    \label{fig:3var3support}
    \includegraphics[width=0.45\linewidth]{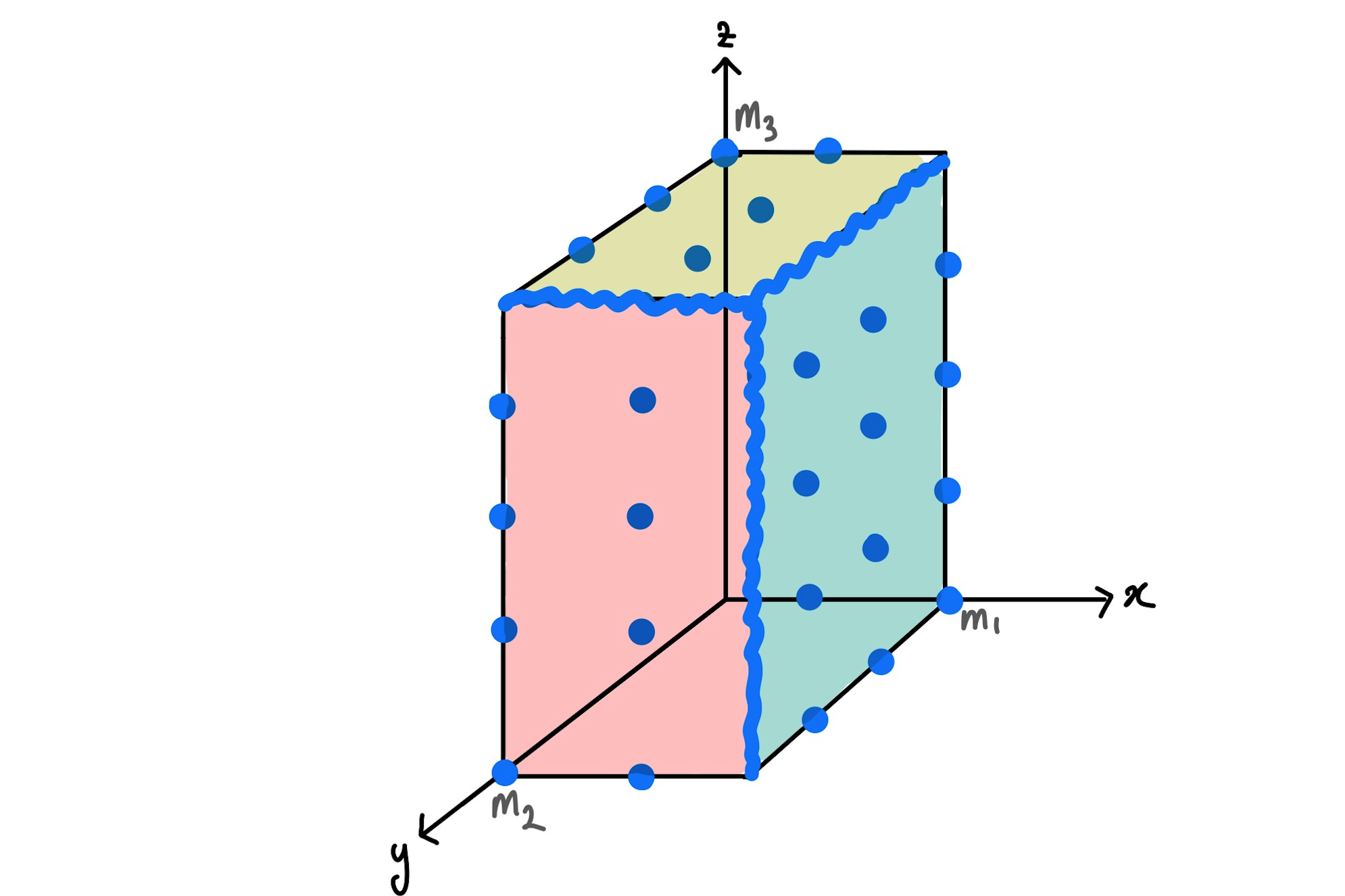}
\end{figure}
\end{example}

For a dominant ideal with two generators, again with lcm = $x^2y^3z^4$,
each generator must have maximal exponent in one or two variables.

\begin{example}\label{ex:dominant3-with2gensyz}
The first family is partitioned by the following 3 footprints:
\begin{itemize}
\item$\{z,y\}$ %\verb|[z,  y]|:
i.e. the first generator is of type $x^2y^3z^k$ and the second is of type $x^2y^jz^4$, both maximal in $x$, thus no generator is dominant in $x$.  There are $3\cdot4$ ideals, see Fig.~\ref{fig:3var2support_1} left;
\item$\{z,x\}$ %\verb|[z,  x]|:
both generators maximal in $y$.
There are $2\cdot4$ ideals, Fig.~\ref{fig:3var2support_1} middle;
\item$\{y,x\}$ %\verb|[y,  x]|:
both generators maximal in $z$.  There are $2\cdot3$ ideals,
Fig.~\ref{fig:3var2support_1} right.
\end{itemize}
There are $12 + 8 + 6$ such ideals.
The choices of the generators are represented by the blue dots on the two coloured edges:
\begin{figure}[H]
    \centering
    \includegraphics[width=0.9\linewidth]{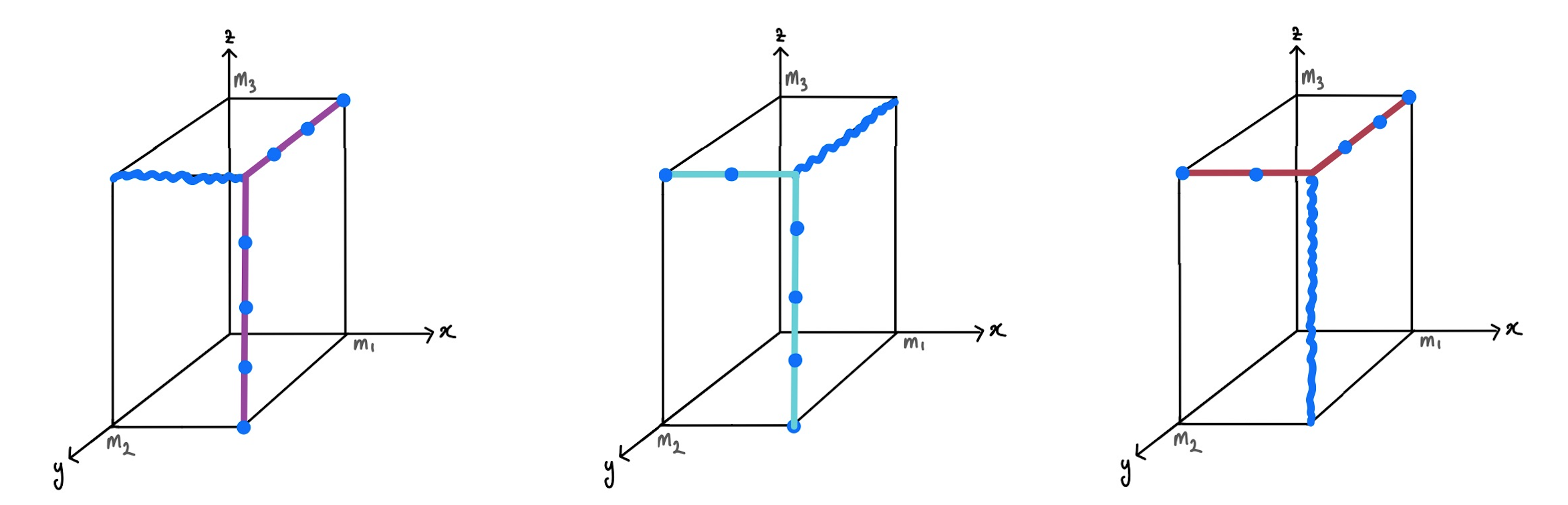}
    \caption{}\label{fig:3var2support_1}
    
\end{figure}
\end{example}

\begin{example}\label{ex:dominant3-with2gensxyz}
The second family for 2 generators is partitioned by the following 3 footprints:
\begin{itemize}
    \item  $\{x, \; yz\}$: on generator is of type $x^iy^3z^4$, and is dominant in $y$ and $z$, and the other is of type $x^2y^jz^k$, dominant in $x$.  
    There are $2\cdot(3\cdot 4)$ ideals, Fig.~\ref{fig:3var2support_2} left;
    \item  $\{y, \; xz\}$:
    There are $3\cdot(2\cdot4)$ ideals, Fig.~\ref{fig:3var2support_2} middle;
    \item  $\{z, \; xy\}$:
    There are $4\cdot(2\cdot3)$ ideals, Fig.~\ref{fig:3var2support_2} right.
\end{itemize}
There are $24 + 24 + 24 = 72$ such ideals.
The choices of the generators are represented by the blue dots, the first on the edge, the second on the face:
\begin{figure}[H]
    \includegraphics[width=0.9\linewidth]{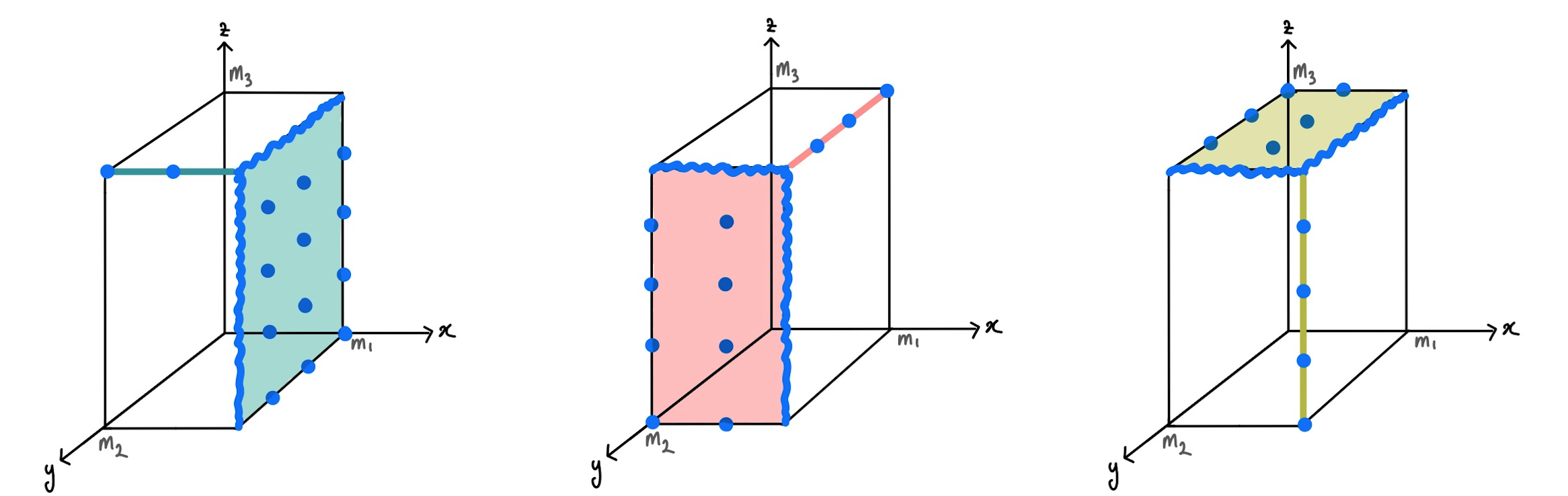}
 \caption{}\label{fig:3var2support_2}
\end{figure}
\end{example}

\begin{example}
We can compute (see Section~\ref{sec:compute-all-dominant-lcm}) the partitions of the dominant ideals with lcm=$x^2y^3z^4$
in two ways. 
The finer partition by the footprint:
\begin{lstlisting}[language=cocoa]
  record[count := 576, footprint := [y*z,  x*z,  x*y]],
  record[count := 24, footprint := [z,  x*y]],
  record[count := 24, footprint := [y,  x*z]],
  record[count := 12, footprint := [z,  y]],
  record[count := 8, footprint := [z,  x]],
  record[count := 6, footprint := [y,  x]],
  record[count := 24, footprint := [x,  y*z]],
  record[count := 1, footprint := [1]]
\end{lstlisting}
or the coarser partition just counting of the variables with \textbf{l}ower/\textbf{m}aximal exponents for each generator:
\begin{lstlisting}[language=cocoa]
  record[LowOrMax := [l^2*m,  l^2*m,  l^2*m], count := 576],
  record[LowOrMax := [l*m^2,  l^2*m], count := 72],
  record[LowOrMax := [l*m^2,  l*m^2], count := 26],
  record[LowOrMax := [m^3], count := 1]
\end{lstlisting}
\end{example}

We are now ready for giving the general formula which enumerates the segments in the partition of all dominant ideals in 3 variables with a given lcm.

\begin{theorem}\label{thm:dominant-lcm-3}
	 The number of dominant ideals in $\KK[x,y,z]$ with $\thelcm=x^{m_1}y^{m_2}z^{m_3}$ is
	 $$
     1 \;
     +m_1^2m_2^2m_3^2 \;
     +m_1m_2 +m_1m_3 +m_2m_3 \;
     + 3m_1m_2m_3 .
	 $$
\end{theorem}

\begin{proof}
Clearly, if $I$ has 1 minimal generator, then $I = ( \thelcm)$. 

%% -- 3 generators n=3 ----------------
$I$ has 3 minimal generators then it must have the following form
$$
( x^{m_1}{y^{j_1}z^{k_1}}, \; {x^{i_2}y^{m_2}}z^{k_2}, \; x^{i_3}y^{j_3}z^{m_3} )
\\\text{ with } 
0\le i_2,i_3< m_1, \quad
0\le j_1,j_3< m_2, \quad 0\le k_1,k_2< m_3.
$$
We label this situation as \verb|[l^2*m,  l^2*m,  l^2*m]|, 
one \textbf{m}aximal power for each generator, 
with one footprint $\{yz, \; xz, \; xy\}$.
For example, for the first generator, the exponents of $y,z$  may vary in their respective ranges, giving $m_2\cdot m_3$ possible choices for the first generator.

Therefore there are $(m_2m_3)(m_1m_3)(m_1m_2) = m_1^2m_2^2m_3^2$ 
dominant ideals with 3 minimal generators and lcm $=\thelcm$.
(see Example~\ref{ex:dominant3-with3gens}).

%% -- 2 generators n=3 ----------------
We now have to consider the ideals with 2 minimal generators, recalling they must be dominant and with lcm =$\thelcm$.
We distinguish two families of sub-cases.

First family: \verb|[l*m^2,  l*m^2]|,  with possible  footprints
%\verb|[y,  x]|, \verb|[z,  x]|, or \verb|[y,  z]|.
$\{y,x\}, \{z,x\}, \{y,z\}$.
\\
Consider $\{y,x\}$, i.e. the ideals of the form 
$$( x^{m_1}\underline{y^j}z^{m_3}, \underline{x^i}y^{m_2}z^{m_3} ) 
\text{\quad with \quad}
0\leq i< m_1 \text{\; and \;} 
0\leq j< m_2,$$
so there are $m_1m_2$ of them. 
Similarly there are $m_1m_3$  ideals labelled \verb|[z, x]|,
and $m_2m_3$  ideals of type \verb|[y, z]|.
In total $m_1m_2+ m_1m_3+m_2m_3$.
See Example~\ref{ex:dominant3-with2gensyz}.

Second family: \verb|[l*m^2,  l^2*m]|, with possible footprints
$\{z, xy\}, \{x, yz\},\text{ or }\{y, xz\}$.
\\
Consider $\{z, xy\}$, representing the ideals of the form 
$$( x^{m_1}y^{m_2}{z^k},\quad{x^{i}y^i}z^{m_3} )  
\text{\quad with \quad}
0\leq i < m_1,\quad 0\leq j< m_2,\quad 0\leq k< m_3.$$
There are $m_1 m_2 m_3$ of them. 
Similarly there are $m_1m_2m_3$  ideals with footprint
$\{x, yz\}$ and  $\{y, xz\}$.
In total $3\,m_1m_2m_3$.
See Example~\ref{ex:dominant3-with2gensxyz}.
\end{proof}

%% 4 variables -------------------------------------------
\begin{theorem}\label{thm:dominant-lcm-4}
	 The number of dominant ideals in $\KK[x_1,\dots,x_4]$ with $\thelcm=x_1^{m_1}\cdots x_4^{m_4}$ is
\begin{gather*}
            1+(m_1m_2m_3m_4)^3+\sum_{i<j}m_im_j+\\
3(m_1m_2m_3+m_1m_2m_4+m_1m_3m_4+m_2m_3m_4)+\\
7m_1m_2m_3m_4+(m_1m_2m_3)^2+(m_1m_2m_4)^2+
  (m_1m_3m_4)^2+(m_2m_3m_4)^2+ \\ 
6(m_1m_2m_3m_4)^2+ 
3(m_1m_2^2m_3^2m_4^2+m_1^2m_2m_3^2m_4^2+m_1^2m_2^2m_3m_4^2+m_1^2m_2^2m_3^2m_4).
\end{gather*}
 \end{theorem}

\begin{proof}
These are the cases to consider:
\begin{lstlisting}[language=cocoa]
/*4gens*/  [l^3*m,  l^3*m,  l^3*m,  l^3*m]
/*3gens*/  [l^2*m^2,  l^2*m^2,  l^3*m]
           [l^2*m^2,  l^2*m^2,  l^2*m^2]
           [l^2*m^2,  l^3*m,  l^3*m]
/*2gens*/  [l^2*m^2,  l^2*m^2]
           [l*m^3,  l^2*m^2]
           [l*m^3,  l*m^3]
           [l*m^3,  l^3*m]
/*1gen*/   [m^4]
\end{lstlisting}

There is only the ideal $\langle\thelcm\rangle$ with 1 generator, and this gives the 1 in the formula. 

%% -- 4 generators n=4 ----------------
The ideals with 4 minimal generators have exactly one maximal power each, so their labels are \verb|l^3*m,  l^3*m,  l^3*m,  l^3*m|, with footprint
$\{x_2x_3x_4,\;  x_1x_3x_4,\;  x_1x_2x_4,\;  x_1x_2x_3\}$.

Thus there are  $(m_2 m_3 m_4)(m_1 m_3 m_4)(m_1 m_2 m_4)(m_1 m_2 m_3) = 
(m_1m_2m_3m_4)^3 $ ideals with 4 generators.

Consider now the ideals with 2 generators.
There are 4 cases:
\begin{itemize}
\item \verb|[l*m^3,  l*m^3]|, more precisely there are 6 footprints: %subcases for the lower exponents
 \\\phantom.\hfil
 $
  \{x_2,  x_1\},    \;\;
  \{x_3,  x_1\} ,   \;\;
  \{x_4,  x_1\}  ,  \;\;
  \{x_3,  x_2\}   , \;\;
  \{x_4,  x_2\}    ,\;\;
  \{x_4,  x_3\}    .
  $
  \\
    For example, for the footprint $\{x_2,  x_1\} $ the two generators of the ideals are
\\      $( x_1^{m_1}x_2^{j}x_3^{m_3}x_4^{m_4}, \;
    x_1^{i}x_2^{m_2}x_3^{m_3}x_4^{m_4})$, 
    with $0{\le} i{<}m_1$ and $0{\le} j{<}m_2$. 
    Hence, there are $m_1m_2$ of them. 

    Similarly for the other footprints, 
    then there are $\sum_{i<j}m_im_j$ such dominant ideals. 

    \item   \verb|[l*m^3,  l^2*m^2]| % count := 12], 
more precisely there are 12 footprints
\\\phantom.\hfil
$  \{x_1,  x_2x_3\} ,   \;
  \{x_1,  x_2x_4\}   , \;
  \{x_1,  x_3x_4\}    ,\;
  \{x_2,  x_1x_4\}  ,  \;
  \{x_2,  x_1x_3\}  ,  \;
  \{x_2,  x_3x_4\}  ,  $
  \\\phantom.\hfil
$  \{x_3,  x_1x_2\} ,   \;
  \{x_3,  x_1x_4\}  ,  \;
  \{x_3,  x_2x_4\}  ,  \;
  \{x_4,  x_1x_2\}  ,  \;
  \{x_4,  x_1x_3\}  ,  \;
  \{x_4,  x_2x_3\}  .  \;
  $
  
    For example, with footprint $\{x_1,  x_2x_3\}$ the dominant ideals 
\\    $( x_1^{i}x_2^{m_2}x_3^{m_3}x_4^{m_4}, \;\;
    x_1^{m_i}x_2^{j}x_3^{k}x_4^{m_4})$, 
    with $0{\le} i{<}m_1$,\; $0{\le} j{<}m_2$,\;  $0{\le} k{<}m_3$ and there are $m_1m_2m_3$ of them.         
And similarly for the other footprints, 
    then there are 
    $ 3(m_1m_2m_3\;+m_1m_2m_4\;+m_1m_3m_4\;+m_2m_3m_4)$ 
    such dominant ideals. 
    
    \item \verb|[l^2*m^2,  l^2*m^2]|  %, count := 3],
more precisely there are 3 footprints
%subcases for the lower exponents
\\\phantom.\hfil      $\{x_1x_4,  x_2x_3\} ,   \quad
    \{x_2x_4,  x_1x_3\} ,   \quad
\{x_3x_4,  x_1x_2\} .   $
\\
   For example, for $\{x_1x_4,  x_2x_3\} $
   the dominant ideals are of the form
   \\
    $( x_1^{m_1}x_2^jx_3^k x_4^{m_4},\; 
    x_1^ix_2^{m_2}x_3^{m_3}x_4^\ell )$, 
        with $0{\le} i{<}m_1$,\; $0{\le} j{<}m_2$,\;  $0{\le} k{<}m_3$,
        \;  $0\le \ell{<}m_4$
and there are $m_1m_2m_3m_4$ of them.  
So, there are $ 3m_1m_2m_3m_4$ such dominant ideals.

    \item \verb|[l*m^3,  l^3*m]| %, count := 4]
more precisely there are 4 footprints
\\\phantom.\hfil      $
  \{x_1,  x_2x_3x_4\},      \quad
  \{x_2,  x_1x_3x_4\} ,   \quad
\{x_3,  x_1x_2x_4\}    ,\quad
  \{x_4,  x_1x_2x_3\}   . 
$
\\  For example, for $\{x_1,  x_2x_3x_4\}  $
   the dominant ideals are of the form
   \\
    $( x_1^ix_2^{m_2}x_3^{m_3} x_4^{m_4},\; 
    x_1^{m_1}x_2^jx_3^kx_4^\ell) $, 
        with $0{\le} i{<}m_1$,\; $0{\le} j{<}m_2$,\;  $0{\le} k{<}m_3$,
        \;  $0{\le} \ell{<}m_4$
    then, in total, there are $ 4m_1m_2m_3m_4$ such dominant ideals.
\end{itemize}

Lastly, consider the ideals with 3 generators. There are 3 cases:

\begin{itemize} %%%%%%%%%%%%%%%%%%%%%%%%%%%%%%%%%%%%%%%%%%%%%%%%
\item \verb|[l^2*m^2,  l^2*m^2,  l^2*m^2]|  %, count := 4],
more precisely there are 4 footprints
  \\\phantom.\hfil
$ 
  \{x_2x_3,  x_1x_3,  x_1x_2\},    
  \{x_3x_4,  x_1x_4,  x_1x_3\} ,   
  \{x_2x_4,  x_1x_4,  x_1x_2\}  ,  
  \{x_3x_4,  x_2x_4,  x_2x_3\}   . 
     $
\\ 
For example, for $\{x_2x_3,  x_1x_3,  x_1x_2\}   $ we have the ideals
\\\phantom.\hfil
$( 
x_1^{m_1}x_2^{j_1}x_3^{k_1}x_4^{m_4}, \quad
x_1^{i_2}x_2^{m_2}x_3^{k_2}x_4^{m_4}, \quad
x_1^{i_3}x_2^{j_3}x_3^{m_3}x_4^{m_4})$
\\
with $ 0{\le}i_2,i_3{<} m_1,\;\;0{\le}j_1,j_3{<} m_2,\;\;0{\le}k_1,k_2{<} m_3 $,
so there are $m_1^2m_2^2m_3^2$ of them. 

Then we have $(m_1m_2m_3)^2+(m_1m_2m_4)^2+  (m_1m_3m_4)^2+(m_2m_3m_4)^2$ dominant ideals of this form.

  \item \verb|[l^2*m^2,  l^2*m^2,  l^3*m]| %, count := 12],
more precisely there are 12 footprints
\\\phantom.\hfil
$
  \{x_3x_4,  x_2x_4,  x_1x_2x_3\} ,   \;
  \{x_3x_4,  x_1x_4,  x_1x_2x_3\}  ,  \;
  \{x_2x_4,  x_1x_4,  x_1x_2x_3\}   , \;
  \{x_2x_3,  x_1x_3,  x_1x_2x_4\}    ,\;
$
\\\phantom.\hfil
$
  \{x_3x_4,  x_2x_3,  x_1x_2x_4\},    \;
  \{x_3x_4,  x_1x_3,  x_1x_2x_4\} ,   \;
  \{x_2x_4,  x_1x_2,  x_1x_3x_4\}  ,  \;
  \{x_2x_3,  x_1x_2,  x_1x_3x_4\}   , \;
 $
\\\phantom.\hfil
$ 
  \{x_2x_4,  x_2x_3,  x_1x_3x_4\} ,   \;
  \{x_1x_4,  x_1x_3,  x_2x_3x_4\}  ,  \;
  \{x_1x_4,  x_1x_2,  x_2x_3x_4\}   , \;
  \{x_1x_3,  x_1x_2,  x_2x_3x_4\}    .\;
 $

For example, for $ \{x_3x_4,  x_2x_4,  x_1x_2x_3\}  $
we have the ideals
\\\phantom.\hfil
$(
x_1^{m_1}x_2^{m_2}x_3^{k_1}x_4^{\ell_1},\quad
x_1^{m_1}x_2^{j_2}x_3^{m_3}x_4^{\ell_2},\quad
x_1^{i_3}x_2^{j_3}x_3^{k_3}x_4^{m_4} )
$
\\
with $ 0{\le}i_3{<} m_1,\;\;
0{\le}j_2,j_3{<} m_2,\;\;
0{\le}k_1,k_3{<} m_3,\;\;
0{\le}\ell_1,\ell_2{<} m_4
$,
so there are $m_1^1m_2^2m_3^2m_4^2$ of them, and the same number is obtained by all 3 footprints in the first column, i.e. those with only one $x_1$.

Then in total there are $3(m_1m_2^2m_3^2m_4^2+m_1^2m_2m_3^2m_4^2+m_1^2m_2^2m_3m_4^2+m_1^2m_2^2m_3^2m_4)$ dominant ideals of this kind. 

  \item \verb|[l^2*m^2,  l^3*m,  l^3*m]| %, count := 6],
more precisely there are 6 footprints
\\\phantom.\hfil
$  \{x_1x_2,  x_2x_3x_4,  x_1x_3x_4\}   \quad
  \{x_1x_3,  x_2x_3x_4,  x_1x_2x_4\}    \quad
  \{x_1x_4,  x_2x_3x_4,  x_1x_2x_3\}    
$
\\\phantom.\hfil
$  \{x_2x_4,  x_1x_3x_4,  x_1x_2x_3\}    \quad
  \{x_3x_4,  x_1x_2x_4,  x_1x_2x_3\}    \quad
  \{x_2x_3,  x_1x_3x_4,  x_1x_2x_4\}    
$
\\
For example, for $\{x_2x_4,  x_1x_3x_4,  x_1x_2x_3\} $
we have $m_1^2m_2^2m_3^2m_4^2$, and the same for the other footprints.
Then in total there are $6m_1^2m_2^2m_3^2m_4^2$ such dominant ideals. 
\end{itemize}
Since we listed all possible cases this concludes the proof.
\end{proof}

The last formula we compute is for $n=5$.
We follow the same partition and counting strategy
as in the previous proofs.

%% 5 variables -------------------------------------------
\begin{theorem}\label{thm:dominant-lcm-5}
   The number of dominant monomial ideals in $\mathbb{K}[x_1,\dots,x_5]$ with lcm $\thelcm=x_1^{m_1}\cdots x_5^{m_5}$ is
   \begin{gather*}
       1 + (m_1 m_2 m_3 m_4 m_5)^4 \\
       +\sum_{i<j} m_im_j
       \quad+3\!\sum_{i<j<k} m_im_jm_k    
       \quad+7\!\sum_{i<j<k<\ell} m_im_jm_km_\ell 
       \quad       +15\;m_1m_2m_3m_4m_5 
       \\
       + 10(m_1m_2m_3m_4m_5)^{3} 
       +\sum_{i<j<k<\ell} (m_im_jm_km_\ell)^3 
       %\textcolor{blue}{+6(m_1m_2m_3m_4m_5)^3}
        \\
       +6(m_1m_2m_3m_4m_5)^2\!\sum_{i<j<k<\ell} (m_im_jm_km_\ell)  
      \;\;+4(m_1m_2m_3m_4m_5)\!\sum_{i<j<k<\ell}(m_im_jm_km_\ell)^2  
      \\
       +25\;(m_1m_2m_3m_4m_5)^2 
              +\sum_{i<j<k} (m_im_jm_k)^2 
              +6\!\sum_{i<j<k<\ell}  (m_im_jm_km_\ell)^2  
              \\
       +9\;(m_1m_2m_3m_4m_5)\!\sum_{i<j<k} (m_im_jm_k) 
       \;+18\;(m_1m_2m_3m_4m_5)\!\sum_{i<j<k<\ell} (m_im_jm_km_\ell)  
       \\
       +3\;(
       m_1((m_2m_3m_4)^2 +(m_2m_3m_5)^2 +(m_2m_4m_5)^2 +(m_3m_4m_5)^2 )\\ \quad
       +m_2((m_1m_3m_4)^2 +(m_1m_3m_5)^2 +(m_1m_4m_5)^2 +(m_3m_4m_5)^2 )\\ \quad
       +m_3((m_1m_2m_4)^2 +(m_1m_2m_5)^2  +(m_1m_4m_5)^2  +(m_2m_4m_5)^2 )\\ \quad
       +m_4((m_1m_2m_3)^2 +(m_1m_2m_5)^2 +(m_1m_3m_5)^2 +(m_2m_3m_5)^2 )\\ \quad\quad
       +m_5((m_1m_2m_3)^2 +(m_1m_2m_4)^2 +(m_1m_3m_4)^2 +(m_2m_3m_4)^2) \quad
       )
   \end{gather*}
\end{theorem}
\begin{proof}
These are the cases to consider
\begin{lstlisting}[language=cocoa]
/*5gens*/  [l^4*m,  l^4*m,  l^4*m,  l^4*m,  l^4*m] 
/*4gens*/  [l^3*m^2,  l^4*m,  l^4*m,  l^4*m] 
           [l^3*m^2,  l^3*m^2,  l^4*m,  l^4*m] 
           [l^3*m^2,  l^3*m^2,  l^3*m^2,  l^3*m^2]  
           [l^3*m^2,  l^3*m^2,  l^3*m^2,  l^4*m]  
/*3gens*/  [l^2*m^3,  l^2*m^3,  l^2*m^3]  
           [l^2*m^3,  l^2*m^3,  l^3*m^2]  
           [l^2*m^3,  l^2*m^3,  l^4*m]  
           [l^2*m^3,  l^3*m^2,  l^3*m^2]  
           [l^2*m^3,  l^3*m^2,  l^4*m] 
           [l^2*m^3,  l^4*m,  l^4*m]  
           [l^3*m^2,  l^3*m^2,  l^3*m^2]  
           [l^3*m^2,  l^3*m^2,  l^4*m] 
/*2gens*/  [l*m^4,  l*m^4] 
           [l*m^4,  l^2*m^3]  
           [l*m^4,  l^3*m^2]  
           [l*m^4,  l^4*m]  
           [l^2*m^3,  l^2*m^3] 
           [l^2*m^3,  l^3*m^2]       
/*1gen*/   [m^5] 
\end{lstlisting}
%Fix the monomial $x_1^{m_1}\cdots x_5^{m_5}$. 
%  [m^5]  %count := 1],

Firstly, the easy cases.
There is one dominant ideal generated by 1 monomial, i.e. $\langle \thelcm \rangle $. 

%  [l^4*m,  l^4*m,  l^4*m,  l^4*m,  l^4*m]% count := 1],
The dominant ideals generated by 5 monomials have footprint
\\\phantom.\hfil
$\{x_2x_3x_4x_5, \qquad x_1x_3x_4x_5, \qquad x_1x_2x_4x_5, \qquad x_1x_2x_3x_5, \qquad x_1x_2x_3x_4\}$
\\i.e.\hfil
  $(
  x_1^{m_1}x_2^*x_3^*x_4^*x_5^*,  \;
  x_1^*x_2^{m_2}x_3^*x_4^*x_5^*,  \;
  x_1^*x_2^*x_3^{m_3}x_4^*x_5^*,  \;
  x_1^*x_2^*x_3^*x_4^{m_4}x_5^*,  \;
  x_1^*x_2^*x_3^*x_4^*x_5^{m_5}  )$,
  \\
and there are  $ (m_1\cdots m_5)^4$ of them. 

Consider now the dominant ideals with 2 generator,
with the following subcases:
\begin{itemize}
    \item \verb|[l*m^4,  l*m^4]|,  %count := 10],
    e.g. the footprint $\{x_1, x_2\}$ giving $m_1m_2$ ideals,    
%If $(a,b)=(1,1)$
then in total there are $\sum_{i<j}m_im_j$ such ideals.

    \item \verb|[l*m^4,  l^2*m^3]|  %count := 30],
    e.g.  the footprints $\{x_1, x_2x_3\}$, $\{x_2, x_1x_3\}$, or $\{x_3, x_1x_2\}$,
    each giving $m_1m_2m_3$ ideals,  
    %If $(a,b)=(1,2)$ 
    then there are in total $3\sum_{i<j<k} m_im_jm_k$ such ideals.
    
    \item \verb|[l*m^4,  l^3*m^2]|  %count := 20],|
    %If $(a,b)=(1,3)$ 
      e.g. $\{x_1, x_2x_3x_4\}$, $\{x_2, x_1x_3x_4\}$, 
      $\{x_3, x_1x_2x_4\}$,
      or $\{x_4, x_1x_2x_3\}$,
    giving $m_1m_2m_3m_4$ ideals,    
    then in total $ 4\sum_{i<j<k<\ell}m_im_jm_km_\ell$ such ideals.
    
    \item \verb|[l*m^4,  l^4*m]|  %count := 5], 
    %If $(a,b)=(1,4)$
    e.g. $\{x_1, x_2x_3x_4x_5\}$, $\{x_2, x_1x_3x_4x_5\}$, 
and so on,
    giving     in total $5\;m_1 m_2 m_3 m_4 m_5$ such ideals.

    \item  \verb|[l^2*m^3,  l^2*m^3]|  %count := 15],
    %If $(a,b)=(2,2)$ 
    e.g. the 3 footprints $\{x_1x_2, x_3x_4\}$, $\{x_1x_3, x_1x_4\}$, $\{x_2x_3, x_1x_4\}$ 
    give $m_1m_2m_3m_4$ ideals each,  
    and in total $3\sum_{i<j<k<\ell}m_im_jm_km_\ell$ such ideals.
    
    \item \verb|[l^2*m^3,  l^3*m^2]|  %count := 10],
    %If $(a,b)=(2,3)$ 
    e.g. $\{x_1x_2, x_3x_4x_5\}$, $\{x_1x_3, x_1x_4x_5\}$, and so on,
    giving in total $10\;m_1m_2m_3m_4m_5$ such ideals.
\end{itemize}

We pass now to the dominant ideals with 3 generators.

\begin{itemize}

    \item \verb|[l^2*m^3,  l^2*m^3,  l^2*m^3]|  %count := 10], Anna: mancava
    e.g. $\{x_3x_4,  x_1x_4,  x_1x_3\}$ gives $(m_1 m_3 m_4)^2$ ideals.
    \\
    In total $\sum_{i<j<k} (m_im_jm_k)^2$ such ideals.

    \item \verb|[l^2*m^3,  l^2*m^3,  l^3*m^2]|  %count := 60],
%|If $(a,b,c)=(2,2,3)$ 
e.g. $\{x_1x_3, \; x_1x_2, \; x_2x_3x_4\}$ gives $(m_1 m_2 m_3)^2m_4$ ideals.
    Then there are in total
$+3\;(
       m_1((m_2m_3m_4)^2 +(m_2m_3m_5)^2 +(m_2m_4m_5)^2 +(m_3m_4m_5)^2 ) + ... )$ such ideal.
       %    $10(m_1 m_2 m_3 m_4^2 m_5^2+\cdots + m_1^2m_2^2m_3m_4m_5)$ dominant ideals.

    \item \verb|[l^2*m^3,  l^2*m^3,  l^4*m]|  %count := 30],
    %If $(a,b,c)=(2,2,4)$ 
e.g. $\{x_1x_3,  x_1x_2,  x_2x_3x_4x_5\}$,
$\{x_2x_3,  x_1x_2,  x_2x_3x_4x_5\}$, and
$\{x_2x_3,  x_1x_2,  x_1x_3x_4x_5\}$
each giving $(m_1 m_2 m_3)^2m_4m_5$ ideals.
\\
Then in total there are 
    $ 3\;m_1m_2m_3m_4m_5 \sum_{i<j<k}m_im_jm_k$ such ideals.

    \item  \verb|[l^2*m^3,  l^3*m^2,  l^3*m^2]|  %count := 90],
    %If $(a,b,c)=(2,3,3)$
    e.g. these footprints
    $\{x_1x_2, \; x_1x_3x_5, \; x_2x_3x_4\}$,\\
    $\{x_1x_2, \; x_2x_3x_5, \; x_1x_3x_4\}$,
    $\{x_1x_3, \; x_1x_3x_5, \; x_2x_2x_4\}$,
    $\{x_1x_3, \; x_2x_3x_5, \; x_1x_2x_4\}$,\\
    $\{x_2x_3, \; x_1x_2x_5, \; x_1x_3x_4\}$,
    $\{x_2x_3, \; x_1x_3x_5, \; x_1x_2x_4\}$,
each giving $(m_1m_2m_3)^2m_4m_5 $  ideals.
Then in total there are \;
    $6\;m_1m_2m_3m_4m_5 \sum_{i<j<k}m_im_jm_k$ such ideals.

Moreover, there are other types of ideals for this label,
e.g. with footprints\\
    $\{x_1x_2, \; x_1x_3x_4, \; x_2x_3x_4\}$,
    $\{x_1x_3, \; x_1x_2x_4, \; x_2x_3x_4\}$,
    $\{x_1x_4, \; x_1x_2x_3, \; x_2x_3x_4\}$,\\
    $\{x_2x_3, \; x_1x_2x_4, \; x_1x_3x_4\}$,
    $\{x_2x_4, \; x_1x_2x_3, \; x_1x_3x_4\}$,
    $\{x_3x_4, \; x_1x_2x_3, \; x_1x_2x_4\}$,
giving each $(m_1m_2m_3m_4)^2$  ideals.
Then in total 
    $ 6 \sum_{i<j<k<\ell}(m_im_jm_km_\ell)^2$ such ideals.

    \item \verb|[l^2*m^3,  l^3*m^2,  l^4*m]|  %count := 60],
    %If $(a,b,c)=(2,3,4)$
    e.g. for the 12 footprints
    $\{\underline{x_1x_2},\;  x_1\underline{x_3x_4}, \; x_2\underline{x_3x_4}x_5\}$,\\
    $\{\underline{x_1x_2}, \; x_2\underline{x_3x_4},\;  x_1\underline{x_3x_4}x_5\}$,
    $\{x_1x_3, \; x_1x_2x_4,\;  x_3x_2x_4x_5\}$,
    $\{x_1x_3, \; x_3x_2x_4,\;  x_1x_2x_4x_5\}$, ... \;
     there are $(m_1m_2m_3m_4)^2m_5$ ideals each.
     \\
     Then in total 
    $ 12\; m_1m_2m_3m_4m_5\sum_{i<j<k<\ell}(m_im_jm_km_\ell)^2$ such ideals.
 
    \item \verb|[l^2*m^3,  l^4*m,  l^4*m]|  %count := 10],
    %If $(a,b,c)=(2,4,4)$ 
    e.g. for the footprint
    $\{x_1x_2,\;  x_2\underline{x_3x_4x_5}, \; x_1\underline{x_3x_4x_5}\}$,\\
    there are $(m_1\cdots m_5)^2$ ideals. \;
     Then in total 
    $ 10\; (m_1m_2m_3m_4m_5)^2$ such ideals.
    
    \item \verb|[l^3*m^2,  l^3*m^2,  l^3*m^2]|  %count := 30]
    %If $(a,b,c)=(3,3,3)$ 
e.g. for the 6 footprints
$\{x_3\underline{x_1x_4},\;x_2\underline{x_1x_4}, \; \underline{x_2x_3}x_5\}$, \;
  $\{x_4\underline{x_1x_3},\;x_2\underline{x_1x_3}, \; \underline{x_2x_4}x_5\}$, 
     there are $(m_1m_2m_3m_4)^2m_5$ ideals each.
     \\
     Then in total 
    $ 6\; m_1m_2m_3m_4m_5\sum_{i<j<k<\ell}(m_im_jm_km_\ell)^2$ such ideals.

    \item \verb|[l^3*m^2,  l^3*m^2,  l^4*m]|  %count := 15],
    %If $(a,b,c)=(3,3,4)$
    e.g. for 
    $\{x_1\underline{x_4x_5}, \; x_1\underline{x_2x_3}, \; \underline{x_2x_3}\, \;\underline{x_4x_5}\}$
    there are 
    $(m_1\cdots m_5)^2$ ideals.
         Then in total 
         $15\; (m_1m_2m_3m_4m_5)^2$

\end{itemize}

%% --- 4 gens -----------------------
And finally, the dominant ideals with 4 generators:
 \begin{itemize}
     \item   \verb|[l^3*m^2,  l^3*m^2,  l^3*m^2,  l^3*m^2]|  %count := 5],
     %For $(3,3,3,3)$ 
e.g. for the footprint \\
  $\{x_2x_3x_4,  x_1x_3x_4,  x_1x_2x_4,  x_1x_2x_3\}$
there are $(m_1m_2m_3m_4)^3$ ideals each.
\;
     Then in total 
    $\sum_{i<j<k<\ell}(m_im_jm_km_\ell)^3$ such ideals.

     \item   \verb|[l^3*m^2,  l^3*m^2,  l^3*m^2,  l^4*m]|  %count := 20],
     %For $(3,3,3,4)$ 
e.g. for the 4 footprints \\
$\{x_2x_3x_4, \; x_1x_2x_4, \; x_1x_2x_3, \; x_1x_3x_4 \underline{x_5}\}$\;
$\{x_2x_3x_4, \; x_1x_2x_4,\;  x_1x_2x_3\underline{x_5}, \;  x_1x_3x_4 \}$\; ...,\\
there are $(m_1m_2m_3m_4)^3m_5$ ideals each.
\\
     Then in total 
    $ 4\; m_1m_2m_3m_4m_5\sum_{i<j<k<\ell}(m_im_jm_km_\ell)^2$ such ideals.
 
     \item   \verb|[l^3*m^2,  l^3*m^2,  l^4*m,  l^4*m]|  %count := 30],
     %For $ (3,3,4,4)$ 
e.g. for the 6 footprints\\
$\{x_2x_3x_4,  x_1x_2x_3,  x_1x_3x_4\underline{x_5},  x_1x_2x_4\underline{x_5}\}$,\;
$\{x_2x_3x_4,  x_1x_2x_3\underline{x_5},  x_1x_3x_4,  x_1x_2x_4\underline{x_5}\}$,\; ...\\
there are $(m_1m_2m_3m_4)^3m_5^2$ ideals each.
\\
     Then in total 
    $ 6\; (m_1m_2m_3m_4m_5)^2\sum_{i<j<k<\ell}m_im_jm_km_\ell$ such ideals.
 
     \item   \verb|[l^3*m^2,  l^4*m,  l^4*m,  l^4*m]|  %count := 10],
e.g. for the footprint\\
$\{\underline{x_1x_2x_3}, \; x_2x_3x_4x_5, \; x_1x_3x_4x_5, \; x_1x_2x_4x_5\}$
there are $(m_1m_2m_3m_4m_5)^3$ ideals.
\\ Then in total 
    $ 10\; (m_1m_2m_3m_4m_5)^3$ such ideals. \qedhere
     %For $(3,4,4,4)$ there are $10 (m_1\cdots m_5)^3$ dominant ideals. \qedhere
 \end{itemize}
 \end{proof}

\subsection{Explicitly computing all dominant ideals with a given lcm}\label{sec:compute-all-dominant-lcm}

In this section
we provide a  description of our CoCoA implementation. 
At the following link \url{https://sites.google.com/view/bigatti-data} the interested reader will find our CoCoA-5 package,  \sloppy\verb|DominantIdeals.cpkg5 |\allowbreak (which is part of the forth-coming release CoCoA-5.4.2), together with the files containing all examples.

In particular, we start by describing an algorithmic improvement that is fundamental for listing all dominant monomial ideals with a given lcm, especially as the number of variables grows.

In $\KK[x,y,z]$ the implementation first makes a list of candidate dominant ideals with a given lcm, making all ideals of the form
$(
x^{m_1}y^{i_1}z^{i_2},
x^{j_1}y^{m_2}z^{j_2},
x^{k_1}y^{k_2}z^{m_3}
)$.
Please note that ideal 
$\thelcm=x^{m_1}y^{m_2}z^{m_3}$
is in all the 3 lists \texttt{Lx, Ly, Lz}
and the ideal $\langle \thelcm, \thelcm, \thelcm\rangle$ has only one minimal generator.
Moreover the ideals
$(
x^{m_1}y^{m_2}z^0,
x^{m_1}y^{m_2}z^{j_2},
x^{k_1}y^{k_2}z^{m_3}
)$
have two minimal generators and are produced in multiple copies.
Then makes the ideal minimally generated, and select those with lcm = $\thelcm$ and dominant.
Finally it removes the duplicates.

(for clarity, we just show the most significant lines):
\begin{lstlisting}[language=cocoa]
define AllDominantIdeals3(R, m) -- m is the list of exponents of the LCM
  (...)
  Lx := [ x^m[1]*y^i[1]*z^i[2] | i in (0..m[2])><(0..m[3])];
  Ly := [ x^j[1]*y^m[2]*z^j[2] | j in (0..m[1])><(0..m[3])];
  Lz := [ x^k[1]*y^k[2]*z^m[3] | k in (0..m[1])><(0..m[2])];
  candidates := [IdealOfMinGens(ideal(L)) | L in CartesianProduct(Lx,Ly,Lz) ];
  (...)
  DominantIdeals := [ I In candidates
		     | lcm(gens(I))=LCMTerm and IsDominant(I) ];
  (...)
  DominantIdealsNoDup := MakeSetIdealList(DominantIdeals);
  (...)
  return DominantIdealsNoDup;    
\end{lstlisting}

Although set-based removal of duplicates may seem straightforward, the procedure is computationally expensive when the list is very long; in fact, the equality test for monomial ideals is easy, but non-trivial. 

For example, with \verb|LCMExps := [5,2,4];| we get 8100 candidates.
Removing the duplicates has quadratic complexity, making $l\cdot(l{-}1)$ equality tests between pairs from~$l$ monomial ideals, and it strongly affects this apparently easy procedure.  So we need be careful and avoid pointlessly long list.

Therefore, with 4 variables $x,y,z,w$, instead of taking all candidates, we take some intermediate steps.
We start by making the list \verb|candidate2| of all monomials ideals with 1 or 2 generators taken from \verb|Lx, Ly|, and remove the duplicates.
Then, in \verb|candidate3| we add to each of them a generator from \verb|Lz|, and remove the duplicates,
then we select those which are dominant (non necessarily with the given lcm), and put them in \verb|dominant3|.
Finally we make the list \verb|candidate4| by adding to each ideal in \verb|dominant3| each monomial in \verb|Lw|, we select those which are dominant and with the given lcm, and then remove the duplicates.

\begin{remark}
    The function \texttt{MakeSetIdealList} takes as input a (long) list of monomial ideals and returns a deduplicated list. 
To overcome this challenge, we developed a recursive divide-and-conquer strategy. 

In our implementation, the monomial ideals in the input list are minimally generated, the routine internally chooses a monomial $t$, and then splits the input list into two disjoint sublists, \verb|with_t| and \verb|without_t| based on whether or not $t$ is a minimal generator, hence a duplicate cannot appear across the sublists. After recursively deduplicating each sublist, the results are merged to form the final output. 
\end{remark}

\subsection{Computing the formula for counting all dominant ideals with given lcm}\label{sec:computing-formula-all-dominant-lcm}

The function \texttt{DominantIdealHistogramFootprint} is designed to compute and count the footprints of each ideal by identifying the variables for which the corresponding exponents are not maximal. 
When the input exponent vector consists solely of 1's (that is, when all monomial generators are squarefree), the sum of the monomials labeled by the resulting footprints exactly coincides with the counting formula for dominant ideals as 
described in implemented in the function \texttt{DominantFormula4} (with the substitution $m \mapsto x$), as stated in Theorem~\ref{thm:dominant-lcm-3}, \ref{thm:dominant-lcm-4}, \ref{thm:dominant-lcm-5}. This result not only confirms the theoretical expectations but also provides a computational tool to produce the formula for $n\ge5$.

%++++++++++++++++++++++++++++++++++

\subsection{Open problems about monomial ideals with a given lcm}\label{sec:lcm-ideals}

So far we computed all dominant monomial ideals with a given lcm up to 5 variables and we understood the idea behind such a computation.  However, the problem becomes quickly intractable as soon as the number of variable increases as we explain more in details in the following.

\begin{remark}%
 To understand distributions of dominant ideals we first need to understand how to generate or enumerate all monomials with a given lcm. One potential strategy is to fix the number of minimal generators $g\leq n$, and consider all sets of $g$ monomials whose least common multiple is $\thelcm$. For each such candidate set $S$, one can check whether $S$ forms the minimal generating set of a monomial ideal and whether it satisfies the dominance condition. However, this approach quickly becomes computationally intractable. Even for small values, such as $n=3$ and $\thelcm=x^4y^4z^4$, the number of possible generating sets whose lcm is $\thelcm$ is large. For instance, to get this number we would include many ideals of the form $(x^4,y^4z^4,t)$, for a monomial $t\in \KK[x,y,z]_{\leq 12}$. 

 This suggests that while fixing the lcm imposes a natural constraint, the space of monomial ideals having the same lcm remains too rich and complex. A combinatorial description or enumeration of dominant ideals within this constraint is thus both challenging and interesting, and likely requires additional structural constraints or combinatorial tools to be analyzed effectively. We leave a detailed investigation of such enumerative aspects for future research.
\end{remark}

A similar analysis applies if one focuses on understanding the same problem from a probabilistic point of view.

\begin{remark} %
 Instead of exhaustively enumerating all dominant ideals with a fixed least common multiple $\thelcm$, one may ask a more probabilistic question, wondering what is the likelihood that a randomly chosen monomial ideal with lcm $\thelcm$ is dominant. 
 
The key challenge in this direction lies in defining a meaningful and unbiased notion of randomness under the constraint that the minimal generators must have lcm equal to $\thelcm$. Sampling uniformly from the space of all such generating sets is nontrivial, since the space itself combinatorially structured.  We leave a careful formulation of this random model, and a probabilistic analysis of dominant ideals within it, to future work.  
\end{remark}

\section{Distributions of dominant ideals}\label{sec:random sampling}

We now seek a systematic way to understand how often  a monomial ideal may be dominant ``if one takes an ideal at random".  
To this end, we use  formal probabilistic models\footnote{To offer guidance to our readers not familiar with  randomness and  distributions, we recall that during the WICA II meeting in Trento, this collaboration began  with a simple question: what \emph{is} a random monomial ideal? By independently solving the challenge of computing 100 random monomial ideals in $3$ variables and computing the number of minimal generators and the Krull dimension---which the reader is encouraged to try themselves!---the authors of this paper jointly discovered that the \emph{model} by which one randomly generates ideals governs much of the resulting outcome.}  for randomly generated algebras and study the induced distributions of the dominance property. 
A natural first step is to use basic models for generating random monomial ideals, introduced in \cite{dLPSSW}. These models mirror and extend those from random graphs and simplicial complexes in the literature, thereby providing a connection to combinatorics. Let us  recall basic notation for probabilistic models in the algebraic setting of monomial ideals in $S=\mathbb K[x_1,\dots,x_n]$.

Fix an integer $D$ and a parameter $p=p(n,D)$, $0\leq p\leq 1$. Construct a random set of monomial ideal generators $B\subset S$ by including, independently, with probability $p$ each non-constant monomial of total degree at most $D$ in $n$ variables. The resulting random monomial ideal is simply $I=\ideal{B}$, and if $B=\emptyset$, then we let $I=\ideal{0}$. 
The notation for this data-generating process uses $\mathcal B(n,D,p)$ to denote the resulting distribution on the sets of monomials. This distribution on monomial sets induces a distribution on the set of ideals. This is called the \er-type distribution on monomial ideals and is denoted by $\rmidist(n,D,p)$. 

When an ideal $I$ is generated using  a random process just described, we say the ideal $I$ is an observation of the random variable drawn from the \er-type distribution. We denote this using standard notation for distributions and random variables: $$\rmi\sim\rmidist(n,D,p).$$ We may also refer to  $\rmi$ as an \er\ random monomial ideal. 

It is important to note that the set $\rmidist(n,D,p)$ is a parametric family of probability distributions, one for each set of values of the parameters $n$, $D$, and $p$.

\medskip

A series of natural questions arise next. For example: suppose $\rmi$ is a an \er\ random monomial ideal. What is the probability that $\rmi$ is dominant? Conditional on the fact that it is dominant, what is the probability that it is $0$-dimensional? Conditional on the fact that it is $0$-dimensional, what is the probability that it is dominant? There is an infinite list of questions one can ask, of increasing complexity. 
As to why one asks questions about a random monomial ideal, consider the fact that enumerating all monomial ideals in a given number of variables is computationally intractable; thus we generate samples of them randomly and compute the frequencies of dominant ideals within each sample. 

 \medskip 
 In the remainder of this section, we provide some statistics on dominant ideals that are generated from various versions of the \er\ type model. 
The experiments were performed using {\tt Macaulay2}, in particular the functions {\tt isDominant()} (Listing~\ref{M2dominantF}) and  in combinations of various functions from the package {\tt RandomMonomialIdeals.m2} \cite{RMIm2}.

To generate a sample of ideals from $\rmidist(n,D,p)$, for each fixed number of variables $n=3,\dots,6$, we  vary the maximal degree $D=2,\dots,15$. 
For a fixed $D$ in this range, we consider several values for the probability parameter $p$. The natural choice is to let $p\in\left\{ \frac{1}i, i=1,\dots,9 \right\}$, to get a global view of the full range of probabilistic outcomes. However, we also add additional very small values of $p$; these are selected to aim for the sample of ideals that contains ideals of each possible Krull dimension. Specifically,  \cite[Fig.1 and Cor.1.2]{dLPSSW} state that the various powers of $1/D$ stratify the set of random monomial ideals by their Krull dimension, asymptotically almost surely. While we are not in the asymptotic regime, as all values of $n$ and $d$ are small in these simulations, these thresholds are  the best guess for obtaining an expected Krull dimension.  
The values of the probability parameter $p$ in the simulations are thus computed as the midpoints of the dimension threshold intervals: 
$$ 
\frac{D^{-i}-D^{-(i-1)}}{2}, 
 \hbox{ for } i=2,\dots,n.
$$

\begin{figure}[h!]

\hbox to\linewidth{%
    \hfil%
    \includegraphics[scale=0.4]{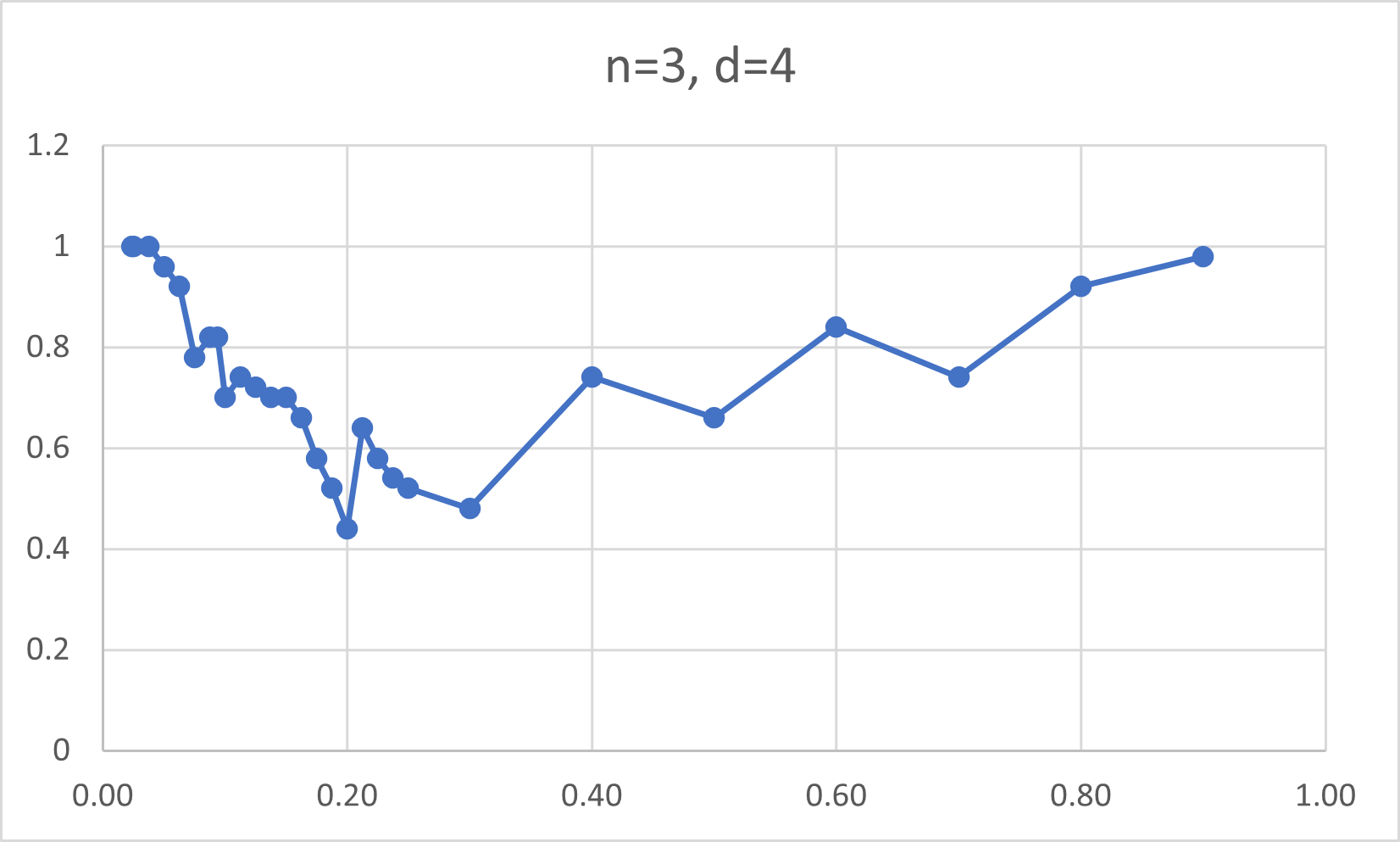}%
    \hfil%
    \includegraphics[scale=0.4]{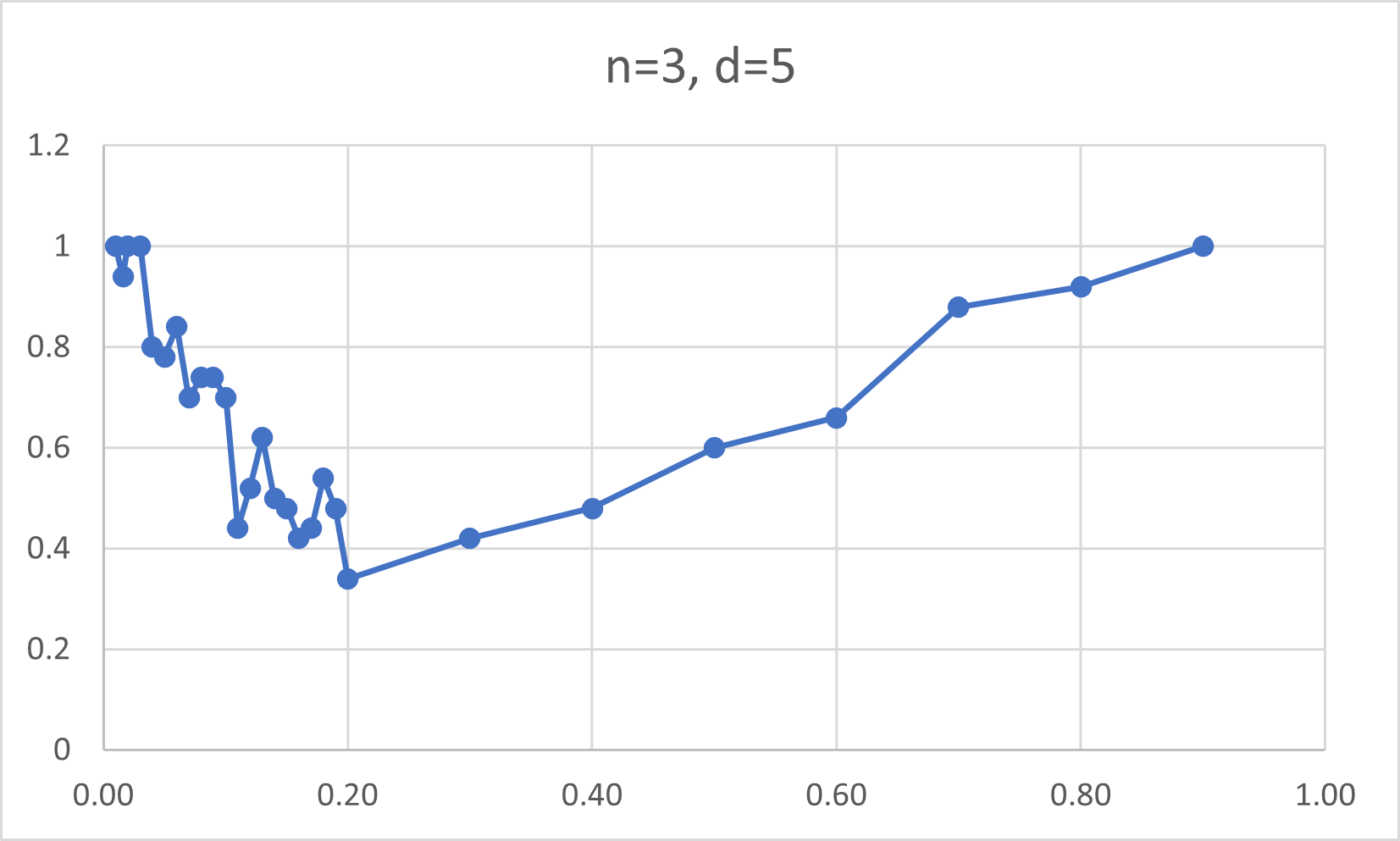}%
    \hfil%
    \includegraphics[scale=0.4]{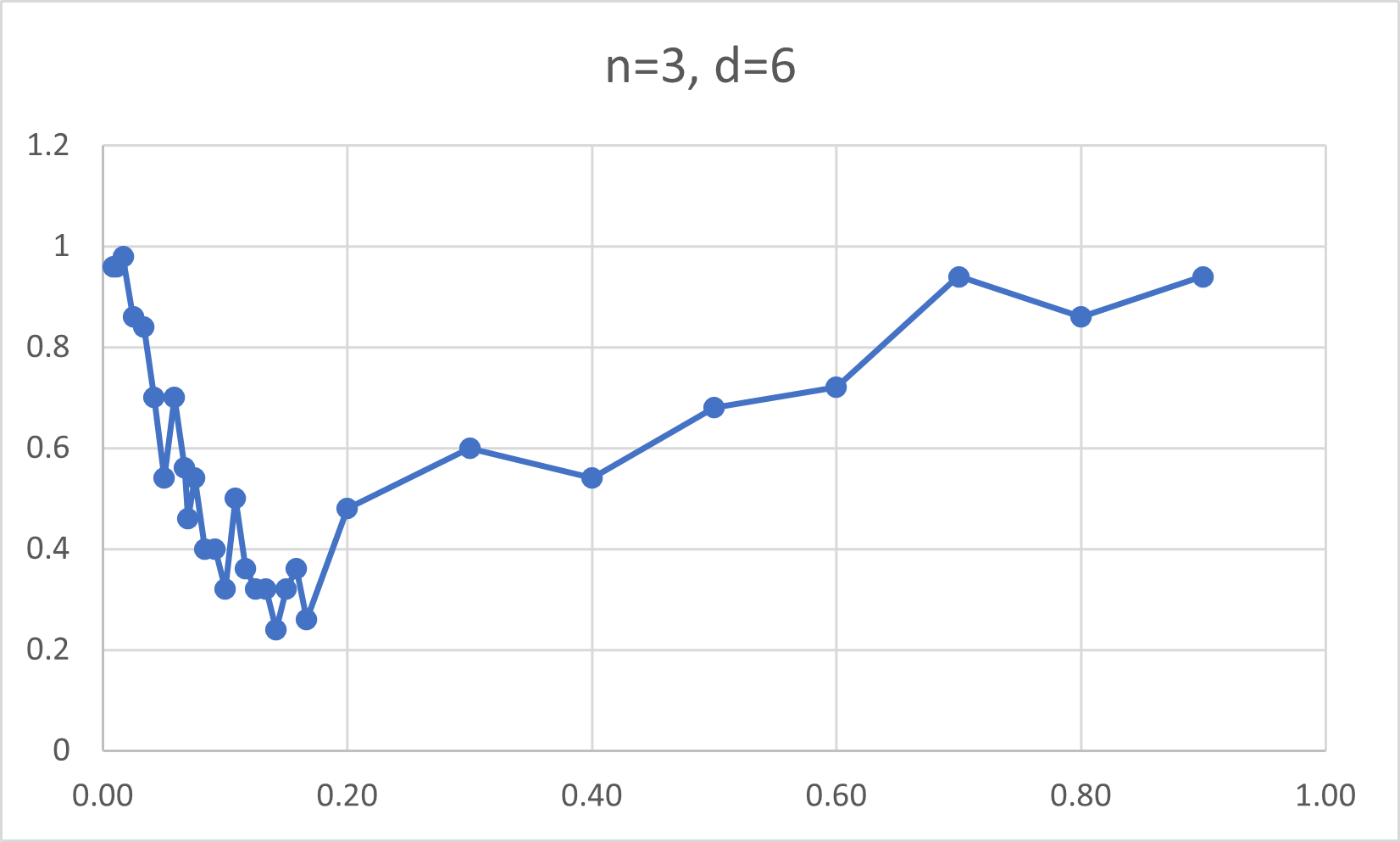}%
    \hfil%
}

\hbox to\linewidth{%
    \hfil%
    \includegraphics[scale=0.4]{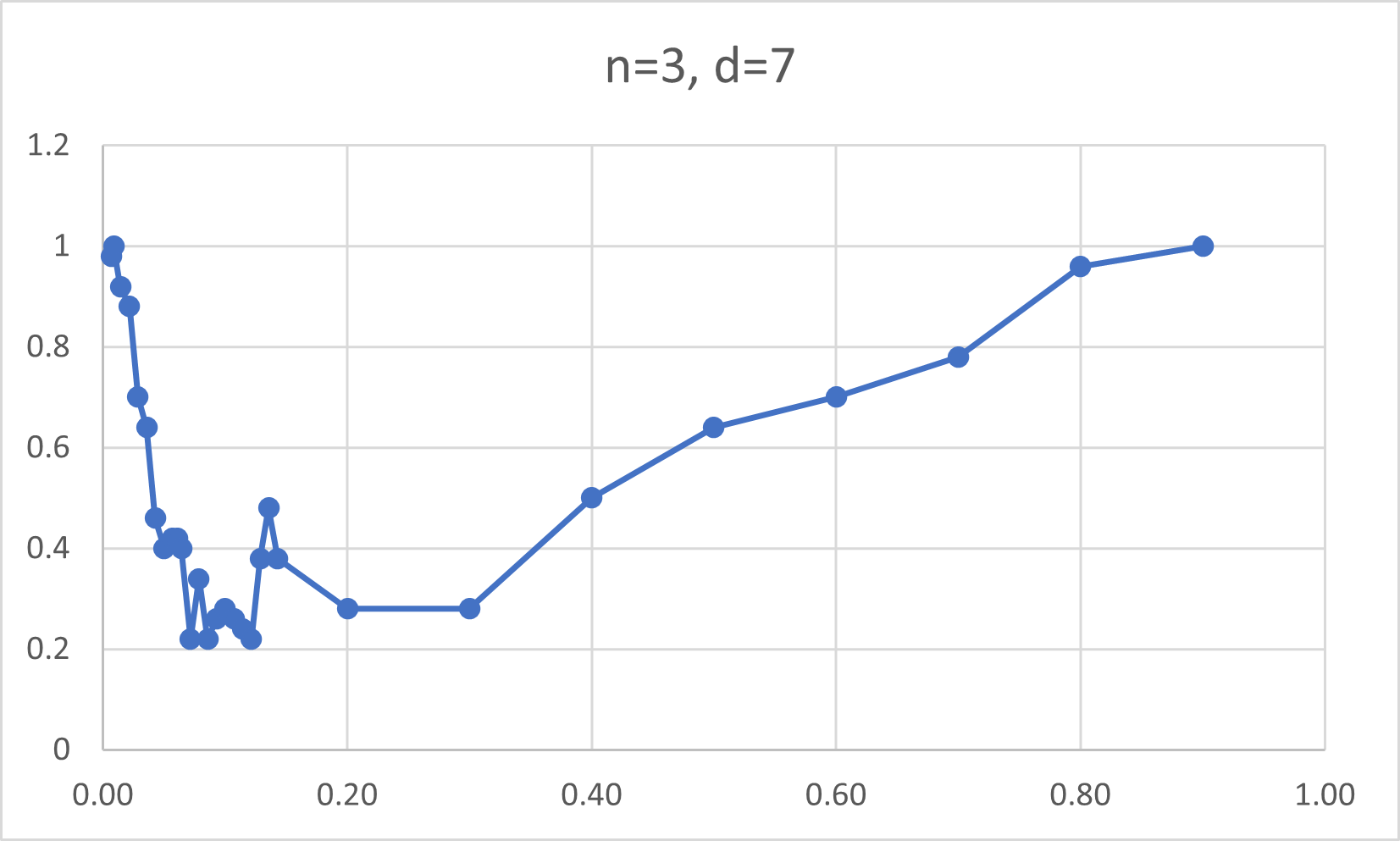}%
    \hfil%
    \includegraphics[scale=0.4]{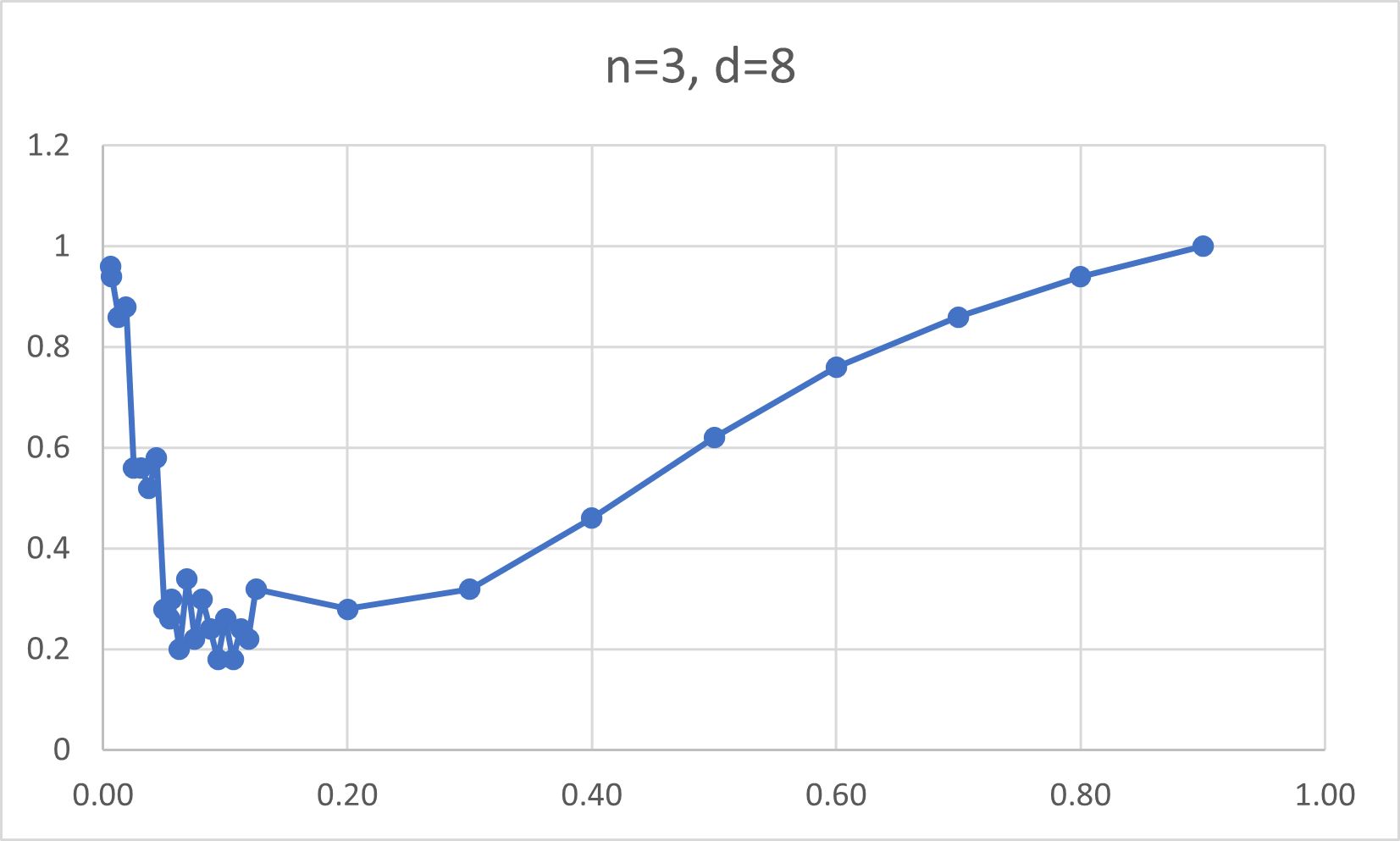}%
    \hfil%
    \includegraphics[scale=0.4]{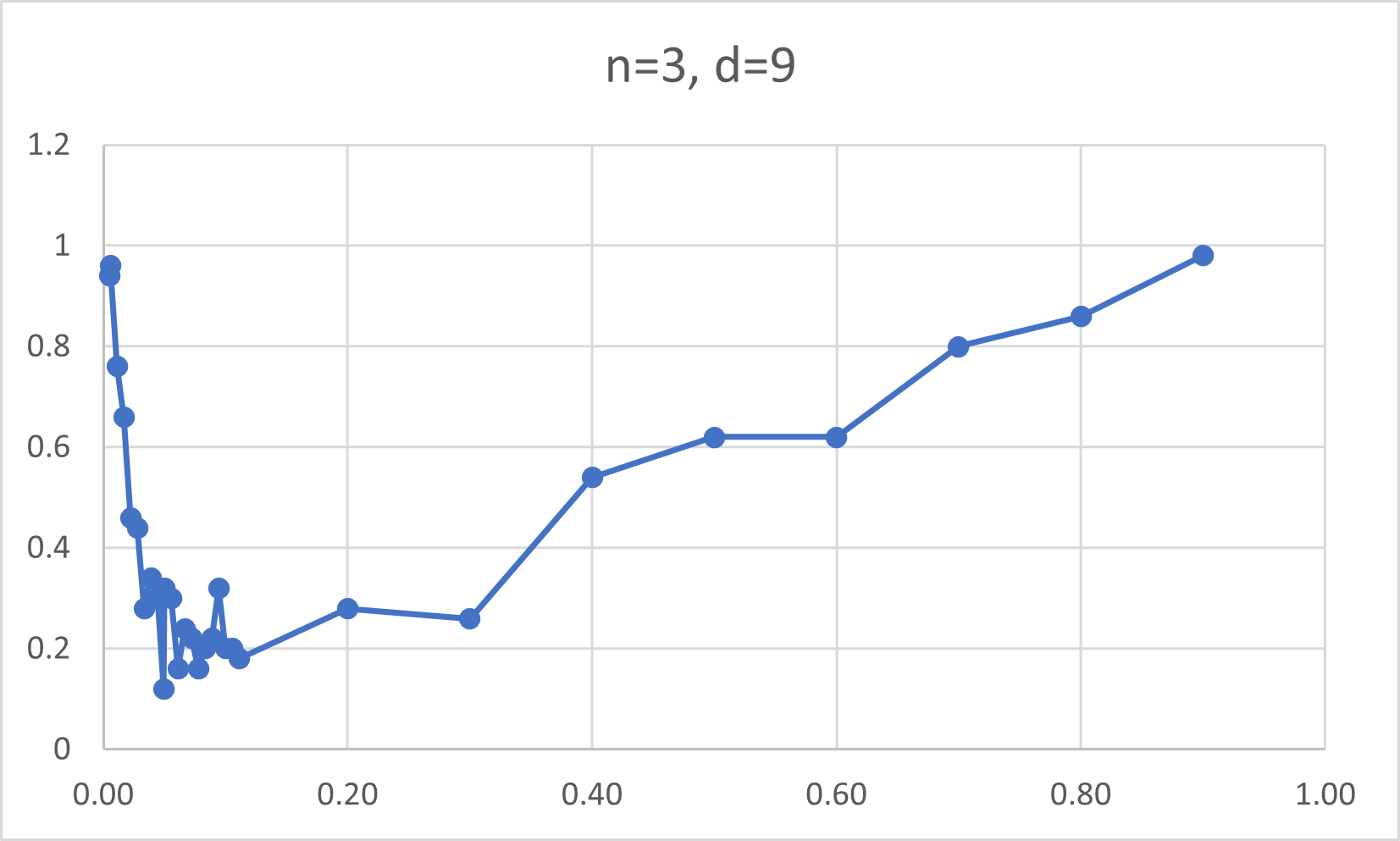}%
    \hfil%
}

\hbox to\linewidth{%
    \hfil%
    \includegraphics[scale=0.4]{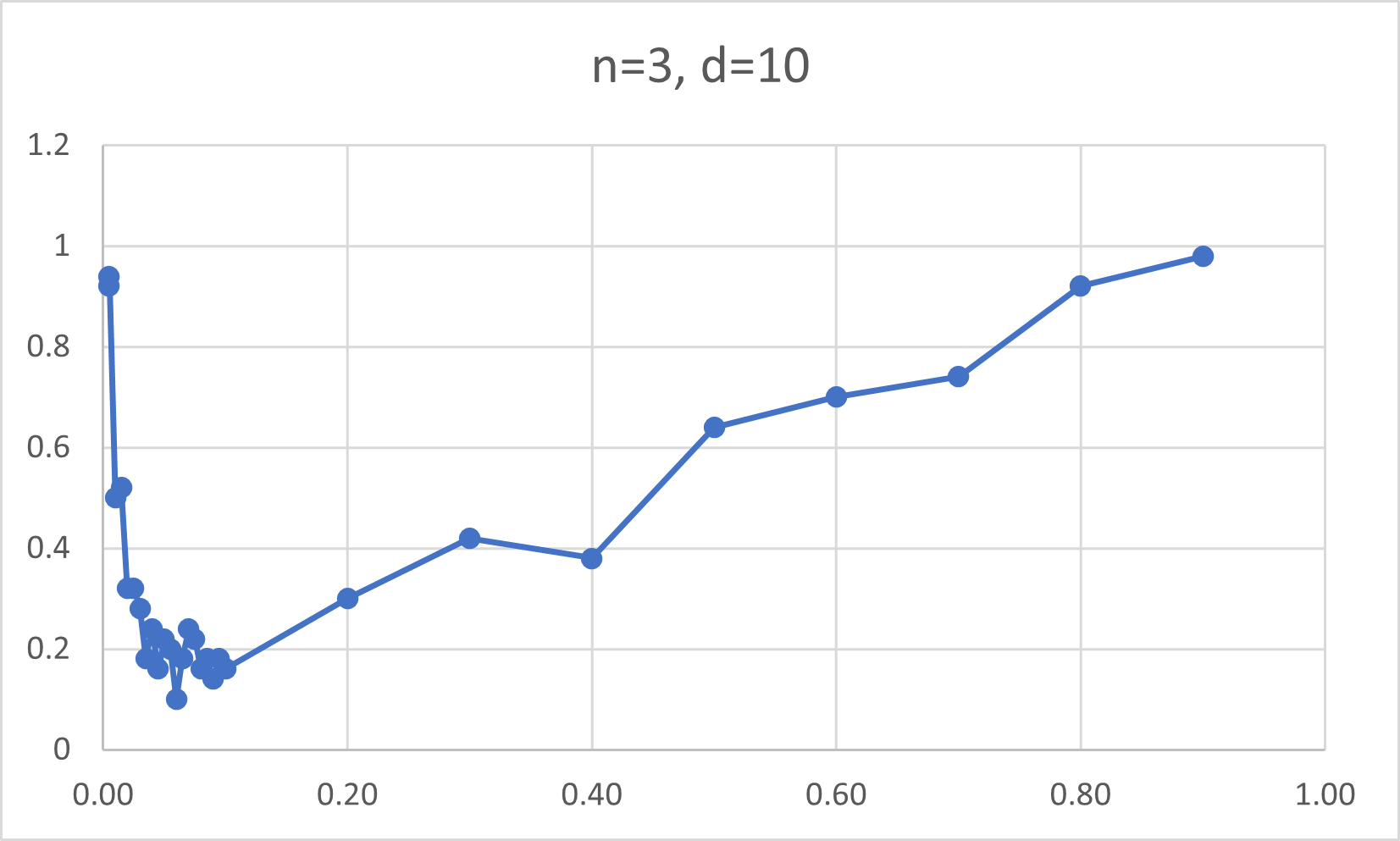}%
    \hfil%
    \includegraphics[scale=0.4]{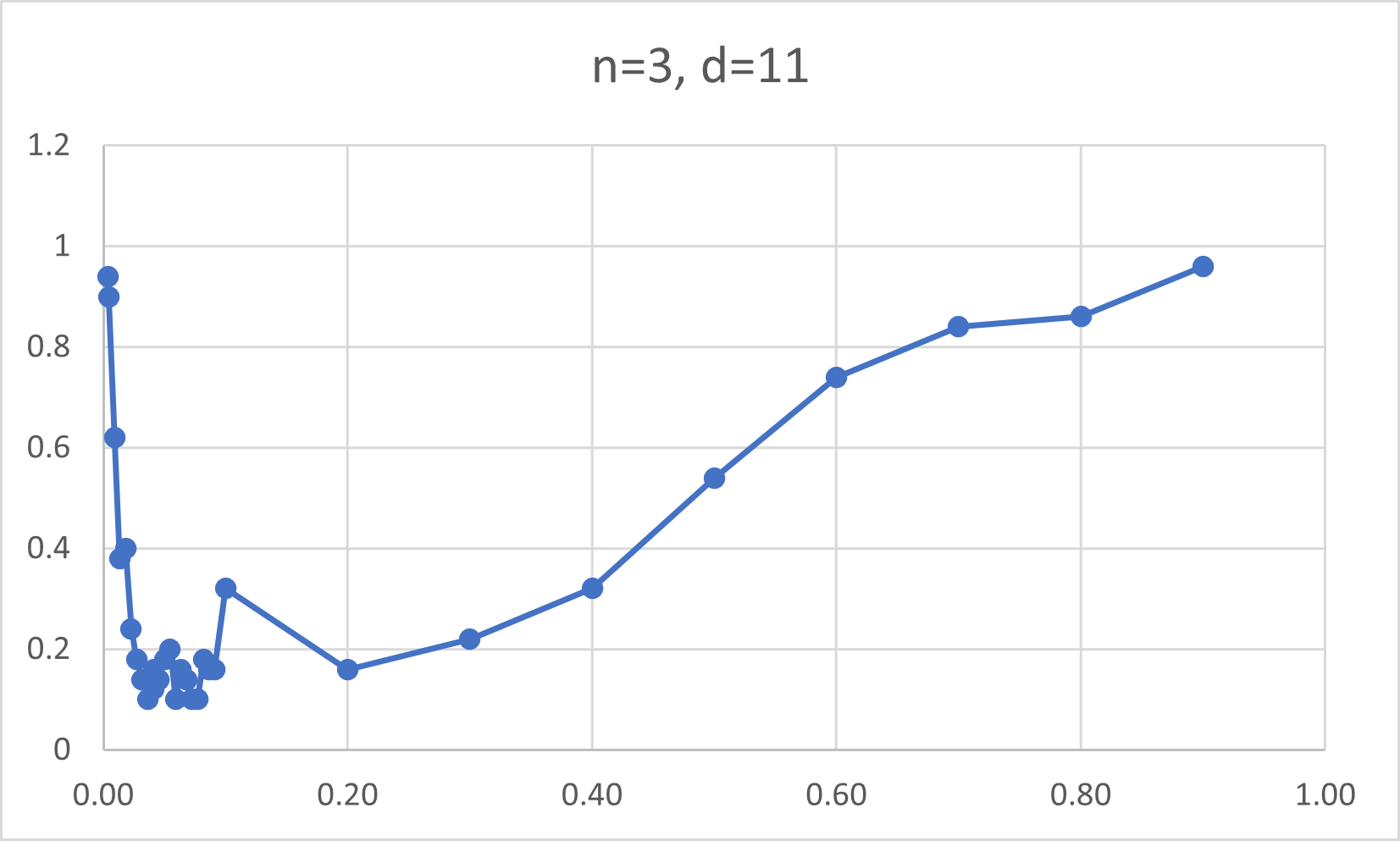}%
    \hfil%
    \includegraphics[scale=0.4]{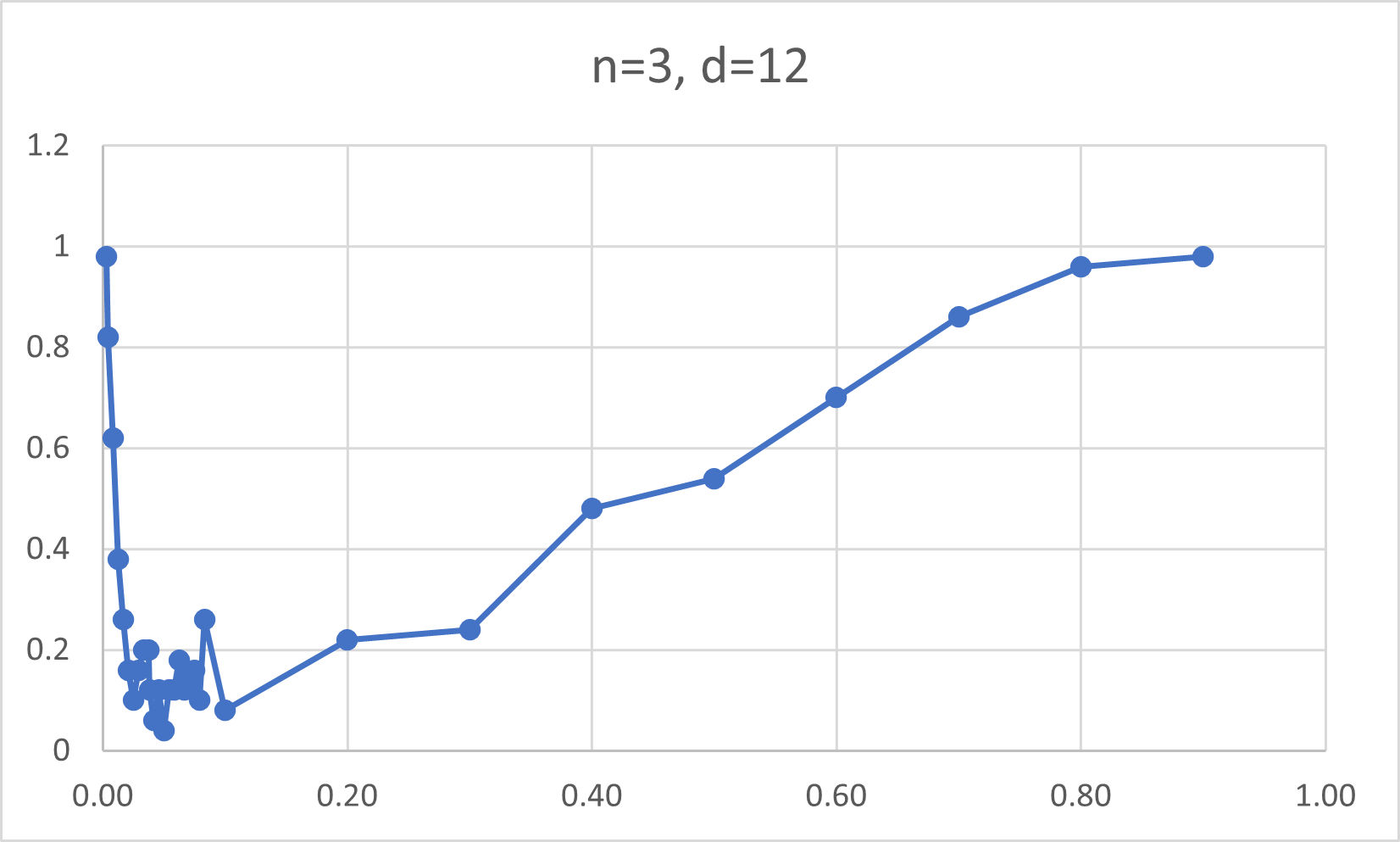}%
    \hfil%
}

\hbox to\linewidth{%
    \hfil%
    \includegraphics[scale=0.4]{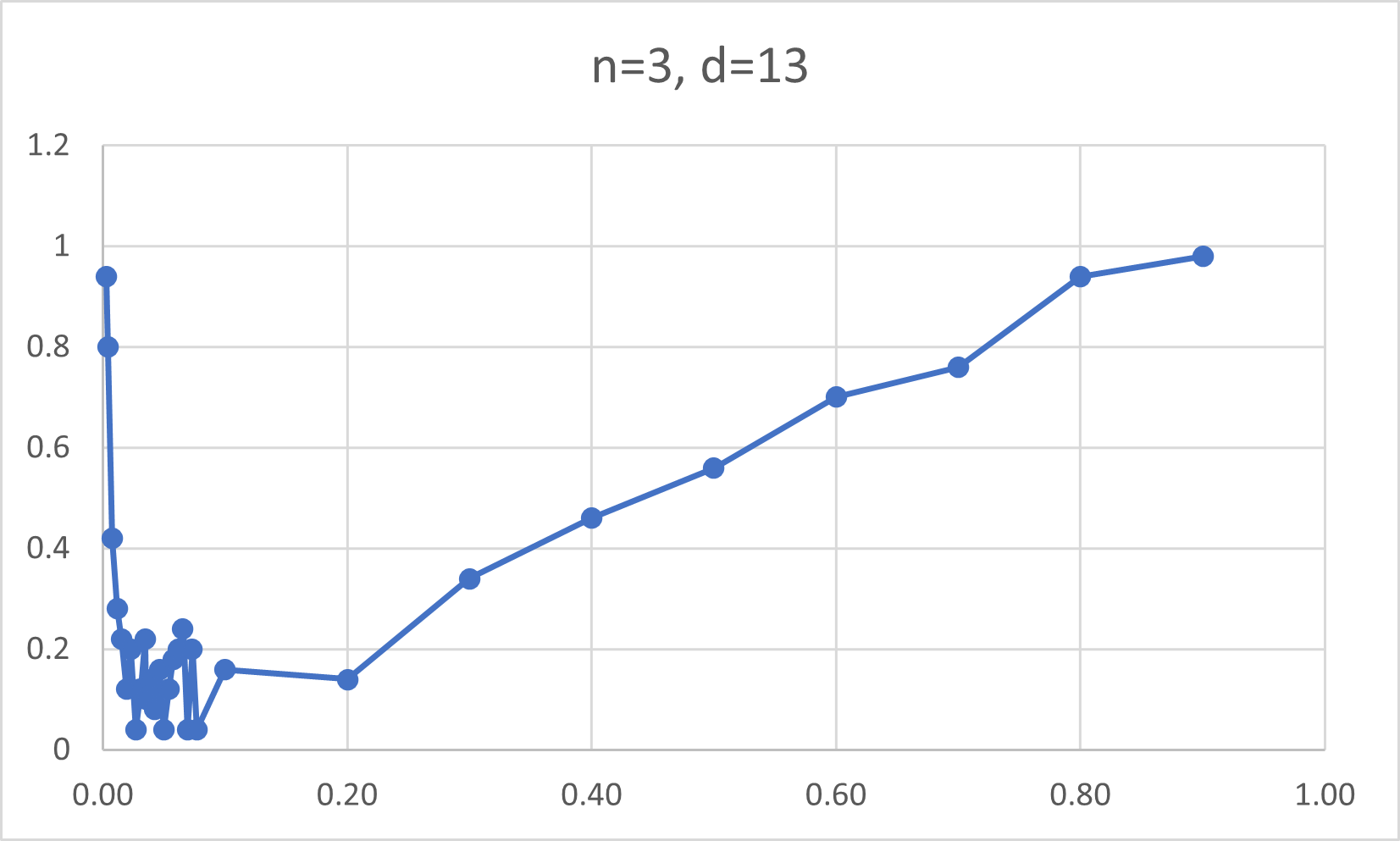}%
    \hfil%
    \includegraphics[scale=0.4]{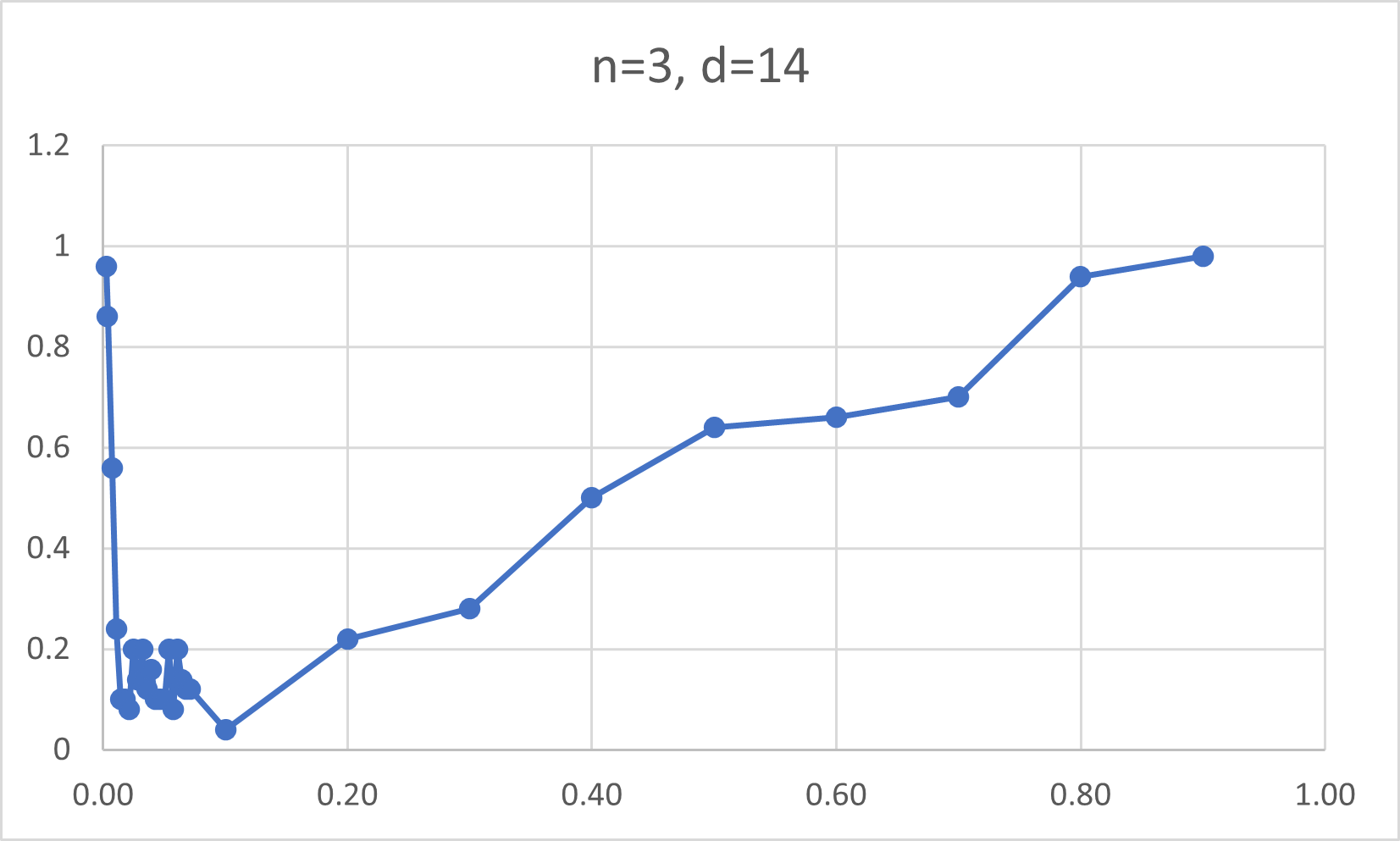}%
    \hfil%
    \includegraphics[scale=0.4]{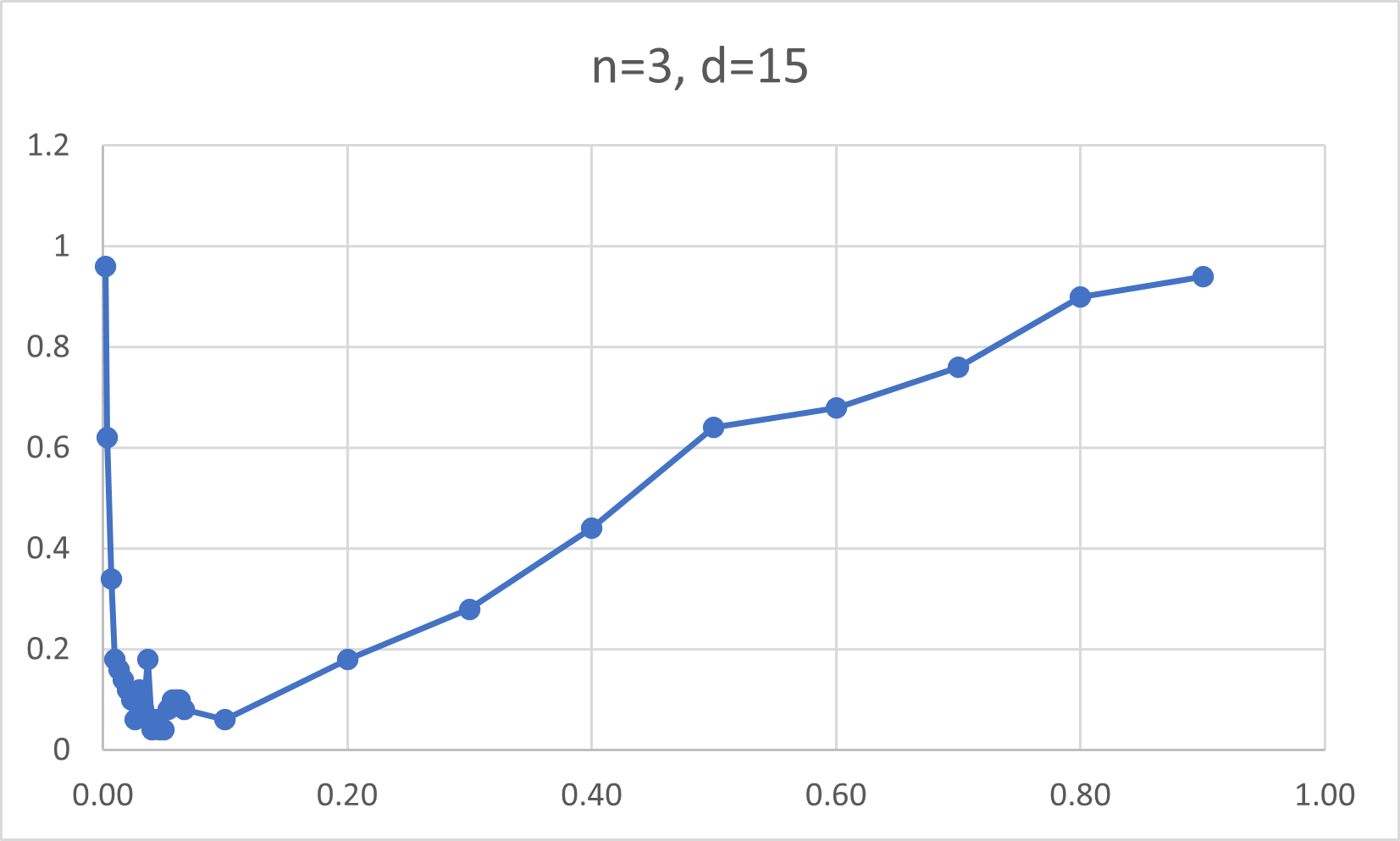}%
    \hfil%
}

  \caption{Frequency of dominant random monomial ideals in $\rmidist{I}(n,D,p)$ for $n=3$, $D=3,\dots,15$ and $p$ taking values in $\{ (D^{-i}-D^{i-1}) /2, \,i=2,3\}\cup \{1/2,1/3,\dots,1/9 \} $. Each data point in the figure represents a sample of size $50$ for a fixed value of $p$.} %We used a sample of size 50 for each $p$.}
    \label{fig:frequency of dominant RMIs}

\end{figure}

Figure~\ref{fig:frequency of dominant RMIs} shows the simulation results. 
\Cref{sec:code for RMI basic model} contains code that can be used to generate the data that produced these plots. 

%\pp{++explain the results and what we read from the results++}
From the plots in Figure~\ref{fig:frequency of dominant RMIs}, we see that all sampled ideals are dominant for $p\approx0$, meaning the higher powers of $1/D$; then the proportion of dominant ideals drops drastically, often to close to $0$. As we increase $p$, the proportion of ideals that are dominant grows closer to $100\%$ again. If $p>1/D$, then the random monomial ideal has Krull dimension zero, asymptotically. In the given samples,  it appears likely that the Krull dimension is indeed zero, and thus each ideal in the sample contains pure power monomials. Each of the pure powers  will then form a dominant subset of the generators. Since the probability of picking any given monomial is close to 1, the pure powers will be low degree, ensuring that additional generators selected in the random set $B$ are not minimal generators. 

The reader will have probably noticed one annoying feature of the simulations so far, namely, that we have no control over \emph{exaclty how many minimal generators} each random monomial ideal has. Thus, we consider a modified model next.

\subsection{Homogeneous random monomial ideals} %Graded model of the ER-type distribution}

The \emph{graded model} (see \cite[Section 5]{dLPSSW}) is a refined version of the $\rmidist(n,D,p)$, in which the probability parameter $p=(p_1,\dots,p_D)$ is a vector,  $0\leq p_i<1$.  This allows for finer control over the number of generators in \emph{each} degree in the set $\{1,\dots,D\}$. 
For all $i\leq D$, the degree-$i$ monomial is selected with probability $p_i$. The corresponding probability distribution is denoted by $\mathcal{I}_{\mathrm{Gr}}(n,D,p)$.

In particular, if only one entry of $p$ is nonzero, say the $j$-th one, then the model generates only \emph{homogenous} ideals, as the only  monomials that are selected are the ones of degree $j\leq d$, and they are chosen independently with probability $p_j$. In this scenario, every generator is minimal by construction.

For the graded model simulation, we vary the number of variables $n=3,\dots,6$ and maximum degree $D=2,\dots,12$. 
For every degree $d=2,\dots,D$, the probability vector is a $d$-tuple $(0,\dots, 0,\alpha)$ where 
$$
\alpha \in \left\{ \frac{d^\ell-d^{\ell-1}}{2},\, \ell=2,\dots,n \right\} \cup \left\{ \frac{1}{9},\frac{1}{8},\dots,1 \right\} \cup \left\{\frac{1}{20d},\dots,\frac{20}{20d}  \right\}.
$$
The code used to generate data in  Figure~\ref{fig:frequency of dominant homogeneous RMIs} can be found in Appendix~\ref{sec:code for RMI homogeneous model}. 
The simulations indicate that the proportion of dominant ideals drops, monotonically, from $1$ to $0$ as the probability of choosing homogeneous monomials grows from $0$ to $1$, for each fixed value of number of variables $n$ and generating degree $d$. 
The data suggests the following, which can likely be stated for general $n$ with a suitable formula for the probability threshold that depends on $n$ and $d$:
\begin{conjecture}
Let  $n=3$. As $d$ increases, when $p>0.1$, the  random homogeneous ideal is dominant with probability tending to $0$. 
\end{conjecture}

\begin{figure}[h!]

\hbox to\linewidth{%
    \hfil%
    \includegraphics[scale=0.25]{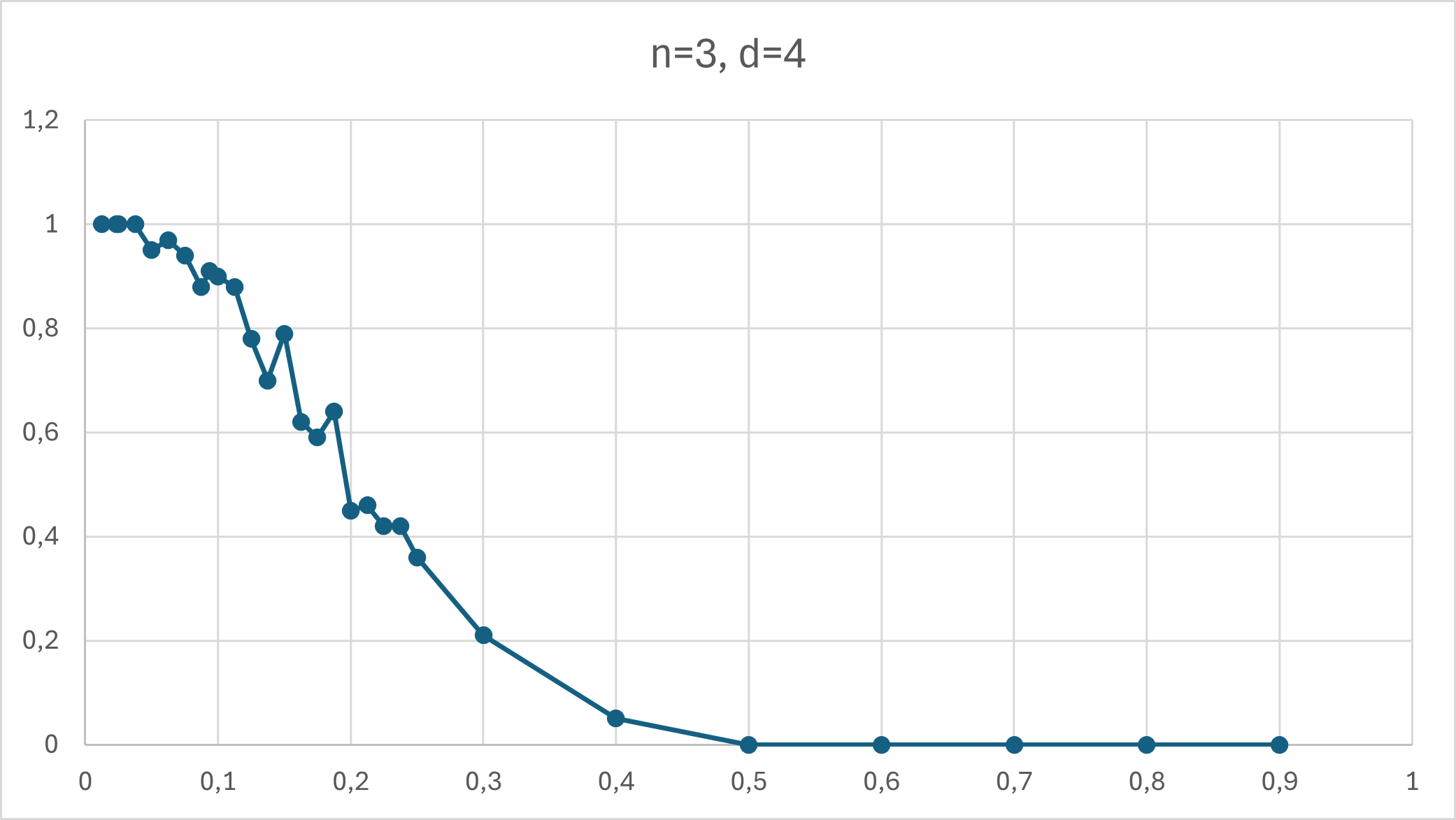}%
    \hfil%
    \includegraphics[scale=0.25]{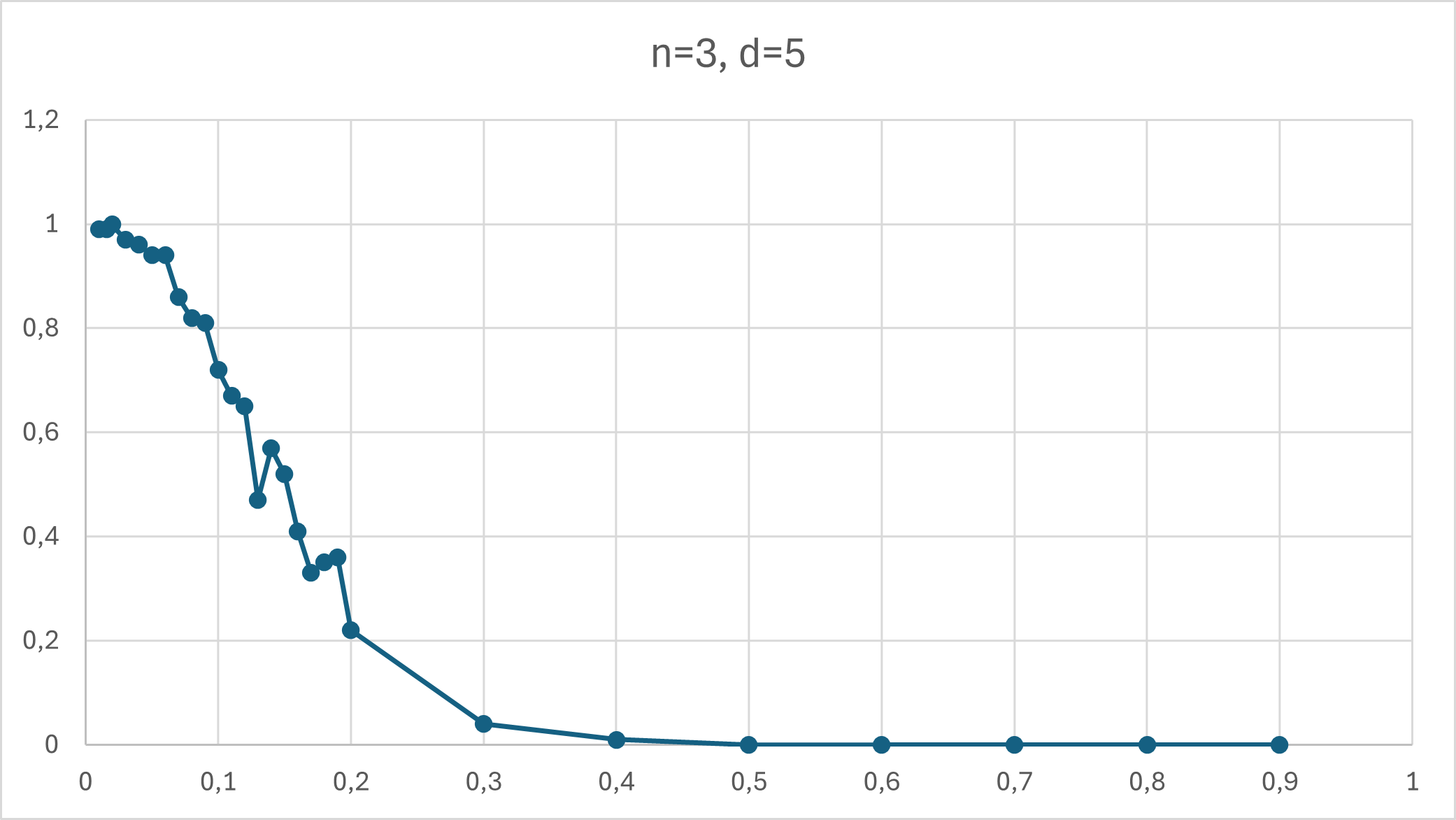}%
    \hfil%
    \includegraphics[scale=0.25]{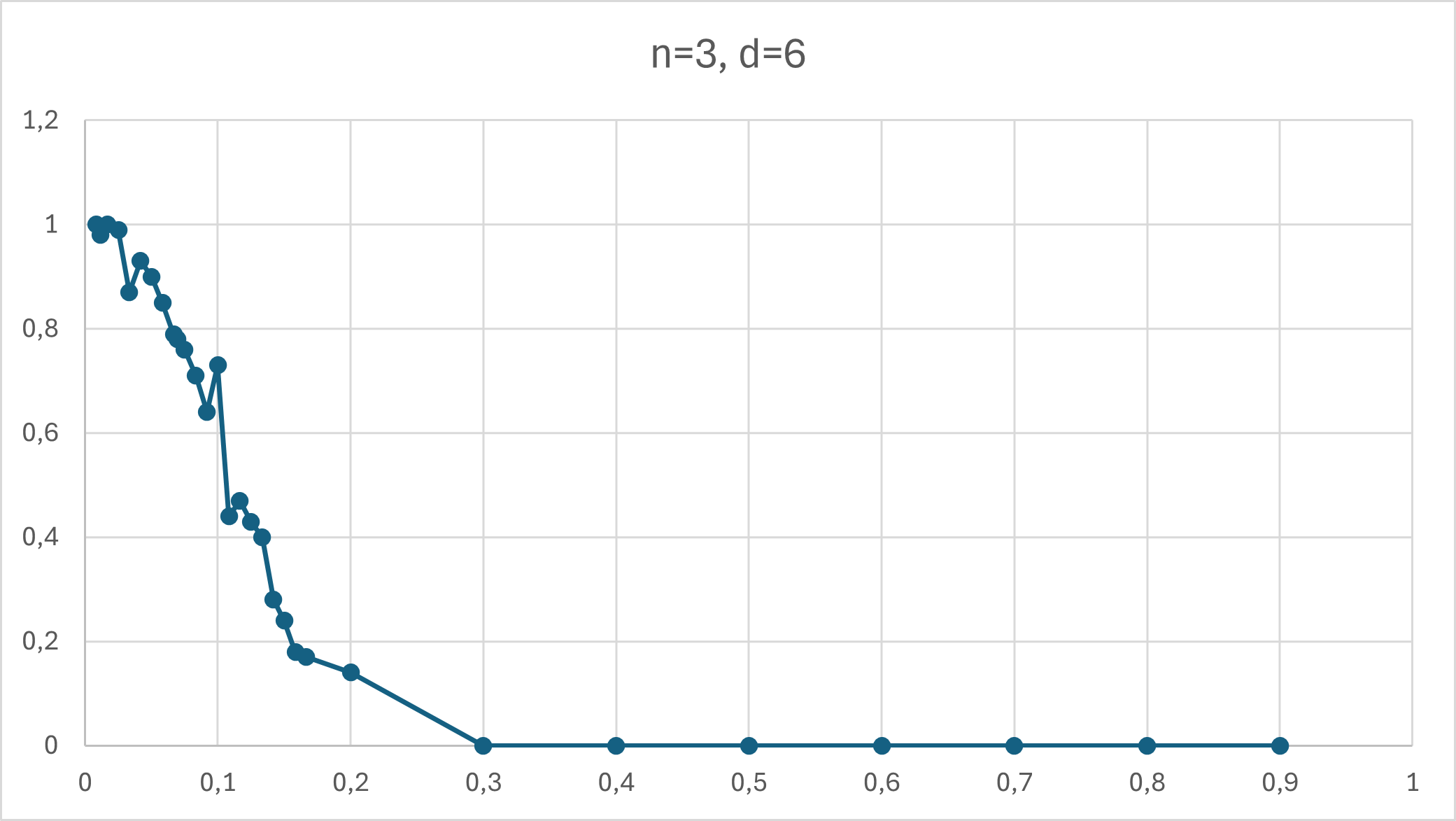}%
    \hfil%
}

\hbox to\linewidth{%
    \hfil%
    \includegraphics[scale=0.25]{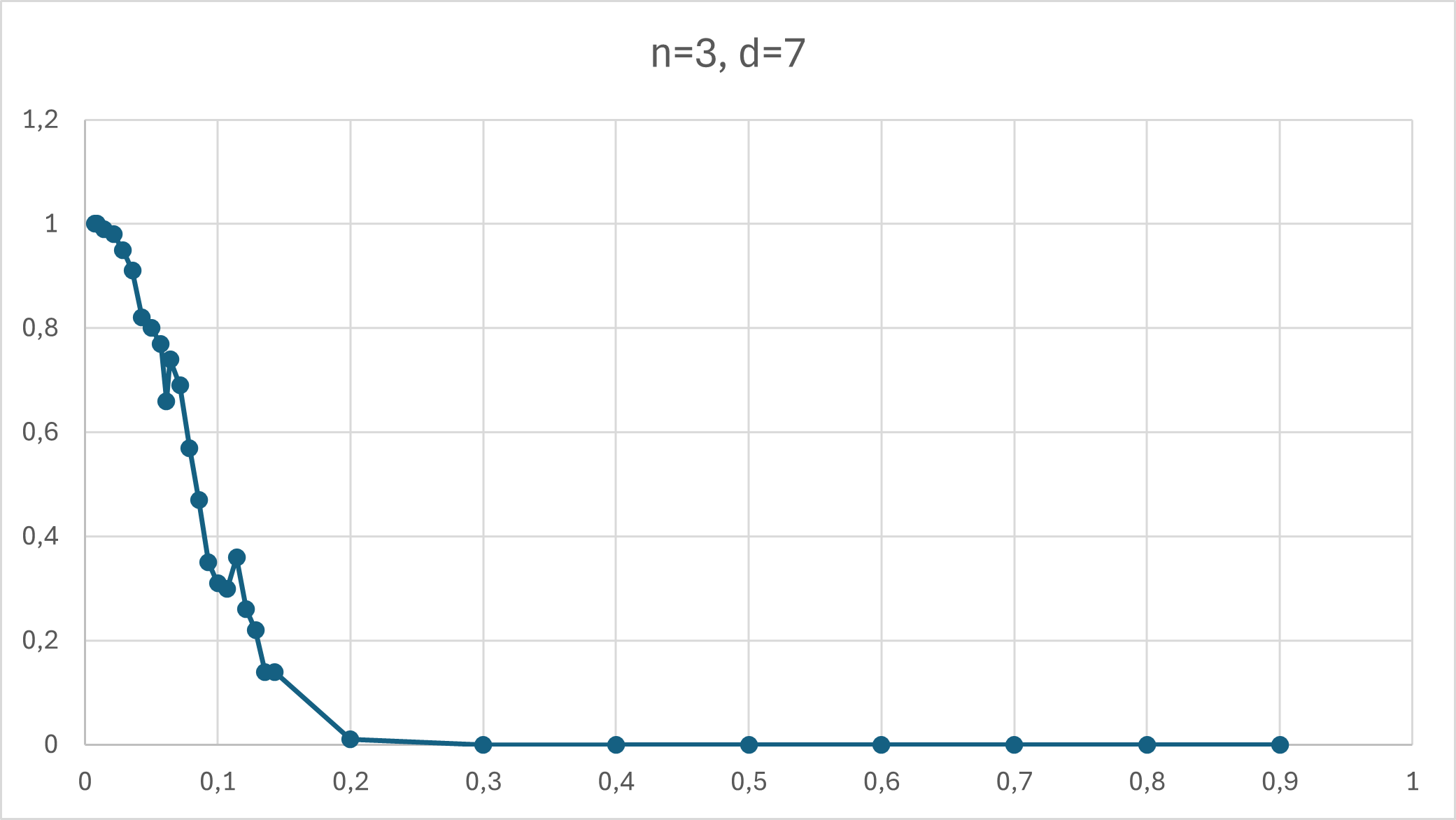}%
    \hfil%
    \includegraphics[scale=0.25]{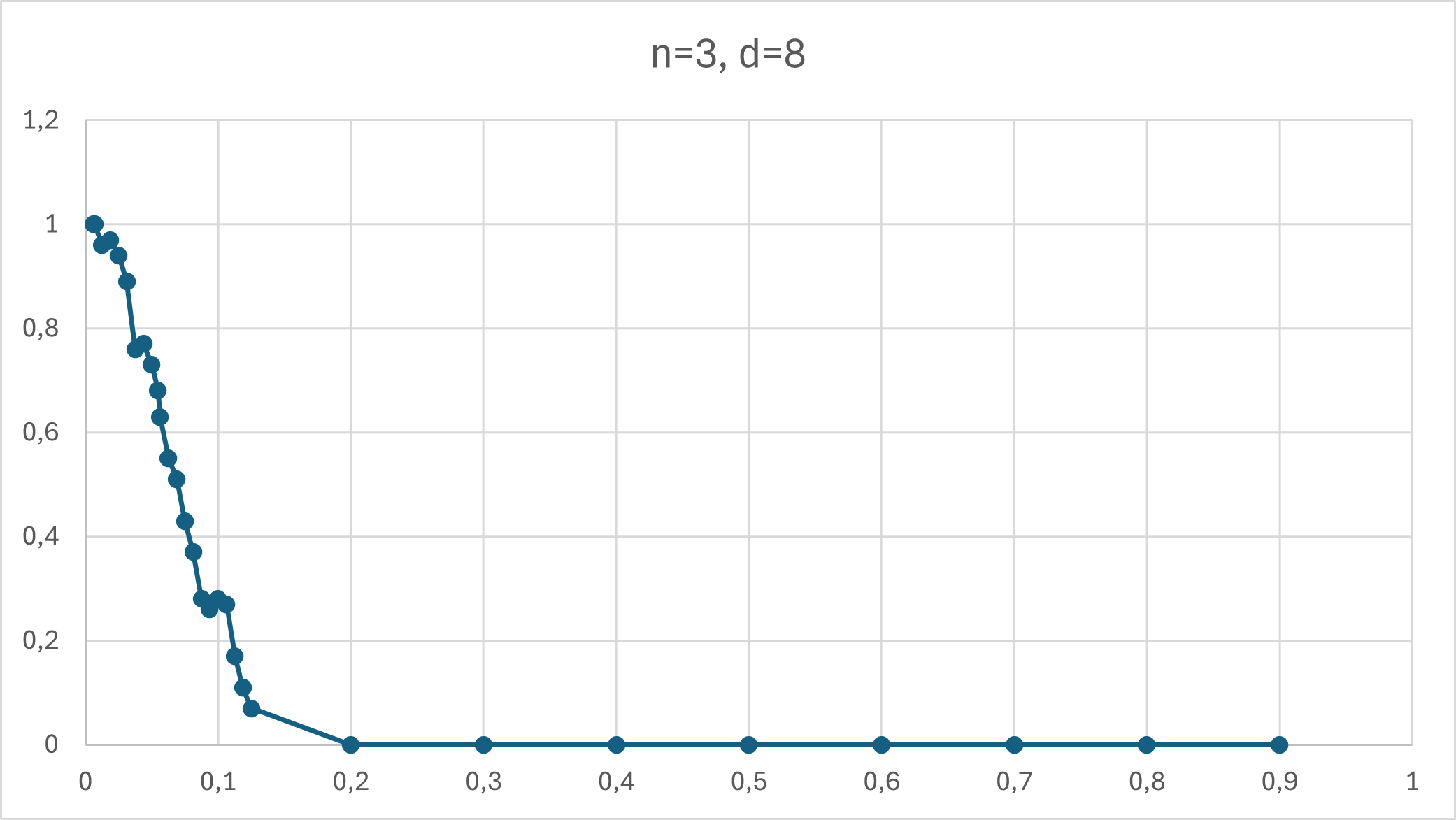}%
    \hfil%
    \includegraphics[scale=0.25]{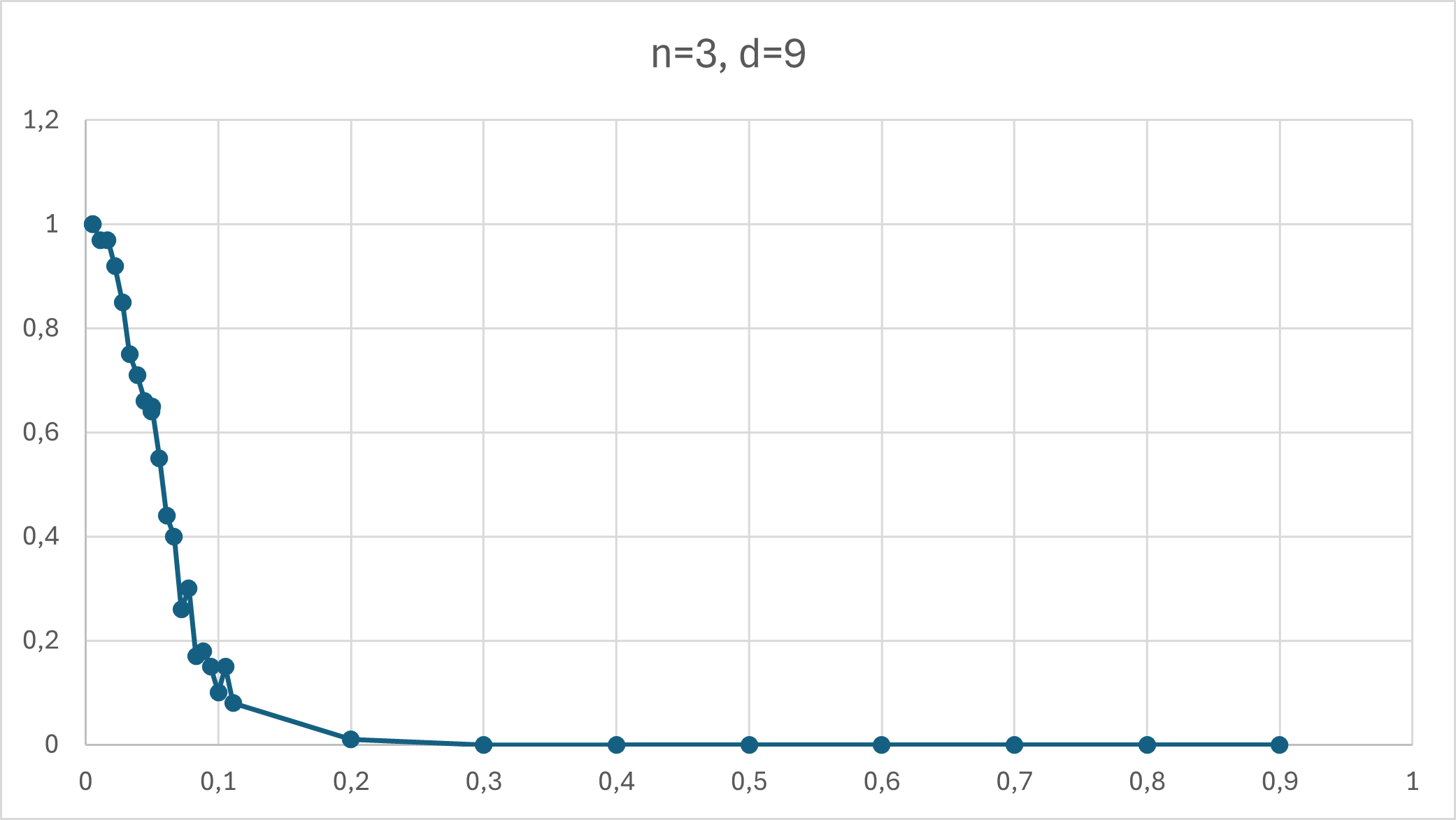}%
    \hfil%
}

\hbox to\linewidth{%
    \hfil%
    \includegraphics[scale=0.25]{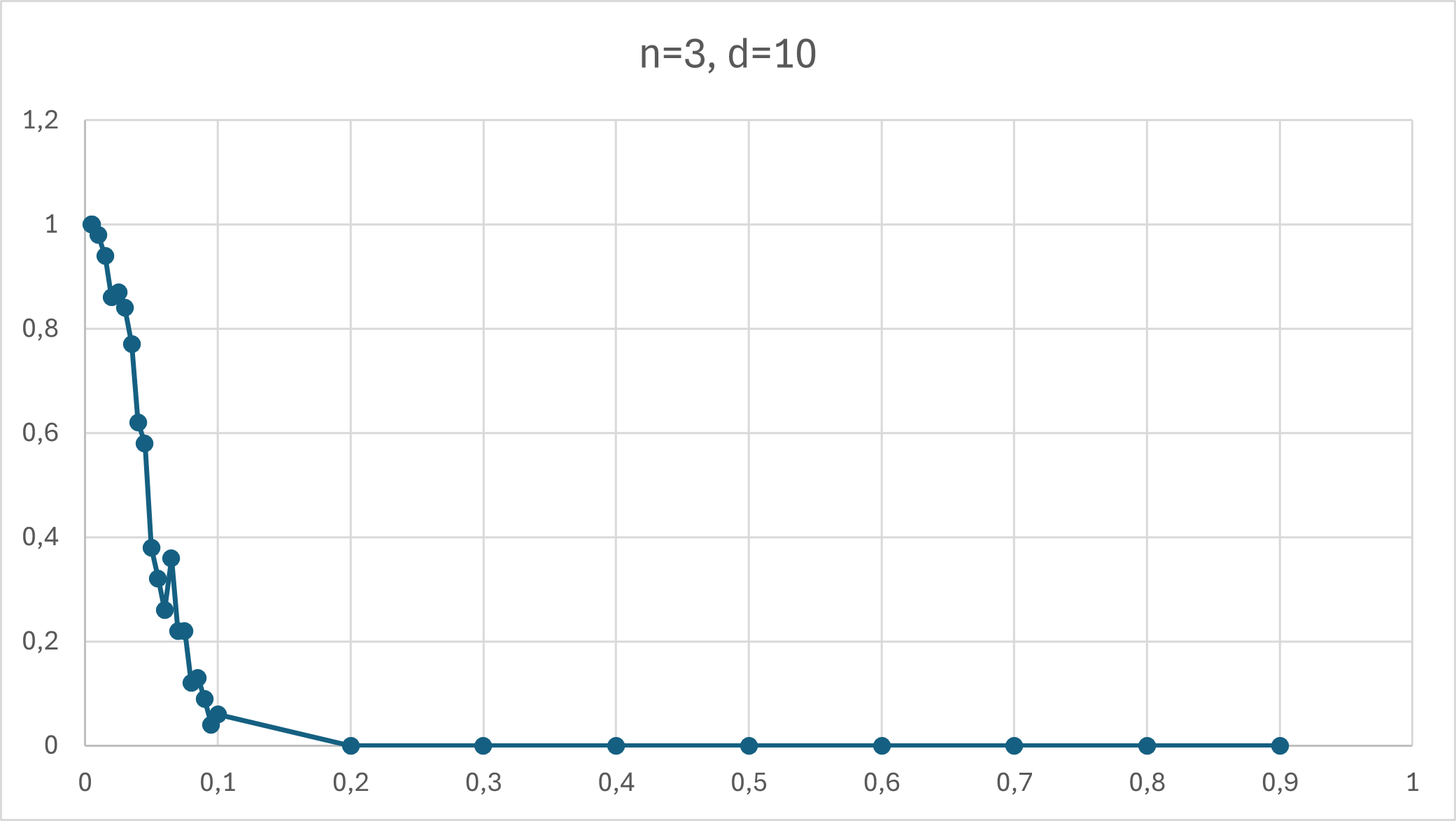}%
    \hfil%
    \includegraphics[scale=0.25]{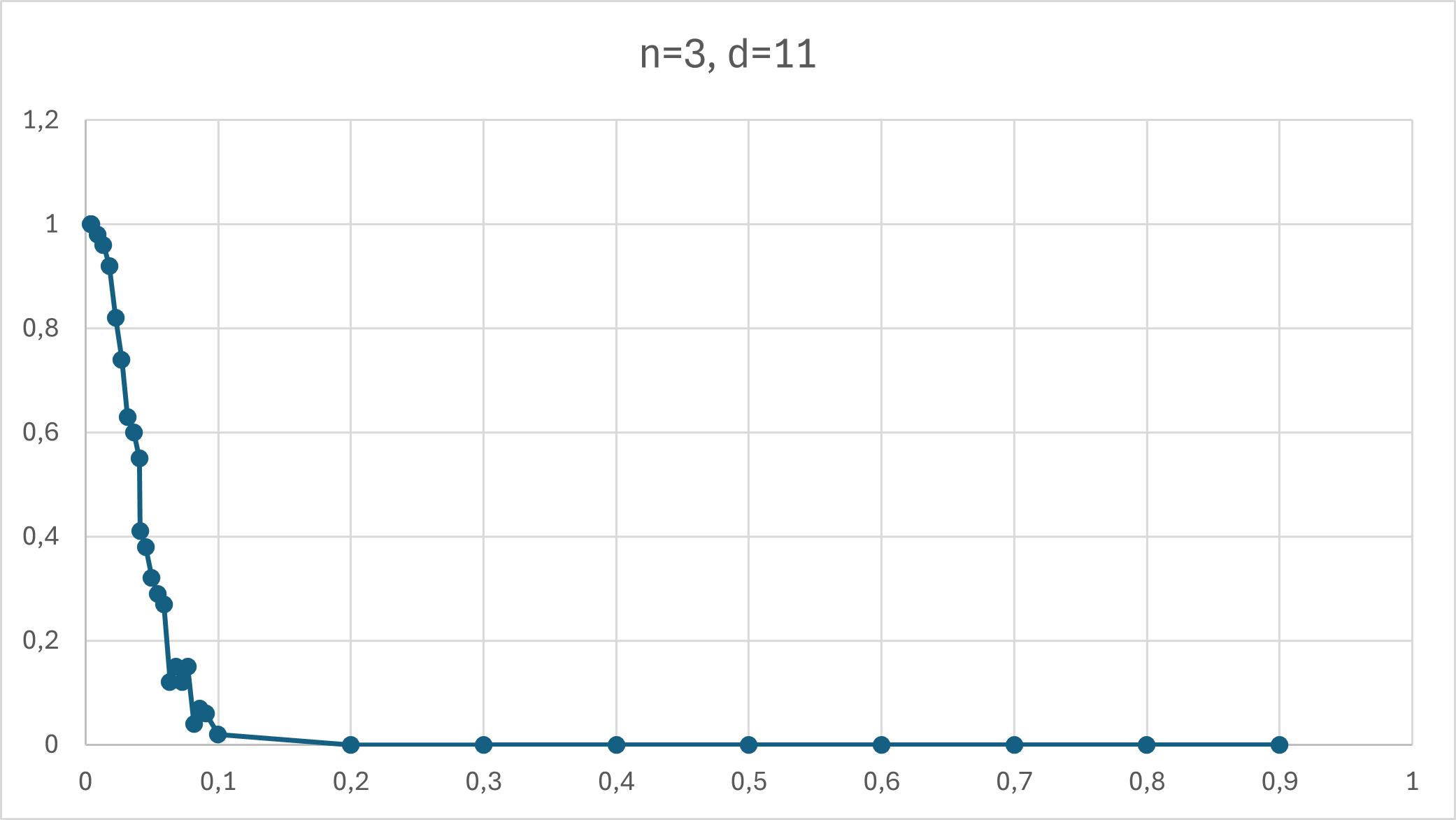}%
    \hfil%
    \includegraphics[scale=0.25]{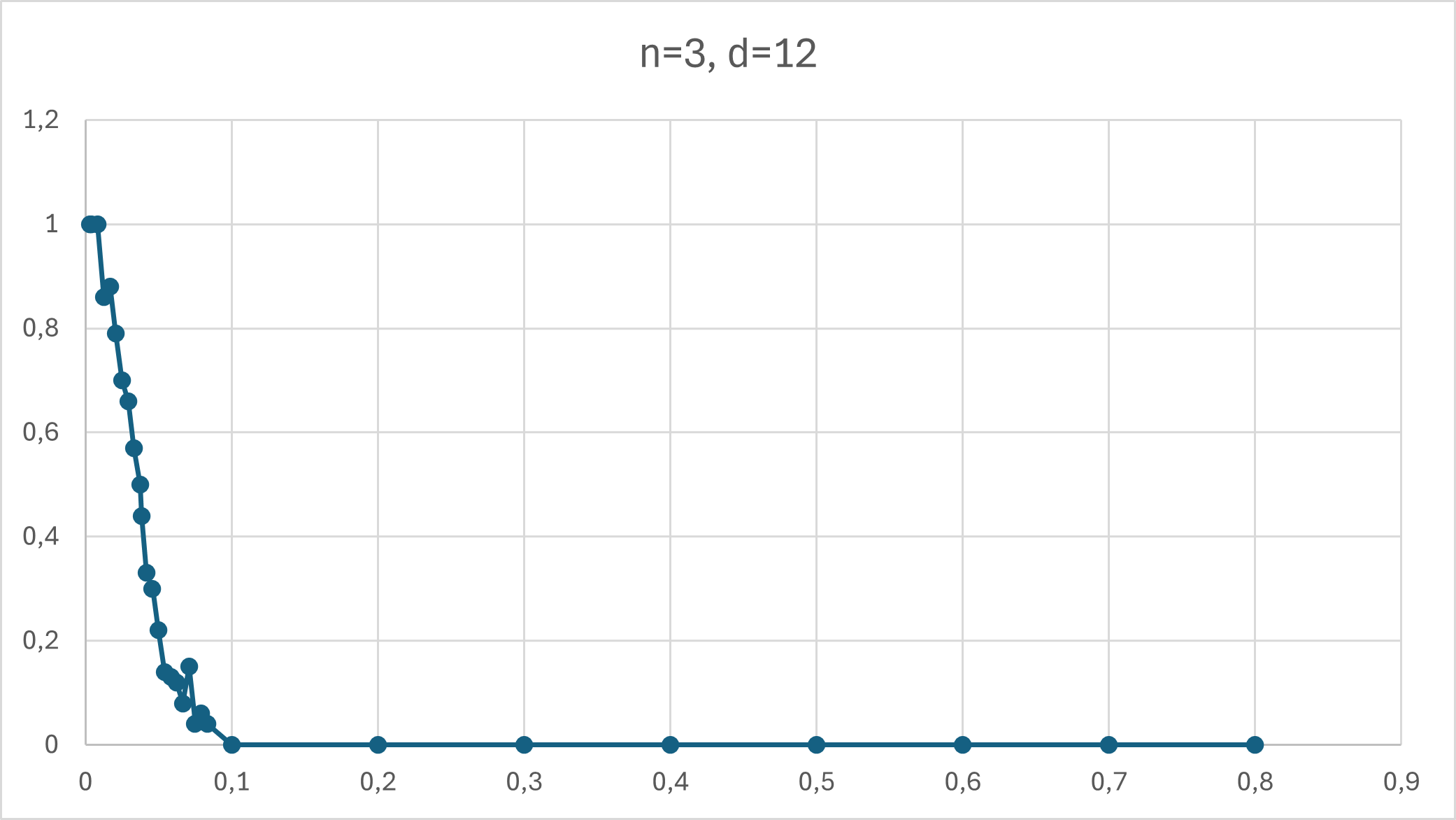}%
    \hfil%
}

  \caption{Frequency of dominant random homogeneous monomial ideals in $\mathcal{I}_{\mathrm{Gr}}(n,d,p)$ for $n=3$, $d=3,\dots,12$ and the nonzero probabilities taking values in $
\left\{ \frac{d^\ell-d^{\ell-1}}{2},\, \ell=2,\dots,n \right\} \cup \left\{ \frac19,\frac18,\dots,1 \right\} \cup \left\{\frac{1}{20d},\dots,\frac{20}{20d}  \right\}
$. %We used a sample of size 100 for each probability.
Each data point in the figure represents a sample of size $100$ for a fixed value of $p=(0,\dots,0,\alpha)$. The horizontal axes are values of $\alpha$.}
    \label{fig:frequency of dominant homogeneous RMIs}

\end{figure}

\medskip 

\begin{figure}[h!]

\hbox to\linewidth{%
    \hfil%
    \includegraphics[scale=0.25]{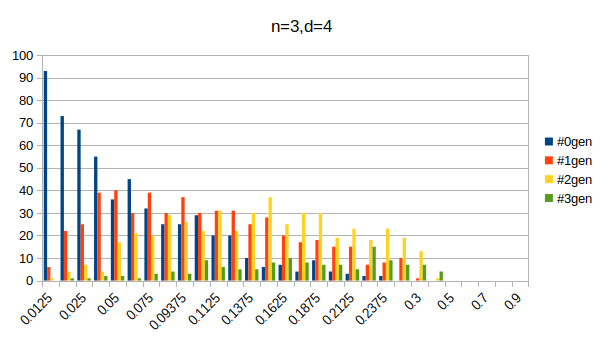}%
    \hfil%
    \includegraphics[scale=0.25]{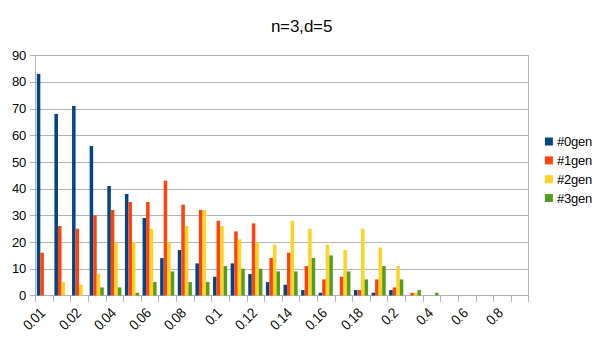}%
    \hfil%
    \includegraphics[scale=0.25]{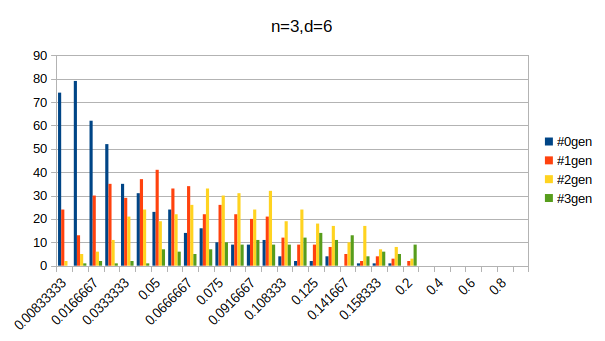}%
    \hfil%
}

\hbox to\linewidth{%
    \hfil%
    \includegraphics[scale=0.25]{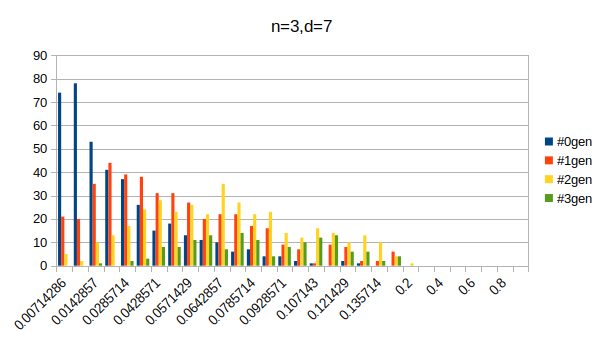}%
    \hfil%
    \includegraphics[scale=0.25]{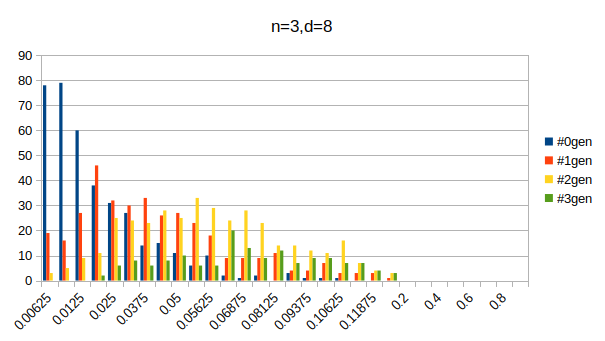}%
    \hfil%
    \includegraphics[scale=0.25]{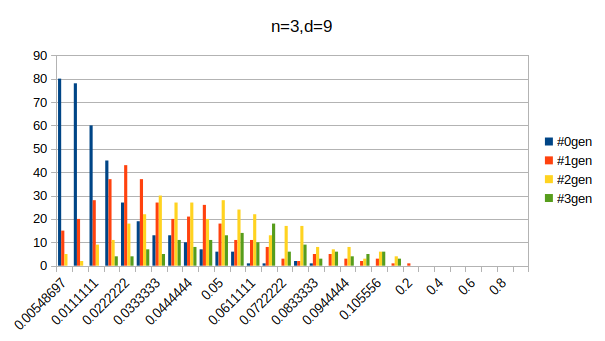}%
    \hfil%
}

\hbox to\linewidth{%
    \hfil%
    \includegraphics[scale=0.25]{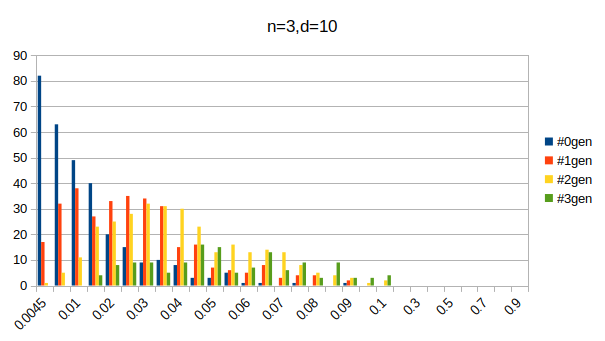}%
    \hfil%
    \includegraphics[scale=0.25]{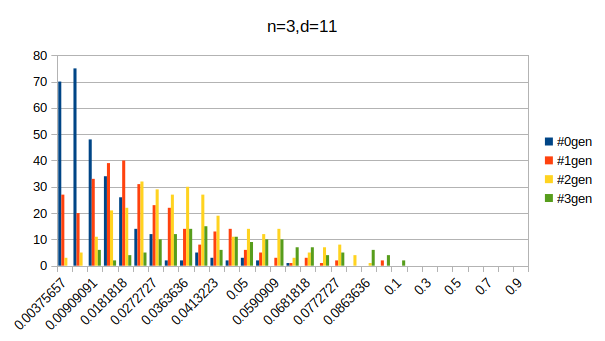}%
    \hfil%
    \includegraphics[scale=0.25]{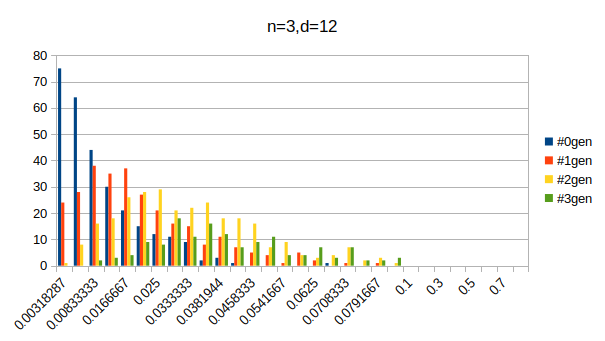}%
    \hfil%
}

    \caption{Frequencies of observing a particular number of minimal generators for the ideals in the data set displayed in \Cref{fig:frequency of dominant homogeneous RMIs}. 
    For each given probability parameter $p$ displayed on the horizontal axes, the color-coded histogram counts the number of random monomial ideals generated from $\mathcal{I}_{\mathrm{Gr}}(n,d,p)$, for $n=3$, $d=3,\dots,12$, with the given number of minimal generators.}
    \label{fig:homogeneous RMIs stratified by numgens}

\end{figure}

It is worth noting, at this point, that while we are controlling the exact degree of the generators of the homogeneous random monomial ideal, the ideals in the sample for a fixed set of parameters $n$, $D$, and $p=(0,\dots,0,\alpha)$  still have a varying \emph{number} of minimal generators. 
Namely, for every fixed value of $n$, $D$, and $\alpha$, how many monomials are chosen in $n$ variables of degree $D$? On average, this number is the total number of degree-$D$ monomials in $n$ variables, times the probability, $\alpha$,  of being selected in the random selection process. But this is only the expected or average number;   any given sample of $100$ ideals with fixed $n$, $D$, and $\alpha$ may contain ideals with less or more minimal generators than this expected number. 
To illustrate this, consider Figure~\ref{fig:homogeneous RMIs stratified by numgens}, where each sample is color-coded by the number of minimal generators. The reader will notice a bell-shape-like histogram for any given color. This corresponds to the Binomial distribution on the  probability that an \er\ random monomial ideal $I$  has a given number of generators.

\subsection{Random dominant ideals with a fixed number of generators}

The previous disucssion leads us to consider an alternative \er\ type model for random monomial ideals; namely, one in which instead of the probability parameter $p$ we use an integer $M$ for the required number of generators; see the method {\tt randomMonomialIdeals} in the {\tt Macaulay2} package \cite{RMIm2}. This model version is denoted by $\mathcal B(n,D,M)$, where $M$ is the fixed number of generators in the random ideal. We are interested again in the graded case, so that $M=(0,\dots,0,g)$ is a vector that specifies the number of generators desired in each degree up to $D$, and we choose only one nonzero entry, $g$.  

\begin{figure}
    \begin{subfigure}[t]{1\textwidth}
        \centering
        \includegraphics[scale=0.25]{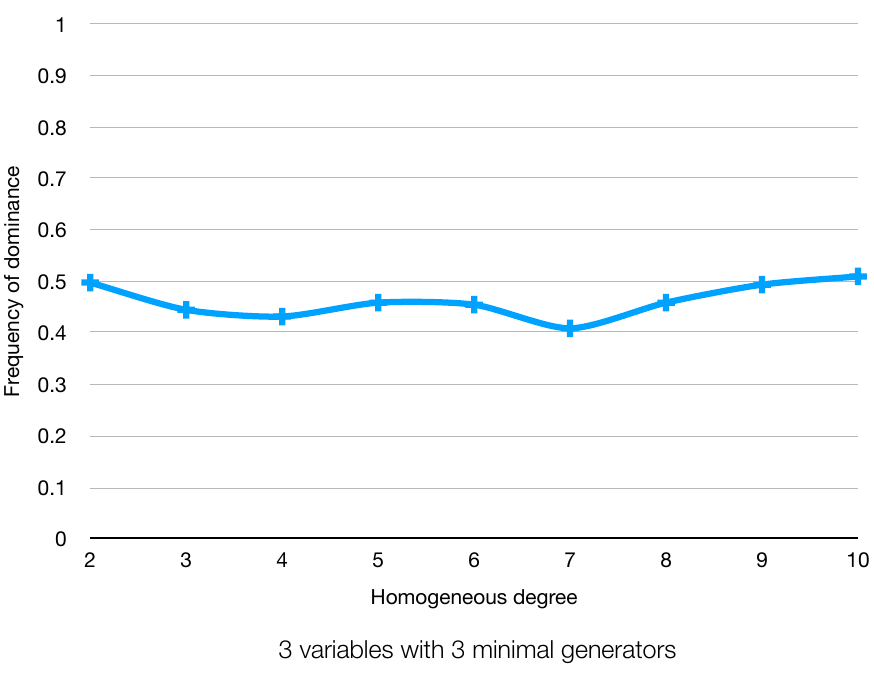}
        \includegraphics[scale=0.25]{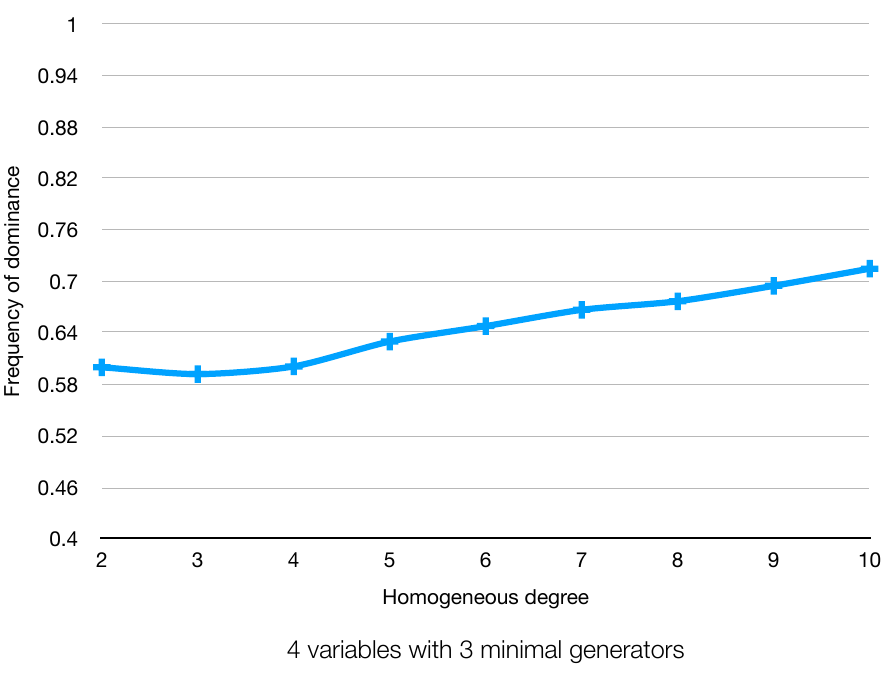}
        \includegraphics[scale=0.25]{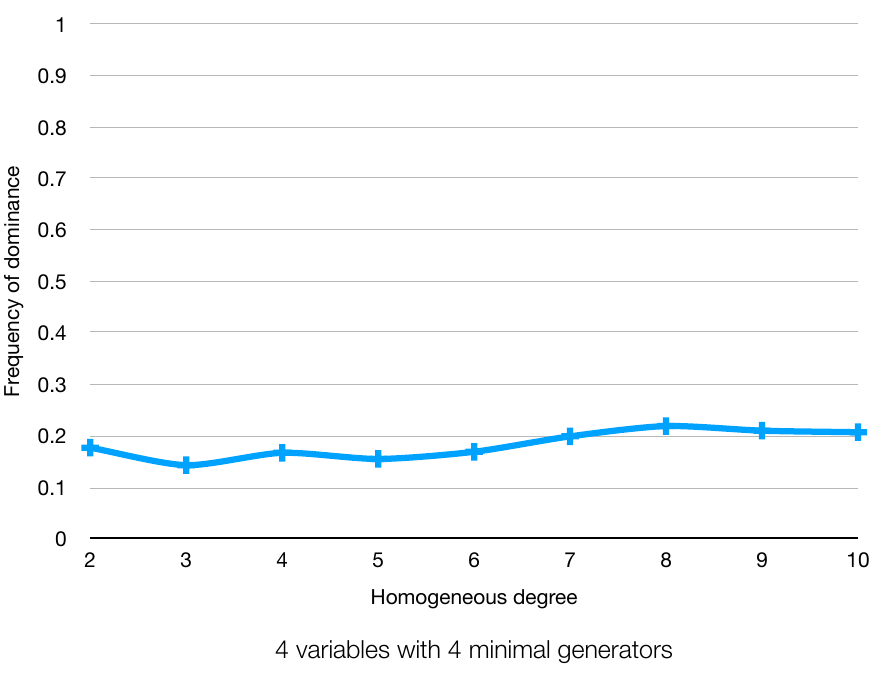}
        \caption{3 and 4 variables}
    \end{subfigure}
    \begin{subfigure}[t]{1\textwidth}
            \centering
        \includegraphics[scale=0.25]{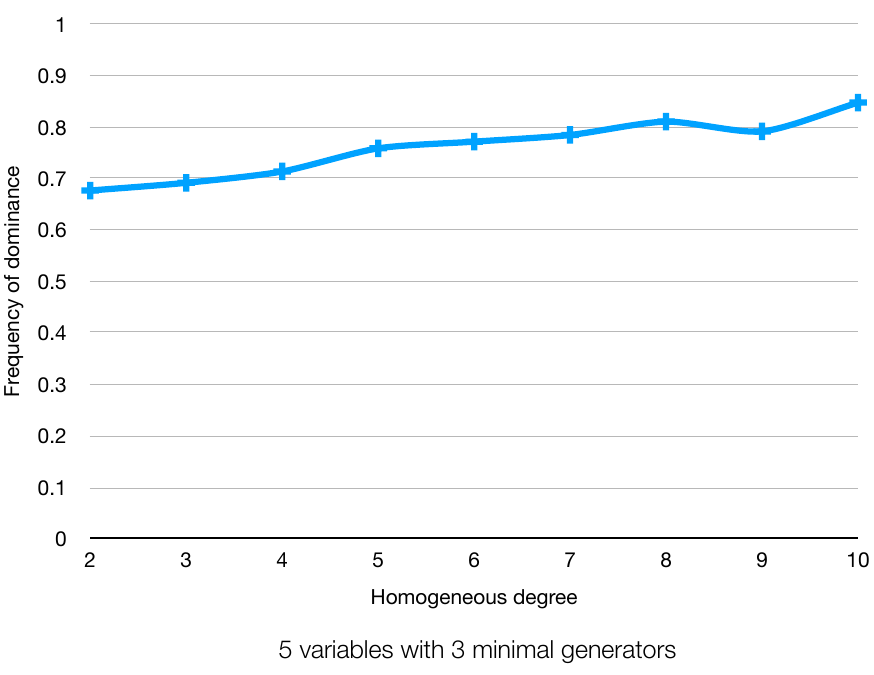}
        \includegraphics[scale=0.25]{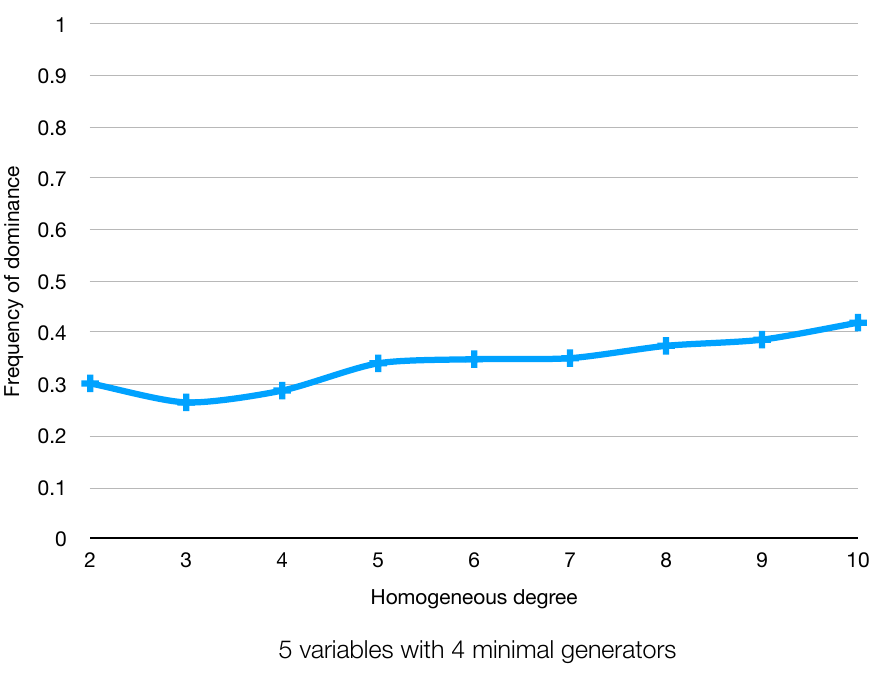}
        \includegraphics[scale=0.25]{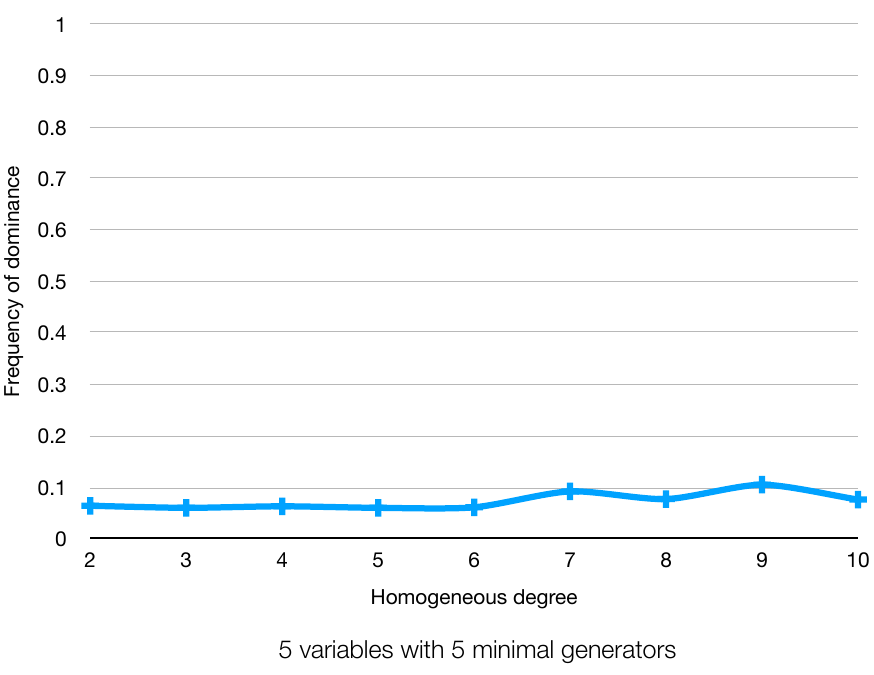}
         \caption{5 variables}
    \end{subfigure}
    \begin{subfigure}[t]{1\textwidth}
        \centering
        \includegraphics[scale=0.25]{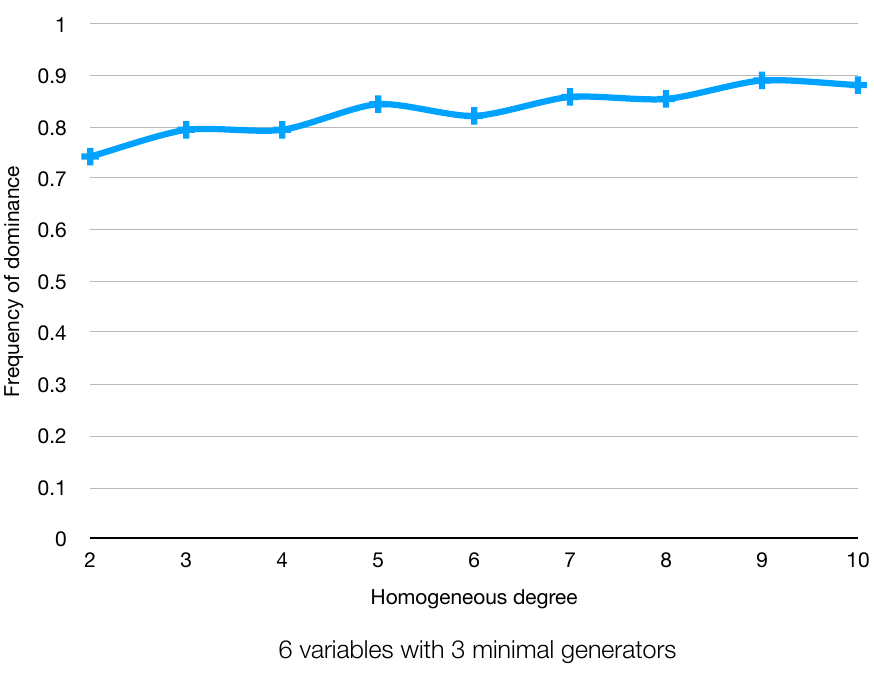}
        \includegraphics[scale=0.25]{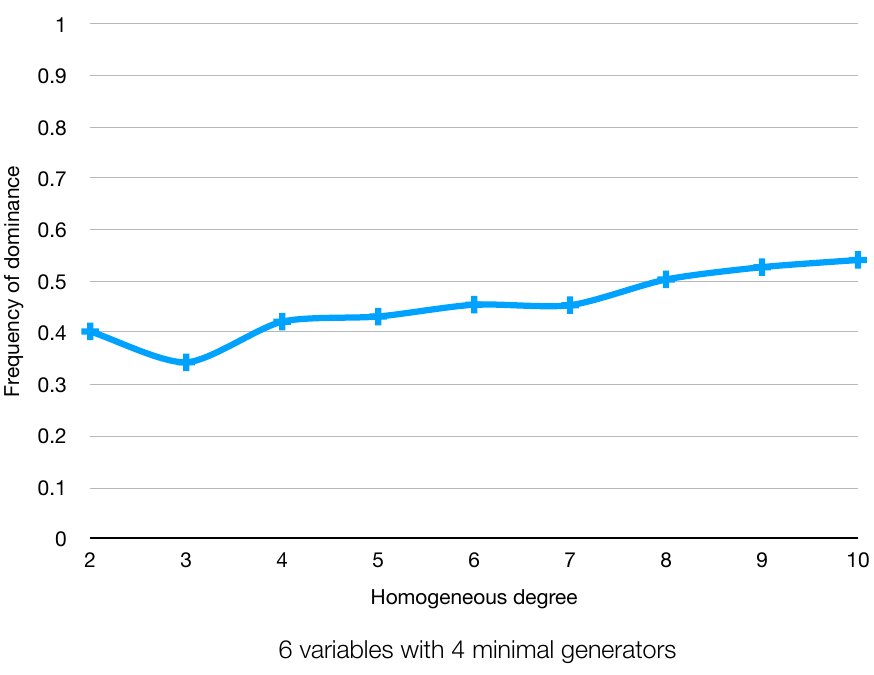}
        \includegraphics[scale=0.25]{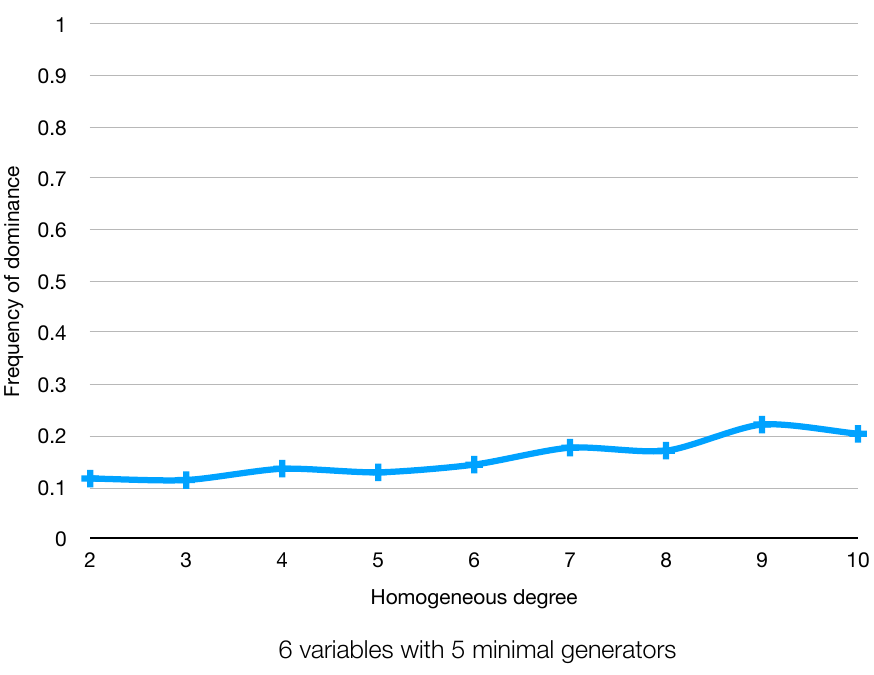}
        \includegraphics[scale=0.25]{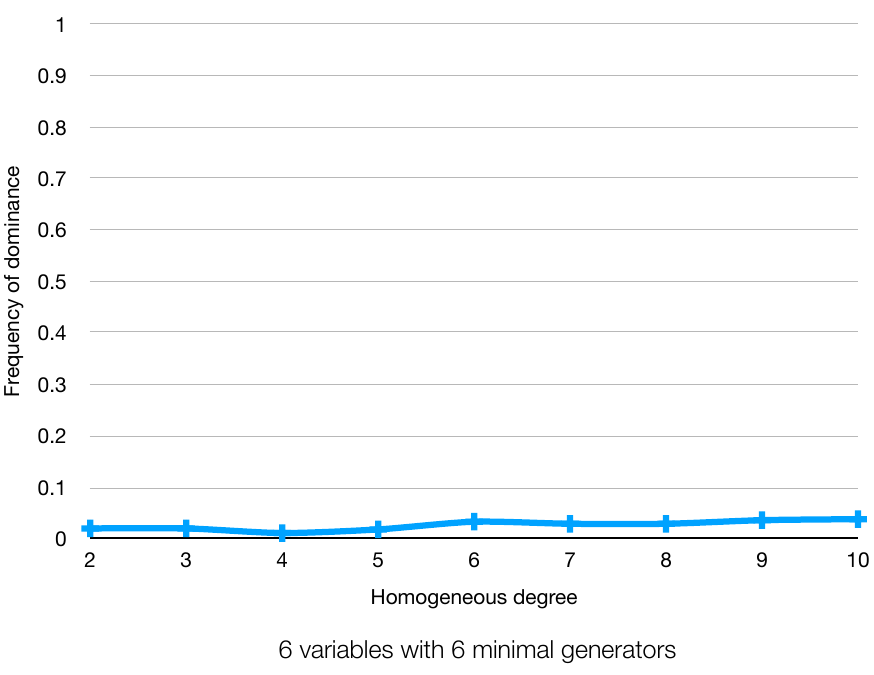}
         \caption{6 variables}
    \end{subfigure}

    \begin{subfigure}[t]{1\textwidth}
        \centering
        \includegraphics[scale=0.25]{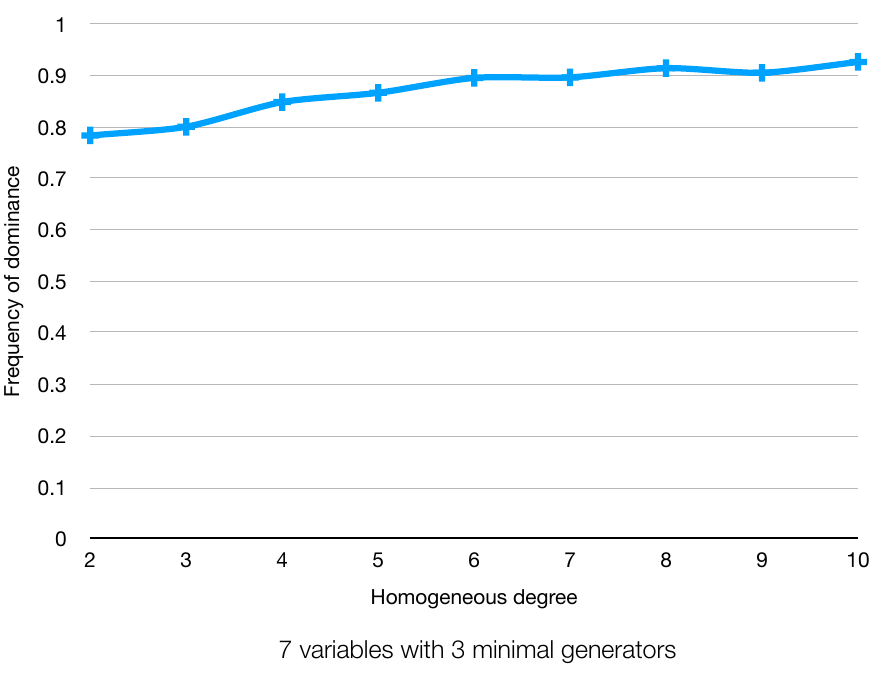}
        \includegraphics[scale=0.25]{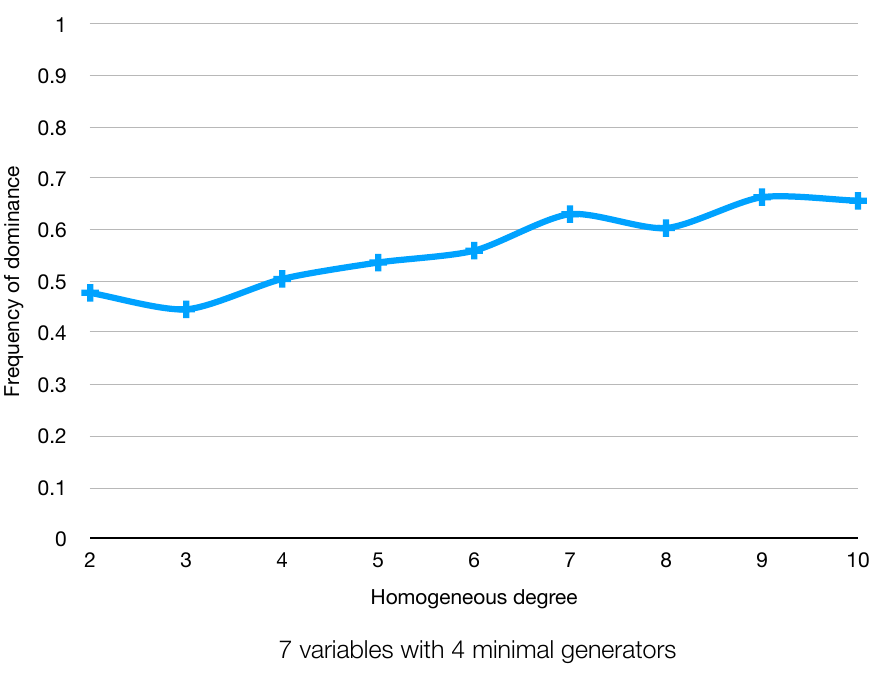}
        \includegraphics[scale=0.25]{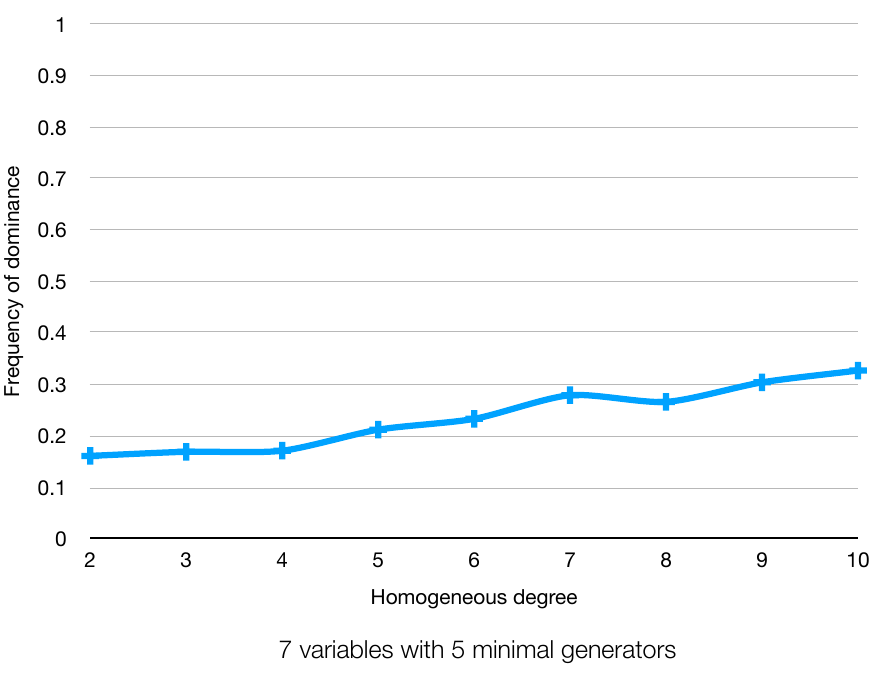}

        \includegraphics[scale=0.25]{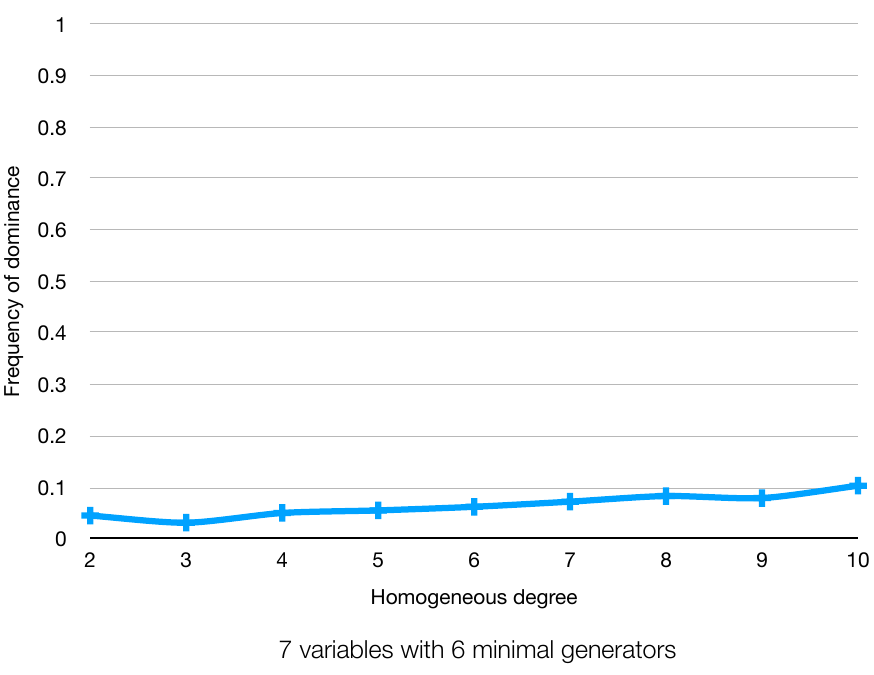}
        \includegraphics[scale=0.25]{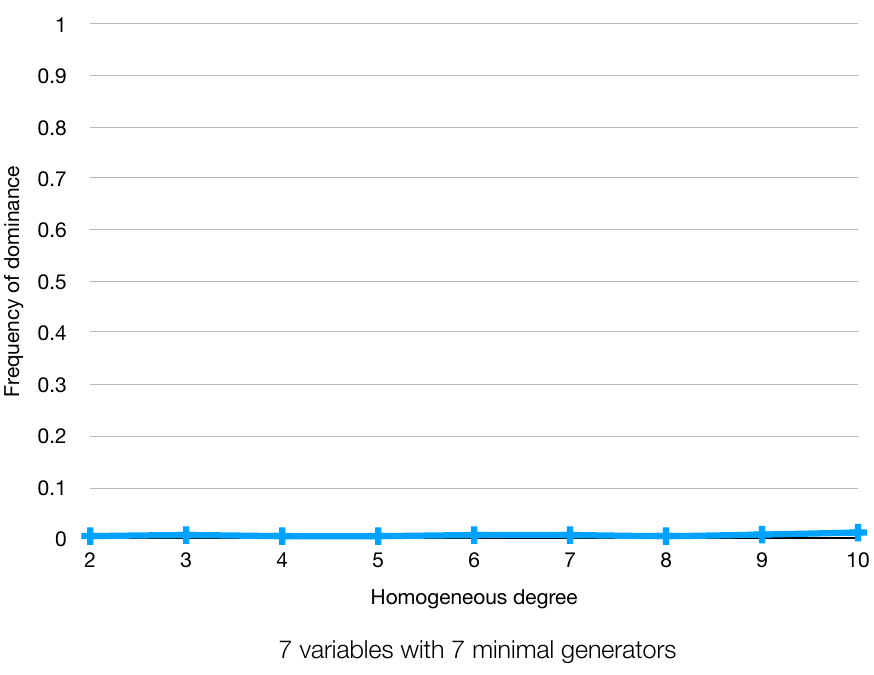}
        \caption{7 variables}
    \end{subfigure}

    \begin{subfigure}[t]{1\textwidth}
        \centering
        \includegraphics[scale=0.25]{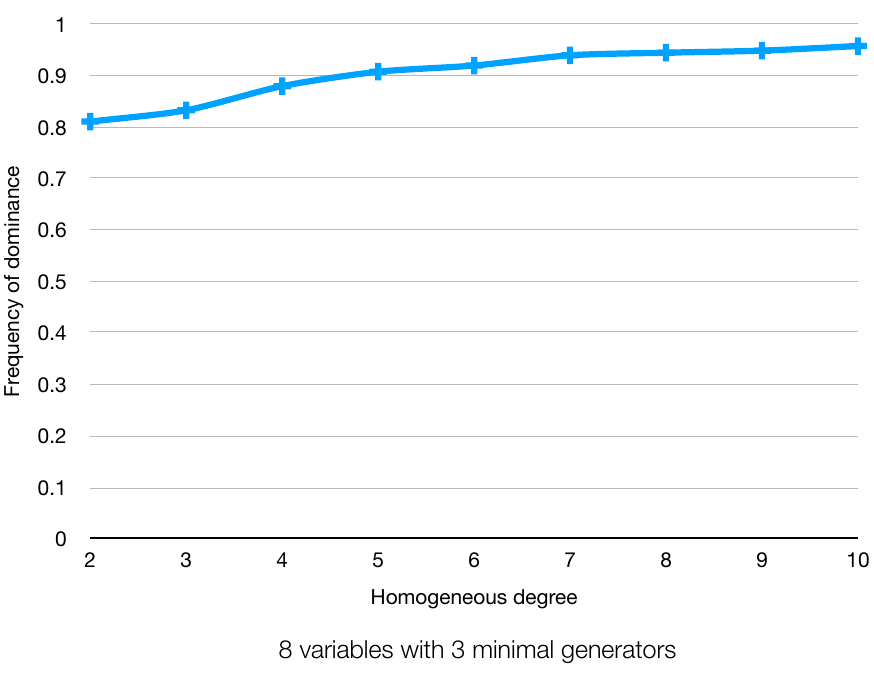}
        \includegraphics[scale=0.25]{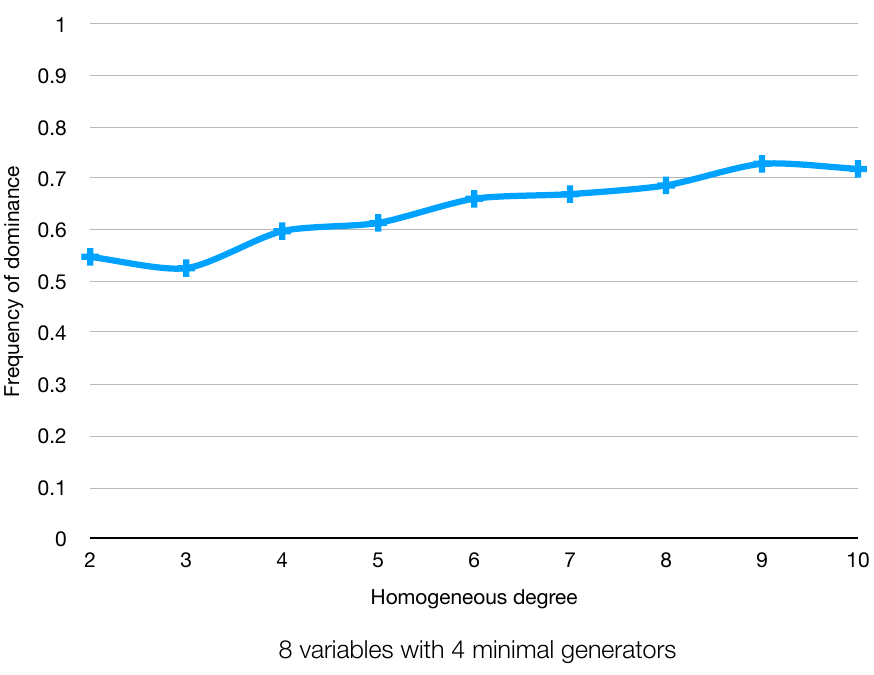}
        \includegraphics[scale=0.25]{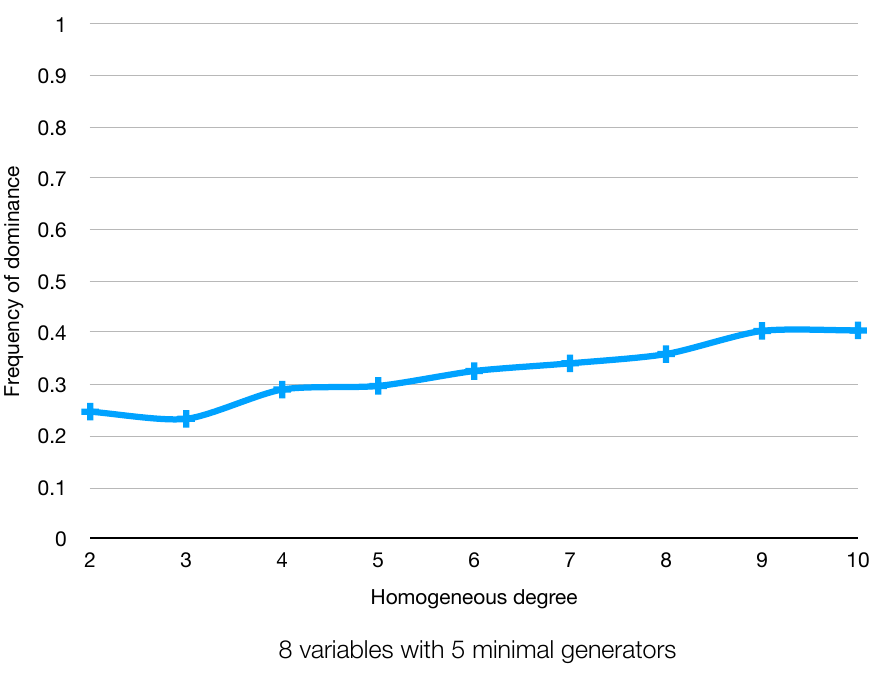}
        \includegraphics[scale=0.25]{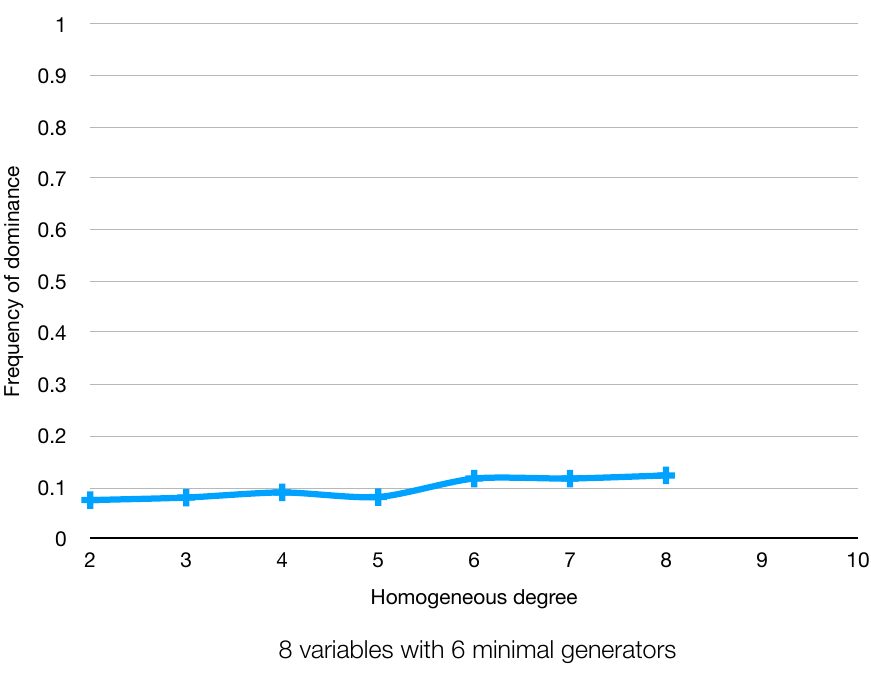}
        \caption{8 variables, from left: 3, 4, 5, and 6 minimal generators. The computation for 7 or 8 generators in the 8-variable case ran out of memory.}
\end{subfigure}    
\caption{Proportion of dominant ideals in a sample of random monomial ideals. Each data point on each graph represents one sample of size $1000$, for a fixed number of variables $n$, fixed generator degree $D$, and fixed number of generators.}
\label{fig: proportion dominant for fixed number of min gens by n}
\end{figure}

With this, we can plot the empirical conditional distribution of the proportion of dominant ideals, given the number of minimal generators of the random homogeneous monomial ideal, with fixed $n$ and $D$. Appendix~\ref{sec:code for RMI homogeneous model with fixed number of generators} contains the code that generates data supporting Figure~\ref{fig: proportion dominant for fixed number of min gens by n}.

 \medskip 
 The reader may now notice that the trends reported in this figure are `expected'; namely, when there are few generators, there is a higher proportion of dominant ideals, then when there are many generators (of course upper bounded by the number of variables). 
In turn, this observation leads to the following pleasing conclusion: 
\emph{
The \er\ type family of models for random monomial ideals is a suitable tool for studying the dominance property.  
}

Therefore, when faced with the challenges of computing examples with many variables and high degrees of generators, rather than attempting an exhaustive enumeration, commutative algebraists can feel comfortable generating small samples of random monomial ideal examples using this model family, and studying the dominance property.

\bibliographystyle{alpha}
 
\bibliography{References.bib}

\appendix

\section{Macualay2 Code for generating samples of random monomial ideals from the basic ER-type model}\label{sec:code for RMI basic model}

In this  set of experiments, consider the following range for the number of variables $n=3,\dots,6$. For each fixed $n$, we  vary the number of maximal degree $d$ of the randomly chosen monomials.  \Cref{code: ex basic model n3} contains an instance of the code we used to perform our experiments in the specific case of $n=3$ variables which work as follows.

For each $d=2,\dots,15$, The set of probabilities we chose for our computations \texttt{PossibleP}
is created by adding the values
$$ 
\frac{d^{-i}-d^{-(i-1)}}{2}, 
 \hbox{ for } i=2,\dots,n 
$$
to the starting set of probabilities $\left\{ \frac{1}i, i=1,\dots,9 \right\}$, removing duplicates if there are any.

Once the set of possible probabilities is created, for each probability \texttt{p} in \texttt{PossibleP} we create a list \texttt{L} of 50 random monomial ideals in $\mathcal{I}(n,d,p)$. The list \texttt{listDom} stores the dominant ones computed via the function \texttt{isDominant}, see \cref{M2dominantF}. The list \texttt{numG} is a tuple containing the number of dominant ideals with $i$ minimal generators for $i=1,\dots,n$ since we also want to count the number of minimal generators of every monomial ideals. Finally, we store in a text file, called \texttt{dominantIdealData$\_3$.txt} in this case, the tuple \texttt{(d,p, length listDom, numG)}.

\begin{lstlisting}[language=code, caption={Example of Macaulay2 code to generate 50 random ideals $I\in \mathcal{B}(n,d,p)$ with $n=3$, $d=2,\dots,15$ and store the number of dominant ones according to the number of minimal generators in a text file, one sample per line.}, label=code: ex basic model n3]
n = 3 
Dd = 15 
D = toList (2..Dd)

for d in D do(
    possibleP = sort toList apply(apply(1..n,x->1/d^x),s->sub(s,RR));
    midP = {};
    for i from 1 to length possibleP-1 do(
    	midP = append(midP,(possibleP_i-possibleP_(i-1))/2);
    	);
    possibleP = unique sort toList join(midP,toList apply(1..9,x->sub(x/10,RR)));
    for p in possibleP do(
        L = randomMonomialIdeals(n,d,p,50);
	listDom = select(L,x->isDominant(x) == true);
	numGens = apply(listDom,x-> #flatten entries mingens(x));
	numG = toList apply(0..n,y->length select(numGens,x->x==y));
        f = openOutAppend "dominantIdealData_3.txt";
        f << (d,p,length listDom,numG) << endl << close;    
    )
)
\end{lstlisting}

\section{Macaulay2 code for generating samples of homogeneous ER-type random monomial ideals}
\label{sec:code for RMI homogeneous model}

In these experiments, the number of variables ranges  $n=3,\dots,6$. For each $n$, we consider varying degree $d=2,\dots,12$. \Cref{code: ex graded model} contains the \texttt{M2} code we used to perform our simulations in the particular case $n=3$. It works as follows. For every degree $d=2,\dots,d$ the list of probabilities considered \texttt{homProb} contains is a $d$-tuple $(0,\dots, 0,\alpha)$ where 
$$
\alpha \in \left\{ \frac{d^\ell-d^{\ell-1}}{2},\, \ell=2,\dots,n \right\} \cup \left\{ \frac{1}{9},\frac{1}{8},\dots,1 \right\} \cup \left\{\frac{1}{20d},\dots,\frac{20}{20d}  \right\}.
$$
For each probability $p$ in \texttt{homProb}, the list \texttt{L} contains 100 random monomial ideals chosen within the distribution $\mathcal{I}_{\mathrm{Gr}}(n,d,p)$ of the graded model. The dominant ideals of \texttt{L} are selected through the function \texttt{isDominant} and stored in the list \texttt{listDom}. The list \texttt{numG} contains the number of dominant ideals sorted by the minimal number of generators that a dominant ideal can possibly have, namely $1,\dots,n$. Lastly, we store  all relevant information given by the tuples $(d,\alpha, \# \texttt{listDom},\texttt{numG})$ in the file \texttt{dominantHomIdealData$\_3$n.txt}, again one sample per line.

\begin{lstlisting}[language=code, caption={Example of Macaulay2 code to generate 100 random homogeneous ideals with $n=3$, $d=2,\dots,12$ and the only nonzero probabilities taken in $
\left\{ \frac{d^\ell-d^{\ell-1}}{2},\, \ell=2,\dots,n \right\} \cup \left\{ \frac{1}{9},\frac{1}{8},\dots,1 \right\} \cup \left\{\frac{1}{20d},\dots,\frac{20}{20d}  \right\}.
$
The text file created stores the number of dominant ideals generated in each sample, along with the sample parameters, one sample per line.}, label=code: ex graded model]
n=3
Dd=12 
for d in D do(
possibleP = sort toList apply(apply(1..n,x->1/d^x),s->sub(s,RR));
    midP = {};
    for i from 1 to length possibleP-1 do(
    	midP = append(midP,(possibleP_i-possibleP_(i-1))/2);
    	);
    possP = unique sort toList join(midP,toList apply(1..9,x->sub(x/10,RR)));
    possibleP2=toList apply(1..20,x->sub(x/(20*d),RR));
    possibleP= unique sort join(possP,possibleP2);
    zeroList= toList (d-1:0.0);
    homProb= apply(possibleP,x->append(zeroList,x));
    for p in homProb do(
        L = randomMonomialIdeals(n,d,p,100);
	listDom = select(L,x->isDominant(x) == true);
	numGens = apply(listDom,x-> #flatten entries mingens(x));
	numG = toList apply(0..n,y->length select(numGens,x->x==y));
        f = openOutAppend "dominantHomIdealData_3n.txt";
        f << (d,p#(length p-1),length listDom,,numG) << endl << close;    
    )
)
\end{lstlisting}

\section{Macaulay2 code for generating samples of homogeneous ER-type random monomial ideals with fixed number of generators}
\label{sec:code for RMI homogeneous model with fixed number of generators}

This section contains the code necessary to produce statistics on the frequency of dominance for a random homogeneous monomial ideal with a fixed number of generators.

More precisely, if we fix integers $n\geq 2$, $g\leq n$, $d\geq 2$, then we want to understand how frequent is for a monomial ideal in $n$ variables, with $g $ minimal generators and of homogenous degree $d$ to be dominant.

To perform our computations we use the following function from the Macaulay2 package % \emph{Random Monomial ideals}  
 \cite{RMIm2}:  \texttt{randomMonomialSets(n,d,deg,N, Strategy=>"Minimal")}. 
  The parameter \texttt{n} represents the number of variables, \texttt{d} is the maximum degree considered, \texttt{deg} is a $d$-tuple of the form $(0,\dots,0,g)$ where $g$ is the number of minimal generators we want. Lastly, the parameter \texttt{N} is for how many random monomial sets we want to select. In this case we are not directly controlling the probability on the generated sets, which in this case is chosen uniformly at random from the collection of all monomials of degree at most $d$.

\Cref{code: hom fixed mingens} contains the \texttt{Macaulay2} code used to perform our computations that needs to be used after loading the package \texttt{RandomMonomialIdeals} and the and the function \texttt{isDominant}.   Unfortunately, the computation stops for $n=8,g=6,d=8$ because the machine used to perform the calculation runs out of memory.
Figure~\ref{fig: proportion dominant for fixed number of min gens by n} shows the results of this computations.
\begin{lstlisting}[language=code, caption={ Macaulay2 code to generate 1000 random homogeneous ideals with $3\leq n\leq 10$ variables, $2\leq g\leq n$ minimal generators of homogenous degree $2\leq  d\leq 10$. The text file created stores the number of dominant ideals generated in each sample, along with the sample parameters, one sample per line.}, label=code: hom fixed mingens]
N = 1000
for n in 3..10 do(
    for NumGen in 3..n do(
    	for d in 2..10 do(
    	    zeroList = toList (d-1:0);
    	    deg = append(zeroList,NumGen);
	    L = randomMonomialSets(n,d,deg,N, Strategy=>"Minimal");
    	    listDom = select(L,x->isDominant(ideal x) == true);
    	    f = openOutAppend "DominantsWithFixedGens1000.txt";
    	    f << (n,NumGen,d,length listDom) << endl << close;
    		);
	);
);
\end{lstlisting}

%\noindent{\bf Authors' addresses}
\smallskip
\small

\noindent Anna Maria Bigatti: Dipartimento di Matematica, Università degli Studi di Genova, Via Dodecaneso 35, Genova, 16146, Italy.\\  {\tt a.m.bigatti@unige.it}
\medskip

\noindent Nursel Erey: Gebze Technical University, Department of Mathematics, Gebze, 41400 Kocaeli, Turkey.\\
 {\tt nurselerey@gtu.edu.tr}
\medskip

\noindent Selvi Kara: Department of Mathematics, Bryn Mawr College, Bryn Mawr, PA, USA, 19010.\\       { \tt  skara@brynmawr.edu}
\medskip

\noindent Augustine O'Keefe: Connecticut College, Department of Mathematics and Statistics, New London, CT, USA 06320.
\\ {\tt aokeefe@conncoll.edu}
\medskip

\noindent Sonja Petrovi\'c: Illinois Institute of Technology, Department of Applied Mathematics, Chicago, IL, USA, 60616.
\\ {\tt sonja.petrovic@iit.edu}
\medskip

\noindent Pierpaola Santarsiero:
Universit\`a di Bologna, Dipartimento di Matematica, Piazza di Porta S. Donato 5, 40126 Bologna, Italy.
\\ {\tt pierpaola.santarsiero@unibo.it}
\medskip

\noindent Janet Striuli: Department of Mathematics, Fairfield University, Fairfield, CT, USA, 06824.
\\ {\tt jstriuli@fairfield.edu}

\end{document}